\documentclass[aap,dvips]{arximspdf}
\usepackage{graphicx}

\doi{10.1214/105051607000000285}
\volume{18}
\issue{2}
\pubyear{2008}
\firstpage{379}
\lastpage{426}

\makeatletter

\def\eqref#1{(\ref{#1})}
\newtheorem{theorem}{Theorem}
\newtheorem{prop}{Proposition}
\newtheorem{Corollaire}{Corollary}
\newtheorem{lemme}{Lemma}
\newproclaim{Remarque}{Remark}
\newproclaim{Exemple}{Example}
\newproclaim{exm}{Example}
\newproclaim{ass}{Assumption}
\makeatother

\begin{document}
\begin{frontmatter}

\title{Recursive computation of the invariant measure of a stochastic
differential equation driven by a~L\'evy~process}
\runtitle{Recursive computation}
\pdftitle{Recursive computation of the invariant measure of a stochastic
differential equation driven by a Levy process}

\begin{aug}
\author[A]{\fnms{Fabien} \snm{Panloup}\corref{}\ead[label=e1]{fpanloup@insa-toulouse.fr}}
\runauthor{F. Panloup}
\affiliation{Universit\'e Paris 6}
\address[A]{Laboratoire de Probabilit\'es\\
\quad  et Mod\`eles Al\'eatoires\\
Universit\'e Paris 6\\
UMR 7599\\
bureau 4D1\\
175 rue du Chevaleret\\
F-75013 Paris\\
France\\
\printead{e1}} 
\end{aug}

\received{\smonth{10} \syear{2005}}
\revised{\smonth{4} \syear{2007}}

%
\begin{abstract}
We study some recursive procedures based on exact
or approximate
Euler schemes with decreasing step to compute the invariant measure of
L\'evy driven SDEs.
We prove
the convergence of these procedures toward the invariant
measure
under weak conditions on the moment of the L\'evy process and on the
mean-reverting of the dynamical system.
We also show that an a.s. CLT for stable processes can be
derived from our main results. Finally, we illustrate our results
by several simulations.
\end{abstract}

%
\pdfkeywords{60H35, 60H10, 60J75, 60F05, Stochastic differential equation,
Levy process, invariant distribution, Euler scheme, almost sure central limit theorem,}
\end{frontmatter}

\section{Introduction}\label{s1}

\subsection{Objectives and motivations}

This paper is devoted to the computation of the
invariant measure (denoted by $\nu$) of ergodic stochastic
processes which obey a stochastic differential equation (SDE)
driven by a L\'evy process. Practically, we want to construct a
sequence of empirical measures $(\bar{\nu}_n(\omega,dx))_{n\ge1}$
which can be recursively simulated and such that
$\bar{\nu}_n(\omega,f)\rightarrow\nu(f)$ a.s. for a range of
functions~$f$ containing bounded
continuous functions.

In the case of Brownian diffusions, some methods have already been
developed by several authors to approximate the invariant measure
(see Section~\ref{back}), but this paper seems to be the first
one that deals with this problem in the case of general L\'evy
driven SDEs. The motivation for this generalization is the study
of dynamical systems that are widely used in modeling. Indeed,
there are many situations where the noise of the dynamical system
is discontinuous or too intensive to be modeled by a Brownian
motion. Let us consider
an example that comes from the fragmentation-coalescence theory. In
situations such as
polymerization phenomenons, when temperature is near to its
critical value, molecules constantly break-up and recombine.
This situation has been modeled by Berestycki \cite{bib17}
through what he terms EFC (Exchangeable Fragmentation-Coalescence) process.
The mass of the dust generated by this process (see \cite{bib17}
for more details) is a solution to a mean-reverting SDE for which
the noise component is driven by a subordinator (an increasing
L\'evy process).
We come back to this example in Section \ref{simulations}.

For other examples of situations where models that use a L\'evy
driven SDE are adapted, we refer to Barndorff-Nielsen et
al. \cite{bib18} for examples in financial modeling (where
ergodic L\'evy driven SDEs are usually used to model the
volatility of a financial market), Protter--Talay \cite{bib19}
for examples in finance, electrical engineering$,\ldots$\ or Deng
\cite{bib20}, who models the spot prices of electricity by a
mean-reverting Brownian diffusion perturbed by a compound Poisson
noise.

\subsection{The stochastic differential equation}

According to the
L\'evy--Khintchine decomposition (for this
result and for introduction to L\'evy processes, see, e.g., Bertoin
\cite{bertoin}, Protter \cite{bib1} or Sato \cite{bib11}),
an $\mathbb{R}^l$-valued L\'evy process $(L_t)$ with L\'evy measure $\pi$
admits the following decomposition: $L_t=\alpha t+\sqrt{Q}
W_t+Y_t+N_t$, where $\alpha\in\mathbb{R}^l$, $Q$ is a symmetric positive
$l\times l$ real matrix, $(W_t)$ is a $l$-dimensional standard
Brownian motion, $(Y_t)$ is a centered $l$-dimensional L\'evy
process with jumps bounded by 1 and characteristic function
given for every $t\ge0$ by
\[
\mathbb{E}\bigl\{e^{i\langle u,Y_t\rangle} \bigr\}
= \exp \biggl[t \biggl(\int_{\{|y|\le1\}}
e^{i\langle u,y\rangle}-1-i\langle u,y\rangle\pi(dy) \biggr) \biggr]
\]
and $(N_t)$ is a compound Poisson process with parameters
$\lambda=\pi(|y|>1)$ and $\mu(dy)=1_{\{|y|>1\}}\pi(dy)/\pi(|y|>1)$
($\lambda$ denotes the parameter for the waiting time between the
jumps of $N$ and $\mu$, the distribution of the jumps). Moreover,
$(W_t)$, $(Y_t)$ and $(N_t)$ are independent L\'evy processes.

Following this decomposition, we consider an $\mathbb{R}^d$-valued càdlàg
process $(X_t)$ solution to the SDE
%
\begin{equation}\label{edss}
dX_{t}=b(X_{t^{-}})\,dt+\sigma(X_{t^{-}})\,dW_{t}+\kappa(X_{t^{-}})\,dZ_{t},
\end{equation}
where $b\dvtx\mathbb{R}^d\mapsto\mathbb{R}^d$,
$\sigma\dvtx\mathbb{R}^d\mapsto\mathbb{M}_{d,l}$
(set of $d\times l$ real matrices) and
$\kappa\dvtx\mathbb{R}^d\mapsto\mathbb{M}_{d,l}$ are continuous with
sublinear growth and $(Z_t)$ is the
sum of the jump components of the L\'evy process: $Z_t=Y_t+N_t$.

In most papers dealing with L\'evy driven SDEs, the SDE reads
$dX_t=f(X_{t^-})\,dL_t$, where $(L_t)_{t\ge0}$ is a L\'evy process.
Here, we separate each part of the L\'evy process because they act
differently on the dynamical system. We isolate the drift term
because it usually produces the mean-reverting effect (which in
turn induces the ergodicity of the SDE). The two other terms are
both noises, but we distinguish them because they do not have the
same behavior.
\begin{Remarque} \label{abouttrunc}
In \eqref{edss} we chose to write
the jump component by
compensating the jumps smaller than 1, but it is obvious that, for
every $h>0$, $(X_t)$ is
also solution to
%
\begin{equation}\label{edssh}
dX_{t}=b^{h}(X_{t^{-}})\,dt
+ \sigma(X_{t^{-}})\,dW_{t}
+ \kappa (X_{t^{-}})\,dZ^{h}_{t}
\end{equation}
with $b^{h}=b+\int_{\{|y|\in(1,h]\}}y\pi(dy)$ if $h>1$,
$b^{h}=b-\int_{\{|y|\in(h,1]\}}y\pi(dy)$ if $h<1$, and
$Z_t^h=Y_t^h+N_t^h$, where the characteristic function of $Y_t^h$
is given for every $t\ge0$ by
\[
\mathbb{E}\bigl\{e^{i\langle u,Y_t^h\rangle}\bigr\}
= \exp \biggl[t \biggl(\int_{\{ |y|\le h\}} e^{i\langle u,y\rangle}
- 1-i\langle u,y\rangle\pi(dy) \biggr) \biggr]
\]
and $(N_t^h)$ is a compound Poisson process with parameters
$\lambda^h=\pi(|y|>h)$ and
$\mu^h(dy)=1_{\{|y|>h\}}\pi(dy)/\pi(|y|>h)$. By this remark, we
want to emphasize that the formulation \eqref{edss} is
conventional and that the coefficient $b$ in \eqref{edss} is
dependent on this choice. We will come back to this remark when we
introduce the assumptions of the main results where we want, on
the contrary, that they be intrinsic (see Remark
\ref{forcederappel}).
\end{Remarque}

Let us recall a result about existence, uniqueness and Markovian
structure of the solutions of \eqref{edss} (see \cite{bib1}).
\begin{theorem}\label{feller}
Assume that $b$, $\sigma$ and $\kappa$ are locally Lipschitz
functions with sublinear growth. Let $(\Omega,\mathcal{F},(\mathcal{F}_t),\mathbb{P})$
be a filtered probability space satisfying the
usual conditions and let $X_0$ be a random variable on
$(\Omega,\mathcal{F},\mathbb{P})$ with values in $\mathbb{R}^d$. Then,
for any
$(\mathcal{F}_t)$-Brownian motion $(W_t)_{t\ge0}$, for any $(\mathcal{F}_t)$-measurable
$(Z_t)_{t\ge0}$ as previously defined, the SDE
\eqref{edss} admits a unique c\`adl\`ag solution $(X_t)_{t\ge0}$
with initial condition $X_0$. Moreover, $(X_t)_{t\ge0}$ is a
Feller and Markov process.
\end{theorem}
\begin{Remarque}
L\'evy driven SDEs are the largest subclass of SDEs driven by
semimartingales such that the solutions have a Markovian
structure. Indeed, a result due to Jacod and Protter (see
\cite{bib13}) shows that, under appropriate conditions on the
coefficients, a stochastic process solution to a homogeneous SDE
driven by a semimartingale is a strong Markov process if, and only
if, the driving process is a L\'evy process.
\end{Remarque}

\subsection{Background on approximation of invariant
measures for Brownian diffusions} \label{back}

This problem has
already been studied by several authors when $(X_t)$ is a Brownian
diffusion, that is, when $\kappa=0$. In \cite{talay} Talay
approximates $\nu(f)$ by $\bar{\nu}_n^\gamma(f)=1/n\sum_{ k=1}^n
f(\bar{X}_{k-1}^\gamma)$, where $(\bar{X}_n^\gamma)_n$ denotes the
Euler scheme with constant step~$\gamma$. Denoting by
$\nu^\gamma$ the invariant distribution of the homogeneous Markov
chain $(\bar{X}_n^\gamma)_n$, he shows that
$\bar{\nu}_n^\gamma\stackrel{n\rightarrow+\infty}{\Longrightarrow}\nu^\gamma$
and that $\nu^\gamma\stackrel{\gamma\rightarrow0}{\Longrightarrow}\nu$,
under some uniform ellipticity and Lyapunov-type stability
assumptions. (A Markov process $(X_t)$ with infinitesimal
generator $A$ satisfies a Lyapunov assumption if there exists a
positive function $\mathcal{V}$ such that $\mathcal{V}(x)\rightarrow+\infty$ and
$\limsup AV(x)=-\infty$ when $|x|\rightarrow+\infty$. Then, $V$ is
called a Lyapunov function for $(X_t)$. Under this assumption,
$(X_t)$ admits a stationary, often ergodic when unique,
distribution. The existence of such a Lyapunov function depends on
the mean-reversion of the drift and on the intensity of the
diffusions term (see, e.g., \cite{borovkov,bib4,hasmin}
and \cite{bib14} for literature on Lyapunov stability).) In this
procedure, $\gamma$ and $n$ correspond to the two types of errors
that the discretization of this long time problem generates.
Practically, one cannot efficiently manage them together. Indeed,
when one implements this algorithm, one sets a positive real
$\gamma$ and then, one approximates the biased target
$\nu^\gamma$. In order to get rid of this problem, Lamberton and
Pag\`es (see \cite{bib2,bib3}) replace the standard Euler scheme
with constant step $\gamma$ with an Euler scheme with decreasing
step~$\gamma_n$. Denoting by $(\bar{X}_n)_{n\ge1}$ this Euler
scheme and by $(\eta_k)_{k\ge1}$ a sequence of weights such that
$H_n=\sum_{k=1}^n \eta_k\stackrel{n\rightarrow+\infty}{\longrightarrow}+\infty$,
they define a sequence of \textit{weighted} empirical measures
$(\bar{\nu}_n)$ and show under some Lyapunov assumptions (but
without ellipticity assumptions) that if $(\eta_n/\gamma_n)$ is
nonincreasing,
\[
\bar{\nu}_n(f)=\frac{1}{H_n}\sum_{k=1}^n \eta_k f(\bar{X}_{k-1})
\stackrel{n\rightarrow+\infty}{\longrightarrow}\nu(f)\qquad \mbox{p.s.},
\]
for every continuous function $f$ with polynomial growth (see
\cite{bib2,bib3} for more details and \cite{bib25} for extensions).
\begin{Remarque} 
These two approaches are
significantly different. Talay's method strongly relies on the
homogeneous Markovian structure of the constant step Euler scheme
and on its classical ``toolbox'' (irreducibility, positive
recurrence, \ldots, see, e.g., \cite{meyn}). Since the Euler
scheme with decreasing step is no longer homogeneous, Lamberton
and Pag\`es develop another method based on stability of Markov
chains and on martingale methods which can be extended to a
nonhomogeneous setting (see \cite{duflo}). This is why they do not
need any ellipticity assumptions on the coefficients.
\end{Remarque}

\subsection{Difficulties induced by the jumps of the L\'evy process}

In this paper we adapt the Lamberton and Pag\`es approach. In order
to obtain some similar results in the case of L\'evy
driven SDEs, one mainly has two kinds of obstacles to overcome.

From a dynamical point of view, the main difficulty comes from the
moments of the jump component. Indeed, by contrast with the case
of Brownian motion, the jump component can have only few moments
(stable processes, e.g.), and it then generates some
instability for the SDE.

The second obstacle appears in the simulation of the Euler scheme.
Actually, only in some very particular cases can the jump
component of a L\'evy process be simulated (compound Poisson
processes, stable processes$,\ldots$). In those cases, the
Euler scheme [that we call \textit{exact} Euler scheme and denote by (\ref{eA})]
can be built by using the true increments of $(Z_t)$. Otherwise, one
has to study some
\textit{approximate} Euler schemes where we replace the increments
of $Z_t$ with some approximations that can be simulated. The
canonical way for approximating the jump component is to truncate
its small jumps. Let $(u_n)_{n\ge1}$ be a sequence of positive
numbers such that $u_n<1$ and
$(u_n)$ decreases to 0 and $(Y^{n})_{n\ge1}$ be the sequence of
c\`adl\`ag processes defined by
\[
Y_t^{n}=\sum_{0<s\le t}\Delta Y_s 1_{\{\Delta
Y_s\in D_n\}}-t\int_{D_n}y\pi(dy)\qquad\forall t\ge0
\]
with $D_n=\{y\in\mathbb{R}^l, u_n<|y|\le1\}$ and $\Delta
Y_s=Y_s-Y_{s^-}$. The process $Y^{n}$ is a compensated compound
Poisson process with parameters $\lambda_n=\pi(D_n)$ and
$\mu_n(dy)=1_{D_n}(y)\pi(dy)/\pi(D_n)$. It converges locally
uniformly in $L^2$ to $Y$, that is, for every $T>0$,
%
\begin{equation}\label{convergencel2}
\mathbb{E}\biggl\{\sup_{0<t\le
T}|Y_t-Y^{n}_t|^2\biggr\}\stackrel{n\rightarrow+\infty}{\longrightarrow}0
\qquad\mbox{(see \cite{bertoin})}.
\end{equation}
We will denote by $Z^{n}$ the process defined by
$Z^{n}=Y^{n}+N$ and by (\ref{eB}) the Euler Scheme built with its
increments. The increments of $Z^n$ can be simulated if
$\lambda_n$ and the coefficient of the drift term can be
calculated, and if $\mu_n$ can be simulated for all $n\in\mathbb{N}$.
This is the case for a broad class of L\'evy processes, thanks to
classical techniques (rejection method, integral approximation$,\ldots$).
If there exists $(u_n)_{n\ge1}$ such that the increments
of $Z^n$ can be simulated for all $n\in\mathbb{N}$, we say that
\textit{the L\'evy measure can be simulated}. However, the
simulation time of $Z_t^n$ depends on the average number of its
jumps, that is, on $\pi(|y|>u_n)t$. When the truncation threshold
tends to 0 (this is necessary to approach the true increment of
the jump component), $\pi(|y|>u_n)$ explodes as soon as the L\'evy
measure is not finite (i.e., as soon as the true jump component is
not a compound Poisson process). It implies that the simulation
time of $Z_t^n$ explodes for a fixed $t$. However, thanks to the
decreasing step, it is possible to adapt the time step $\gamma_n$
and the truncation threshold~$u_n$ so that the expectation of the
number of jumps at each time step remains uniformly bounded.
Following the same idea, it is also possible to choose some steps
and some truncation thresholds so that the average number of jumps
at each time step tends to 0. In this case, approximating the true
component by the preceding compound Poisson process stopped at
its first jump time (the first time when it jumps) can also be
efficient [see Scheme (\ref{eC})].

\subsection{Construction of the procedures}

Let $(\gamma_n)_{n\ge1}$ be a decreasing sequence of
positive real numbers such that
${\lim} \gamma_n=0$ and $\Gamma_n=\sum_{k= 1}^n\gamma_k\rightarrow+\infty$
when $n\rightarrow+\infty$. Let
$(U_n)_{n\ge1}$ be a sequence of i.i.d. square integrable centered
$\mathbb{R}^l$-valued random variables with $\Sigma_{U_1}=I_d$. Finally,
let $(\bar{Z}_n)_{n\ge1}$, $(\bar{Z}^{B}_n)_{n\ge1}$ and
$(\bar{Z}^{C}_n)_{n\ge1}$ be sequences of independent
$\mathbb{R}^l$-valued random variables independent of $(U_n)_{n\ge1}$,
such that,
\[
\bar{Z}_n\stackrel{(\mathbb{R}^l)}{=}Z_{\gamma_n},\qquad
\bar{Z}^{B}_n\stackrel{(\mathbb{R}^l)}{=}Z^n_{\gamma_{n}}
\quad\mbox{and}\quad
\bar{Z}^{C}_n\stackrel{(\mathbb{R}^l)}{=}Z^n_{\gamma_{n}\wedge T^n}
\qquad\forall n\ge1,
\]
with $T^n=\inf\{s>0, |\Delta Z^n_s|>0\}.$ Let $x\in\mathbb{R}^d$. The
Euler Schemes (\ref{eA}), (\ref{eB}) and~(\ref{eC}) are recursively defined by
$\bar{X}_0=\bar{X}^{B}_0=\bar{X}^{C}_0=x$ and for every $n\ge
1$,
{\def\theequation{\Alph{equation}}
\setcounter{equation}{0}
\begin{eqnarray}
\bar{X}_{n+1} &=& \bar{X}_{n}+\gamma_{n+1} b(\bar{X}_{n})
+ \sqrt{\gamma_{n+1}}\sigma(\bar{X}_{n})U_{n+1}+\kappa(\bar{X}_{n})\bar{Z}_{n+1},  
\label{eA}\\
\bar{X}^{B}_{n+1} &=& \bar{X}^{B}_{n}+\gamma_{n+1}
b(\bar{X}^{B}_{n})+\sqrt{\gamma_{n+1}}\sigma(\bar{X}^{B}_{n})U_{n+1}
+ \kappa(\bar{X}^{B}_{n})\bar{Z}^{B}_{n+1}, 
\label{eB}\\
\bar{X}^{C}_{n+1} &=& \bar{X}^{C}_{n}+\gamma_{n+1}
b(\bar{X}^{C}_{n})+\sqrt{\gamma_{n+1}}\sigma(\bar{X}^{C}_{n})U_{n+1}
+ \kappa(\bar{X}^{C}_{n})\bar{Z}^{C}_{n+1}. 
\label{eC}
\end{eqnarray}}%
We set $\mathcal{F}_n=\sigma(\bar{X}_k,k\le n)$,
$\mathcal{F}^{B}_n=\sigma(\bar{X}^{B}_{k},k\le n)$
and $\mathcal{F}^{C}_n=\sigma(\bar{X}^{C}_{k},k\le n)$. Let
$(\eta_k)_{k\in\mathbb{N}}$ be a sequence of positive numbers such that
$H_n=\sum_{k=1}^n\eta_k\rightarrow+\infty$. For each scheme, we
define a sequence of weighted empirical measures by
%
\begin{eqnarray}
\bar{\nu}_n &=&\frac{1}{H_n}\sum_{k=1}^{n}\eta_k\delta_{\bar{X}_{k-1}},\qquad
\nonumber\\[-8pt]
\\[-8pt]
\bar{\nu}^{B}_n&=&\frac{1}{H_n}\sum_{k=1}^{n}\eta_k\delta_{\bar{X}^{B}_{k-1}}
\quad\mbox{and}\quad
\bar{\nu}^{C}_n=\frac{1}{H_n}\sum_{k=1}^{n}\eta_k\delta_{\bar{X}^{C}_{k-1}}.
\nonumber
\end{eqnarray}
For a function $f\dvtx\mathbb{R}^d\mapsto\mathbb{R}$, $(\bar{\nu}_n(f))$ can be
recursively computed (so is the case for the two other schemes).
Indeed, we have $\bar{\nu}_1(f)=f(x)$ and for every $n\ge1$,
\[
\bar{\nu}_{n+1}(f)=\bar{\nu}_n(f)+\frac{\eta_{n+1}}
{H_{n+1}}\bigl(f(\bar{X}_{n+1})-\bar{\nu}_n(f)\bigr).
\]

\textit{Some comments about the approximate Euler schemes}. In
Scheme (\ref{eB}), since $(u_n)$ decreases to 0, we discard fewer and
fewer jumps of the true component when $n$ grows. We will see in
Theorem \ref{principal} below that this is the only condition on
$(u_n)$ for the convergence of $(\bar{\nu}^{B}_n)$. This means
that we only need the law of $Z^n_{\gamma_n}$ to be an
``asymptotically good approximation'' of the law of $Z_{\gamma_n}$.
Yet, as previously mentioned, there is a hidden constraint induced
by the simulation time which is proportional to the average number
$\pi(|y|>u_n)\gamma_n$ of jumps of $Z^n$ on $[0,\gamma_n]$. In
practice, we require
$(\pi(|y|>u_n)\gamma_n)$ to be bounded.

Furthermore, if $\pi(|y|>u_n)\gamma_n\rightarrow0$
(i.e., the average number of jumps at each step tends to 0), we
will see that the first jump of $Z^n$ is all that matters for the
convergence of the empirical measure. This means that Scheme (\ref{eC})
becomes efficient.

\subsection{Notations} \label{notation}

Throughout this paper,
every positive real constant is denoted by $C$ (it may vary from
line to line). We denote the usual scalar product on $\mathbb{R}^d$ by
$\langle\cdot,\cdot\rangle$ and the Euclidean norm by $|\cdot|$. For any
$d\times l$ real matrix $M$, we define
$\|M\|= \sup_{\{|x|\le 1\}}|Mx|/|x|$. For a symmetric
${d\times d}$ real matrix $M$, we
set $\lambda_M=\max(0,\lambda_1,\ldots,\lambda_d)$, where
$\lambda_1,\ldots,\lambda_d$ denote the eigenvalues of $M$. For
every \mbox{$x\in\mathbb{R}^d$},
%
\begin{equation}\label{sym}
Mx^{\otimes2}=x^{*}Mx\le\lambda_M |x|^2.
\end{equation}
We denote by $\mathcal{C}_b(\mathbb{R}^d)$ [resp.
$\mathcal{C}_0(\mathbb{R}^d)$] the set of
bounded continuous functions on $\mathbb{R}^d$ with values in
$\mathbb{R}$ (resp. continuous
functions that go to $0$ at infinity) and $\mathcal{C}^2_K(\mathbb{R}^d)$, the set of
$\mathcal{C}^2$-functions on $\mathbb{R}^d$ with values in $\mathbb{R}$
and compact support. One says that $f$ is a $p$-H\"older function on $E$ with values in
$F$ (where $E$ and $F$ are normed vector spaces) if
\[
[f]_p=\sup_{x,y\in E} \frac{\|f(x)-f(y)\|_F}{\|x-y\|_E^p}<+\infty.
\]
Finally, we say that $V\dvtx\mathbb{R}^d\mapsto\mathbb{R}_+^*$ is an EQ-function
(\textit{Essentially Quadratic} function) if $V$ is $\mathcal{C}^2$,
$\lim V(x)=+\infty$ when $|x|\rightarrow+\infty$,
$|\nabla V|\le C \sqrt{V}$ and $D^2V$
is bounded. [In particular, $V$ given by $V(x)=\rho+Sx^{\otimes2}$,
where $\rho$ is a positive number and $S$ is a definite and positive
symmetric $d\times d$
real matrix, is an EQ-function.] For $p>0$, one checks
that $\|D^2 (V^p)\|\le CV^{p-1}$, that $V^p$ is a $2p$-H\"older function
if $p\le1/2$, and that $V^{p-1}\nabla V$ is a $(2p-1)$-H\"older function if
$p\in(1/2,1)$ (see Lemma \ref{lemmetech}). Hence, $\lambda_p$ and $c_p$ given by
%
\begin{eqnarray}\label{cplamb}
\lambda_p &:= & \frac{1}{2p}\sup_{x\in\mathbb{R}^d}
\lambda_{V^{1-p}D^2(V^p)(x)}
\quad \mbox{and}
\nonumber\\[-8pt]
\\[-8pt]
c_p &:=& \cases{%
\biggl[\dfrac{V^p}{p}\biggr]_{2p}, & \quad if $p\in(0,1/2]$,\cr
[V^{p-1}\nabla V]_{2p-1}, & \quad if $p\in (1/2,1]$}
\nonumber
\end{eqnarray}
are finite positive numbers.

\subsection{Organization of the paper}

The main results (Theorems \ref{principal} and
\ref{principal'}) are stated in Section \ref{section2} and are
proved in Sections \ref{section3},
\ref{section4} and \ref{section5}. First, we focus on the proof
of these theorems for the exact Euler Scheme (\ref{eA}): in Section
\ref{section3} we prove the almost sure tightness of
$(\bar{\nu}_n)$ and in Section \ref{section4} we establish that
every weak limiting distribution of $(\bar{\nu}_n)$ is invariant
for the SDE \eqref{edss}. Second, in Section~\ref{section5} we
point out the main differences which arise in the proofs when
considering the approximate Euler Schemes (\ref{eB}) and (\ref{eC}). In Section
\ref{section7} we show that the almost sure central limit theorem
for symmetric stable processes (see \cite{bib9}) can be obtained
as a consequence of our main theorems. Finally, in Section
\ref{simulations} we simulate the procedure on some concrete
examples.

\section{Main results}\label{section2}

In Theorem \ref{principal} we obtain a result under simple
conditions on the steps and on the weights. In Theorem
\ref{principal'} we show that, under more stringent conditions on
the steps and on the weights, some assumptions on the coefficients
of the SDE can be relaxed. Let us introduce the joint assumptions.
First, we state some assumptions on the moments of the L\'evy
measure at $+\infty$ and 0:
\[\hypertarget{eH}
(\mathrm{H^1_p})\dvtx \int_{|y|>1}\pi(dy)|y|^{2p}<+\infty,\qquad
(\mathrm{H^2_q})\dvtx \int_{|y|\le 1}\pi(dy)|y|^{2q}<+\infty,
\]
where $p$ is a positive real number and $q\in[0,1]$.

Assumption \hyperlink{eH}{$\mathrm{(H^1_p)}$} is satisfied if, and only if,
$\mathbb{E}|Z_t|^{2p}<+\infty$ for every $t\ge0$ (see \cite{bib18},
Theorem 6.1). By the compensation formula (see \cite{bertoin}),
\hyperlink{eH}{$\mathrm{(H^2_q)}$} is satisfied if and only if
$\mathbb{E}\{\sum_{0<s\le t}|\Delta Y_t|^{2q}\}<+\infty$, that is,
if and only if $(Y_t)$ has
locally \textit{$2q$-integrable variation}. We recall that
\hyperlink{eH}{$\mathrm{(H^2_q)}$} is always satisfied for $q= 1$ since
$\int_{\{|y|\le 1\}}|y|^2\pi(dy)<\infty$ for any L\'evy measure $\pi$.

Now, we introduce the Lyapunov assumption on the coefficients of
the SDE and on $\pi$ denoted by \ref{asS}. 
The parameter $a$ specifies the intensity of the mean-reversion. We
denote by $\tilde{b}$ the function defined by
\[
\tilde{b}= \cases{%
b, & \quad if $p\le1/2\le q$,\cr
\displaystyle b-\kappa\int_{\{|y|\le1\}}y\pi(dy), & \quad if $p,q\le1/2$,\cr
\displaystyle b+\kappa\int_{\{|y|> 1\}}y\pi(dy), & \quad if $p>1/2$.}
\]
The function $\tilde{b}$ plays the role of the global drift of the
dynamical system resulting from $b$ and from the jump component
(see Remark \ref{forcederappel} for more precisions). Let
$a\in(0,1]$, $p>0$ and $q\in[0,1]$.
{\def\theass{($\mathrm{S_{a,p,q}}$)}
\begin{ass}\label{asS}
There exists an EQ-function $V$ such that:

1. \textit{Growth control}: $|b|^2\le CV^a$,
\[
\cases{\operatorname{Tr}(\sigma\sigma^*)+\|\kappa\|^{2(p\vee q)}\le C V^{a+p-1},
& \quad if $p< 1$,\cr
\operatorname{Tr}(\sigma\sigma^*)+\|\kappa\|^{2}\le C V^{a},
& \quad if $p\ge1$.}
\]

2. \textit{Mean-reversion}: there exist $\beta\in\mathbb{R}$,
$\alpha>0$ such that,
$\langle\nabla V,\tilde{b}\rangle+\phi_{p,q}(\sigma,\kappa, \pi,\break V)\le
\beta-\alpha V^{a}$, where $\phi_{p,q}$ is given by
\begin{eqnarray*}
&& \phi_{p,q}(\sigma,\kappa,\pi,V)
\\
&&\qquad = \cases{%
c_{p} m_{2p,\pi}\|\kappa\|^{2p}V^{1-p}1_{q\le p}, & \quad if $p<1$,\cr
\lambda_1 \bigl(\operatorname{Tr}(\sigma\sigma^*)+m_{2,\pi}\|\kappa\|^2 \bigr), & \quad if $p=1$,\cr
d_{p}\lambda_p \biggl(\operatorname{Tr}(\sigma\sigma^*)+m_{2,\pi}\|\kappa\|^2+e_{p}
m_{2p,\pi}\dfrac{\|\kappa\|^{2p}}{V^{p-1}} \biggr), & \quad  if $p>1$,}
\end{eqnarray*}
with $m_{r,\pi}=\int|y|^r\pi(dy)$, $d_p=2^{(2(p-1)-1)_+}$, $c_p$
and $\lambda_p$ given by \eqref{cplamb}, and $e_p = {[\sqrt{V}]}^{2(p-1)}_1$.
\end{ass}}%

Assumption \ref{asS}.2 can be viewed as a
discretized version of ``$AV^p\le\beta-\alpha V^{a+p-1}$,''
where $A$ is the infinitesimal generator of $(X_t)$ defined on a
subset $\mathcal{D}(A)$ of $\mathcal{C}^2(\mathbb{R}^d)$ by
%
\begin{eqnarray}\label{1110}
\quad
Af(x)&=& \langle\nabla f,b\rangle(x)+\tfrac{1}{2}\operatorname{Tr}
(\sigma^* D^2 f\sigma)(x)
\nonumber\\[-8pt]
\\[-8pt]
&&{} +\int \bigl(f\bigl(x+\kappa(x)y\bigr)-f(x)-\langle\nabla
f(x),\kappa(x)y\rangle1_{\{|y|\le1\}} \bigr)\pi(dy).
\nonumber
\end{eqnarray}
Furthermore, one can check that if Assumption
\ref{asS} is fulfilled, then there exist
$\bar{\beta}\in\mathbb{R}$ and $\bar{\alpha}>0$ such that
``$AV^p\le \bar{\beta}-\bar{\alpha} V^{a+p-1}$.'' This means that if $V$ is
the function whose existence is required in Assumption
\ref{asS}, then $V^p$ is a Lyapunov function for the
stochastic process $(X_t)$
and for the Euler scheme $(\bar{X}_n)$.

The left-hand side of \ref{asS}.2 is the sum of two
antagonistic components: $\langle\nabla V,\tilde{b}\rangle$ produces
the mean-reverting effect (see Example~\ref{rmy} for concrete cases),
whereas the positive function $\phi_{p,q}$ expresses the noise
induced by the Brownian and jump components. In particular, if the
following tighter growth control condition holds,
%
\begin{equation}\label{hyp}
|b|^2\le CV^a, \qquad
\cases{%
\operatorname{Tr}(\sigma\sigma^*(x))\stackrel{|x|\rightarrow +\infty}{=}o\bigl(V^{a-(1-p)_+}(x)\bigr), \cr
\|\kappa(x)\|^2\stackrel{|x|\rightarrow+\infty}{=}o(V^{\eta_{a,p,q}}(x)),}
\end{equation}
with $\eta_{a,p,q}=(p\vee q)^{-1}(a+p-1)$ if $p\le1$ and
$\eta_{a,p,q}=a$ if $p>1$, then the term $\phi_{p,q}$
becomes negligible and the mean-reversion assumption becomes
\[
\langle\nabla V,\tilde{b}\rangle\le\beta-\alpha V^a.
\]
\begin{Remarque}\label{forcederappel}
If we had chosen to compensate the jumps smaller than $h>0$
rather than $h=1$, the corresponding assumption would have been
\hyperlink{asS}{$\mathrm{(S_{a,p,q}^h)}$}, where
\hyperlink{asS}{$\mathrm{(S_{a,p,q}^h)}$} is
obtained from \ref{asS} by replacing $b$ with $b^h$
and $\tilde{b}$ with $\tilde{b^h}$ defined by
\[
\tilde{b}^h= \cases{%
b^h, & \quad if $p\le1/2\le q$,\cr
\displaystyle b^h-\kappa\int_{\{|y|\le h\}}y\pi(dy), & \quad if $p,q\le1/2$,\cr
\displaystyle b^h+\kappa\int_{\{|y|> h\}}y\pi(dy), & \quad if $p>1/2$.}
\]
One can check that, for every $h>0$,
\hyperlink{asS}{$\mathrm{(S_{a,p,q}^h)}$}${}\Longleftrightarrow{}$\ref{asS}.
This means that these assumptions do not depend on the choice of
the truncation parameter $h$. Indeed, first, it is clear that
\hyperlink{asS}{$\mathrm{(S_{a,p,q}^h)}$}.1${}\Longleftrightarrow{}$\ref{asS}.1.
Second, when $p>1/2$ or \mbox{$p,q\le1/2$},
\hyperlink{asS}{$\mathrm{(S_{a,p,q}^h)}$}.2${}\Longleftrightarrow{}$\ref{asS}.2
because in these cases, $\tilde{b}^h=\tilde{b}$ for every $h>0$.
This can be explained by the existence of a formulation of the SDE
that does not depend on $h$. Actually, when $p>1/2$, we can
rewrite the SDE \eqref{edssh} by replacing $b^h$ with
$\tilde{b}^h$ and, $Z_t^h$ with
$\hat{Z}_t^h=Z_t^h-t\int_{\{|y|>h\}}y\pi(dy)$, that is, we can
compensate the big jumps. Since $(\hat{Z}_t^h)=(Z_t^\infty)$ for
every $h>0$, it follows that $\tilde{b}^h=\tilde{b}\,(=b^\infty)$
for every $h>0$. There also exists an intrinsic formulation when
$p,q\le1/2$ because in this case, we can replace $b^h$ with
$\tilde{b}^h$ and $Z_t^h$ with
$\check{Z}_t^h=Z_t^h+t\int_{\{|y|\le h\}}y\pi(dy)$ (now, we do not
compensate any jumps). Since $(\check{Z}_t^h)=(Z_t^0)$ for every
$h>0$, $\tilde{b}^h=\tilde{b}\,(=b^0)$. These formulations can be
considered as the natural formulations of the dynamical system in
these settings.

When $p\le1/2<q$, there is no intrinsic
formulation of the SDE (even if $\pi$ is symmetrical). Since $b^h$
depends on $h$, it appears that the left-hand side of
\hyperlink{asS}{$\mathrm{(S_{a,p,q}^h)}$}.2 also depends on $h$. However, under
the growth assumption on $\kappa$, one can check that $\langle\nabla
V,b^h\rangle=\langle\nabla V,b\rangle+o(V^a)$ and it follows that
the same
conclusion still holds in this case.
\end{Remarque}

We now state our first main result.
\begin{theorem}\label{principal}
Let $a\in(0,1]$, $p>0$ and $q\in[0,1]$ such
that \hyperlink{eH}{$\mathrm{(H^1_p)}$}, \hyperlink{eH}{$\mathrm{(H^2_q)}$} and
\ref{asS} are satisfied. Suppose that
$\mathbb{E}\{|U_1|^{2(p\vee1)}\}<+\infty$ and that the sequence
$(\eta_n/\gamma_n)_{n\ge1}$ is nonincreasing. Then:
\begin{longlist}[(3)]
\item[(1)] If $p/2+a-1>0$, the sequence $(\bar{\nu}_n)_{n\ge
1}$ is almost surely tight. Moreover, if
$\kappa(x)\stackrel{|x|\rightarrow+\infty}{=}o(|x|)$ and
$\operatorname{Tr}(\sigma\sigma^*)+\|\kappa\|^{2q}\le C V^{{p}/{2}+a-1}$, then
every weak limit of this sequence is an invariant probability for
the SDE \eqref{edss}. In particular, if $(X_t)_{t\ge0}$ admits
a unique invariant probability $\nu$, for every continuous
function $f$ such that $f=o(V^{{p}/{2}+a-1})$,
$\lim_{n\rightarrow\infty}\bar{\nu}_n(f)=\nu(f)$ a.s.
\item[(2)] The same result holds for $(\bar{\nu}^{B}_n)_{n\ge1}$.
\item[(3)] The same result holds for
$(\bar{\nu}^{C}_n)_{n\ge1}$ under the additional condition
%
\begin{equation}\label{condsupplementaire}
\pi({|y|>u_n})\gamma_n\stackrel{n\rightarrow+\infty}{\longrightarrow}0.
\end{equation}
\end{longlist}
\end{theorem}

We present below some examples which fulfill the
conditions of {Theorem }~{\ref{principal}}. In the first we
suppose that the dynamical system has a radial drift term and a
noise generated by a centered jump L\'evy process with a L\'evy
measure close to that of a symmetric stable process. In the
second we suppose that the SDE is only driven by a jump L\'evy
process, but we suppose that it is not centered. This implies that
even if the SDE has seemingly no drift term, a mean-reverting
assumption can be still satisfied.
\begin{Exemple} \label{rmy}
Let $\phi$ and $\psi$ be positive, bounded and continuous
functions on~$\mathbb{R}^d$ such that
\begin{eqnarray*}
\phi(x) &=& \phi(-x) \qquad \forall x\in\mathbb{R}^d,
\\
\underline{\phi} &=& \min_{\mathbb{R}^d}\phi(x)>0
\quad\mbox{and}\quad
\underline{\psi}=\inf_{\{|x|>1\}}\psi(x)>0.
\end{eqnarray*}
Consider $(Z_t)$
defined as in the SDE \eqref{edss} with L\'evy measure $\pi$
given by $\pi(dy)=\phi(y)/|y|^{d+r}\lambda_d(dy)$, where
$r\in(0,2)$.
When $\phi=C>0$, the increments of $(Z_t)$ can be exactly simulated
because $(Z_t)$ is a symmetric
$\mathbb{R}^d$-stable process with order~$r$. In the other cases,
$\bar{Z}_n^{B}$ and $\bar{Z}_n^{C}$ can be simulated by
the rejection method since
the density of $\pi$ is dominated by the density of a Pareto's law.

Let $\rho\in[0,2)$ and $b$ be a continuous function defined by
$b(x)=-\psi(x)x/|x|^\rho$. We consider $(X_t)$
solution to
%
\begin{equation}\label{edsexa}
dX_t=b(X_{t^-})\,dt +\kappa(X_{t^-})\,dZ_t,
\end{equation}
where
$\kappa$ is a continuous function such that $\|\kappa(x)\|^2\le
C(1+|x|^2)^\epsilon$ with $\epsilon\le1$. A natural candidate for
the function $V$ is $V(x)=1+|x|^2$. Indeed, since $\tilde{b}=b$
[because $\phi(y)=\phi(-y)$], one checks that there exists
$\beta\in\mathbb{R}$ such that
\[
\langle\nabla V(x),\tilde{b}(x)\rangle=-2\psi(x)|x|^{2-\rho}\le
\beta-\underline{\psi}V(x)^{1-{\rho}/{2}}.
\]
We set $a:=1-\rho/2$ and
\[
\Delta(r):=\{(p,q)\in (0,+\infty)\times[0,1],
\mbox{\hyperlink{eH}{$\mathrm{(H^1_p)}$} and \hyperlink{eH}{$\mathrm{(H^2_q)}$} hold}\}.
\]
We have $\Delta(r)=(0,r/2)\times
(r/2,1)$. By \eqref{hyp}, for every $(p,q)\in\Delta(r)$,
\ref{asS} is satisfied if $\epsilon<
(p+a-1)/q=(p-\rho/2)/q$. Hence, $(\bar{\nu}_n)$ is tight if there
exists $(p,q)\in\Delta(r)$ such that $p/2+a-1>0$, that is, if
$p>\rho$, and $(2p-\rho)/(2q)>\epsilon$. If, moreover,
$\|\kappa\|^{2q}\le C(1+|x|^2)^{{p}/{2}+a-1}$, that is, if
$(p-\rho)/(2q)\le\epsilon$, every weak limit $\nu$ is invariant for
the SDE \eqref{edsexa}.

It follows that if the invariant distribution $\nu$ is unique,
$\bar{\nu}_n\stackrel{\mathcal{L}}{\Longrightarrow}\nu$ a.s. as soon
as $2\rho<r$ and
$\epsilon<\sup\{(p-\rho)/(2q),(p,q)\in\Delta(r)\}=1/2-\rho/r$.
Furthermore, $\bar{\nu}_n(f)\rightarrow\nu(f)$ a.s. for every
continuous function $f$
satisfying $f(x)\le C(1+|x|)^\theta$ with $\theta\in[0,(r/2-\rho)/2)$.
\end{Exemple}
\begin{Exemple}
Let $\pi$ be a L\'evy measure on $\mathbb{R}$ such that
$\int_{|y|\le 1}|y|\pi(dy)<+\infty$, $\int_{|y|> 1}|y|^{2p}\pi(dy)<+\infty$ with
$p\ge2$ and $\int y\pi(dy)>0$. Let $(Z_t^0)$ be a real L\'evy
process with characteristic function given for every $t\ge0$ by
\[
\mathbb{E}\bigl\{e^{i\langle u,Z_t^0\rangle}\bigr\}
=\exp \biggl[t \biggl(\int\bigl(e^{i\langle u,y\rangle}-1\bigr)\pi(dy)
\biggr) \biggr].
\]
For instance, $(Z_t^0)$ can be a subordinator with no drift term.
We assume that $\kappa(x)=-\psi(x)x/|x|^{\rho}$ with
$\rho\in[0,2)$ and $\psi$ defined as in the preceding example. We
then consider the SDE:
\[
dX_t=\kappa(X_{t^-})\,dZ_t^{0}=b(X_{t^-})\,dt+\kappa(X_{t^-})\,dZ_t,
\]
with $b(x)=\kappa(x)\int_{\{|y|\le1\}}y\pi(dy)$ and
$Z_t=Z_t^0-t\int_{\{|y|\le1\}}y\pi(dy)$. Since $p>1/2$,
$\tilde{b}(x)=b(x)+\kappa(x)\int_{\{|y|>1\}}
y\pi(dy)=\kappa(x)\int y\pi(dy)$. Setting $V(x)=1+x^2$, one checks
that there exists $\beta\in\mathbb{R}$ such that
\[
V'(x)\tilde{b}(x)=-2\psi(x)\int y\pi(dy)|x|^{2-\rho}\le\beta
-\underline{\psi}\int y\pi(dy)V(x)^{1-{\rho}/{2}}.
\]
We set $a=1-\rho/2$. Let $p\ge 2$ and $q\le1/2$ such that
\hyperlink{eH}{$\mathrm{(H^1_p)}$} and \hyperlink{eH}{$\mathrm{(H_q^2)}$} hold. First,
checking that as soon as $\rho>0$, $\|\kappa(x)\|^2= o(1+|x|^2)^a$
when $|x|\rightarrow+\infty$, we derive from \eqref{hyp} that
\ref{asS} is satisfied as soon as $\rho\in(0,2)$.
Second, for every $p\ge2$, $q\le1$ and $a\in(0,1)$, one can check
that $p/2+a-1>0$
and $\|\kappa(x)\|^{2q}\le C(1+|x|^2)^{{p}/{2}+a-1}$.
Hence, Theorem~\ref{principal} applies for every $\rho\in(0,2)$.
\end{Exemple}

The interest of Theorem \ref{principal} lies in the facility with
which it can be put to use in concrete situations.
For instance, in Scheme (\ref{eA}), we only have to take
a sequence $(\gamma_n)_{n\ge1}$ decreasing to $0$,
with infinite sum and $\eta_n=\gamma_n$.
The next theorem (Theorem \ref{principal'}) requires tougher conditions
on the sequences $(\gamma_n)$ and $(\eta_n)$,
but it can be applied to SDEs where the coefficients do not
necessarily verify all conditions of Theorem \ref{principal}. It
broadens the class of SDEs for which
we can find an efficient procedure for the approximation of the
invariant measure.
\begin{theorem}\label{principal'}
Let $a\in(0,1]$, $p>0$ and $q\in [0,1]$ such
that \hyperlink{eH}{$\mathrm{(H^1_p)}$}, \hyperlink{eH}{$\mathrm{(H^2_q)}$} and
\ref{asS} are satisfied. Suppose that
$\mathbb{E}\{|U_1|^{2(p\vee1)}\}<+\infty$. Then:
\begin{longlist}[(3)]
\item[(1)] Let $s\in(1,2]$ satisfying the following additional
conditions when $p>1/2$:
%
\begin{eqnarray}\label{cond1}
\cases{%
s>\dfrac{2p}{2p+(a-1)(2p-1)/p}, & \quad if $\dfrac{1}{2}<p\le1$,\cr
s>\dfrac{2p}{2p+a-1}, & \quad if $p\ge1$.}
\end{eqnarray}
If $p/s+a-1>0$, there exist some sequences $(\gamma_n)_{n\ge1}$
and $(\eta_n)_{n\ge1}$ such that $(\bar{\nu}_n)_{n\ge1}$ is
almost surely tight. Moreover, if
$\kappa(x)\stackrel{|x|\rightarrow+\infty} {=}o(|x|)$ and
$\operatorname{Tr}(\sigma\sigma^*)+\|\kappa\|^{2q}\le C V^{{p}/{s}+a-1}$, then
every weak limit of this sequence is an invariant probability for
the SDE \eqref{edss}. In particular, if $(X_t)_{t\ge0}$ admits a
unique invariant probability $\nu$, for every continuous function
$f$ such that $f=o(V^{{p}/{s}+a-1})$,
$\lim_{n\rightarrow\infty}\bar{\nu}_n(f)=\nu(f)$ a.s.
\item[(2)] The same result holds for $(\bar{\nu}^{B}_n)_{n\ge 1}$.
\item[(3)] The same result holds for
$(\bar{\nu}^{C}_n)_{n\ge1}$ under the additional condition
\eqref{condsupplementaire}.
\end{longlist}
\end{theorem}
\begin{Remarque}
The sequences $(\eta_n)_{n\ge1}$ and $(\gamma_n)_{n\ge1}$
must verify an explicit condition given in Proposition
\ref{tension} below (see Remark \ref{pas2} for a version adapted
to polynomial steps and weights).
\end{Remarque}

In the following example, we consider the same class of
SDEs as in Example~\ref{rmy} in the nonintegrable case (i.e.,
$r\le1$). One can observe that the mean-reversion condition and
the growth condition on $\kappa$ and on the functions $f$ whose
the procedure converges can be relaxed. We also give some explicit
polynomial weights and steps for which Theorem \ref{principal'}
applies in this case.
\begin{Exemple}\label{tpr} 
Let $\rho\in[0,2)$ and $r\in(0,1]$, let $b$ and $\kappa$ be continuous
functions defined as in Example \ref{rmy}. Consider $(X_t)$ solution to
the SDE \eqref{edsexa} and assume that the invariant measure $\nu$
is unique. For $s\in(1,2]$, denote by $(\gamma_{n,s})$ and
$(\eta_{n,s})$ some sequences of steps and weights satisfying
$\gamma_{n,s}=Cn^{-r_1}$, $\eta_{n,s}=Cn^{-r_2}$ with $r_1\le r_2$
and
\[
0<r_1<2\biggl(1-\frac1 s\biggr) \quad\mbox{and}\quad
r_2<1\quad\mbox{or}\quad
0<r_1\le2\biggl(1-\frac1 s\biggr)\quad\mbox{and}\quad r_2=1.
\]
Then, for these choices of steps and weights,
$\bar{\nu}_n\stackrel{\mathcal{L}}{\Longrightarrow}\nu$ {a.s.} as soon
as $s\rho<r$ and $\epsilon\in[0, 1/s-\rho/r)$ (this improves the
condition: $2\rho<r$ and $\epsilon\in[0, 1/2-\rho/r)$ of Example
\ref{rmy}). Furthermore, $\bar{\nu}_n(f)\rightarrow\nu(f)$ {a.s.}
for every continuous function $f$ satisfying $|f(x)|\le
C(1+|x|)^\theta$ with $\theta\in[0,(r/s-\rho)/2)$ (this improves
the condition: $\theta\in[0,(r/2-\rho)/2)$ of Example \ref{rmy}).
\end{Exemple}

\section{Almost sure tightness of $(\bar{\nu}_n(w,dx))_{n\in\mathbb{N}}$}\label{section3}

The main result of this section is Proposition
\ref{tension}. We need to introduce the function $f_{a,p}$
defined for all $s\in (1,2]$ by
%
\begin{equation}\label{puis}
f_{a,p}(s)= \cases{%
s, & \quad if $s\ge2p$,\cr
\dfrac{p+a-1}{{p}/{s}+{(a-1)}/{(2(p\wedge1))}}\wedge s, & \quad if $s<2p$.}
\end{equation}
Assume that $p+a-1>0$. Then, $s\mapsto f_{a,p}(s)$ is a
nondecreasing function which satisfies $f_{1,p}(s)=s$ for all
$p>0$ and $f_{a,p}(2)=2$. Note that $f_{a,p}(s)>1$ if and only if
$s$ satisfies assumption \eqref{cond1}.
\begin{prop}\label{tension}
Let $a\in(0,1]$, $p>0$ and $q\in(0,1]$ such that
\hyperlink{eH}{$\mathrm{(H^1_p)}$}, \hyperlink{eH}{$\mathrm{(H^2_q)}$} and \ref{asS}
are satisfied. Assume that $\mathbb{E}\{|U_1|^{2(p\vee1)}\}<+\infty$ and
$(\eta_n/\gamma_n)_{n\ge1}$ is nonincreasing.
\begin{longlist}[(2)]
\item[(1)] Then,
\[
\sup_{n\ge1} \bar{\nu}_n(V^{{p}/{2}+a-1})<+\infty\qquad \mbox{a.s.}
\]
Consequently, if $\frac{p}{2}+a-1>0,$ the sequence
$(\bar{\nu}_n)_{n\in\mathbb{N}}$ is a.s. tight.
\item[(2)] Let $s\in(1,2)$ such that assumption \eqref{cond1} is
satisfied. Assume that $(\eta_n)_{n\ge1}$ and $(\gamma_n)_{n\ge 1}$ are such that
%
\begin{equation}\label{hop}
\qquad
\biggl(\frac{1}{\gamma_n} \biggl(\frac{\eta
_n}{H_n\sqrt{\gamma_n}} \biggr)^{^{f_{a,p}(s)}} \biggr)
\mbox{ is nonincreasing and }
\sum_{n\ge 1} \biggl(\frac{\eta_n}{H_n\sqrt{\gamma_n}}\biggr)^{^{f_{a,p}(s)}}<+\infty.
\end{equation}
Then, $\sup_{n\ge1}\bar{\nu}_n(V^{{p}/{s}+a-1})<+\infty$
\textit{a.s.} and the sequence $(\bar{\nu}_n)_{n\in\mathbb{N}}$ is a.s.
tight as soon as $p/s+a-1>0$.
\end{longlist}
\end{prop}
\begin{Remarque}\label{pas2} 
If $\gamma_n=Cn^{-r_1}$ and
$\eta_n=C n^{-r_2}$ with $r_1\le r_2$, then assumption~\eqref{hop}
reads
%
\begin{eqnarray}\label{cond2}
r_2 &<& 1 \quad \mbox{and}\quad
0<r_1<\bar{r}_1:=2 \biggl(1-\frac{1}{f_{a,p}(s)} \biggr)
\quad \mbox{or}\quad
\nonumber\\[-8pt]
\\[-8pt]
r_2 &=& 1 \quad \mbox{and}\quad
0<r_1\le\bar{r}_1.
\nonumber
\end{eqnarray}
\end{Remarque}
The proof of Proposition \ref{tension} is organized as
follows: first, in Section \ref{4.1} (see Proposition
\ref{prop9}) we establish a fundamental recursive control of the
sequence $(V^p(\bar{X}_n))$: we show
that
\hyperlink{e15}{$\mathrm{(R_{a,p})}$}: There exist $n_0 \in\mathbb{N}$,
$\alpha'>0$, $\beta'>0$ such that $\forall n\ge n_0$,
%
\begin{equation}\hypertarget{e15}
\mathbb{E}\{V^p(\bar{X}_{n+1})|\mathcal{F}_n\}
\le V^p(\bar{X}_n)+\gamma_{n+1}V^{p-1}(\bar{X}_n)
\bigl(\beta' -\alpha'V^a(\bar{X}_n) \bigr).
\end{equation}
For this step, we rely on Lemma \ref{moment} that provides a
control of the moments of the increments of the jump component in
terms of $p$ and $q$.

Second, in Section \ref{conseqprop} we make use of
martingale techniques in order to derive some consequences from
\hyperlink{e15}{$\mathrm{(R_{a,p})}$}. In Lemma \ref{lemme4} we establish a
$L^p$-control of the Euler scheme with arguments close to
\cite{bib3}. This control is fundamental for the proof of
Corollary \ref{cor} where we show the following property:

$\mathrm{(C_{p,s})}$:\hypertarget{eCps}\ There exist $\rho\in(1,2]$ and a sequence
$(\pi_n)$ of $\mathcal{F}_n$-measurable random variables such that
%
\begin{equation}\label{assumcpsrev}
\sum_{n\ge 1}\biggl(\frac{\eta_n}{H_n\gamma_n}\biggr)^{\rho}
\mathbb{E}\{|V^{{p}/{s}}(\bar{X}_n)-\pi_{n-1}|^{\rho}\}< +\infty.
\end{equation}
This step is used to obtain a $L^\rho$-martingale control (see
proof of Lemma \ref{explication}). We will see in the proof of
Corollary \ref{cor} that the choice of the sequence $(\pi_n)$
depends on $p$ and $q$. In particular, even if $q$ does not appear
in \hyperlink{eCps}{$\mathrm{(C_{p,s})}$}, this assumption indirectly depends on
this parameter. The same remark holds for
\hyperlink{e15}{$\mathrm{(R_{a,p})}$}. In
the following lemma, we show that these two steps are all what we
have to show for the proof of Proposition \ref{tension}.
\begin{lemme}\label{explication}
Let $p>0$, $a\in(0,1]$ and $s\in(1,2]$ such that
\hyperlink{eH}{$\mathrm{(H^1_p)}$}, \hyperlink{e15}{$\mathrm{(R_{a,p})}$} and
\hyperlink{eCps}{$\mathrm{(C_{p,s})}$} are
fulfilled. Assume
that $\mathbb{E}\{|U_1|^{2(p\vee1)}\}<+\infty$ and that
$(\eta_n/\gamma_n)$ is nonincreasing. Then,
%
\begin{equation}\label{conclusionrev}
\sup_{n\ge1}\bar{\nu}_n (V^{{p}/{s}+a-1})<+\infty\qquad \mbox{a.s.}
\end{equation}
and the sequence $(\bar{\nu}_n)_{n\ge1}$ is \textit{a.s.} tight as soon
as $p/s+a-1>0$.
\end{lemme}
\begin{pf}
By a convexity argument (see Lemma 3 of \cite{bib3}), one shows
that \hyperlink{e15}{$\mathrm{(R_{a,p})}$}${}\Longrightarrow{}$\hyperlink{e15}{($\mathrm{R_{a,\bar{p}}}$)}
for all $\bar{p}\in(0,p]$. Hence, for all $s\in(1,2]$, there exists $
n_0\,\in\mathbb{N}$, $\hat\alpha>0$ and $\hat\beta>0$ such that
$\forall k\ge n_0$,
%
\begin{equation}\label{p/s1}
\qquad
\mathbb{E}\{V^{{p}/{s}}(\bar{X}_{k})|\mathcal{F}_{k-1}\}
\le V^{{p}/{s}}(\bar{X}_{k-1})+\gamma_{k}V^{{p}/{s}-1}(\bar{X}_{k-1})
\bigl(\hat\beta-\hat\alpha V^a(\bar{X}_{k-1}) \bigr).
\hspace*{-6pt}
\end{equation}
For $R>0$, set
$\varepsilon(R)=\sup_{\{|x|>R\}}V^{-a}(x)$ and
$M(R)=\sup_{\{|x|\le R\}}V^{p/s -1}(x)$. We have
%
\begin{equation}\label{closeargumentrev}
V^{{p}/{s}-1}(x)\le\varepsilon(R) V^{{p}/{s}+a-1}(x)+M(R).
\end{equation}
Since
$V(x)\rightarrow+\infty$ when $|x|\rightarrow+\infty$ (resp. since
$V$ is bounded on every compact set), $\varepsilon(R)\rightarrow0$
when $R\rightarrow+\infty$ [resp. $M(R)$ is finite for every
$R>0$]. Hence, for every $\varepsilon>0$, there exists
$M_\varepsilon>0$ such that $V^{{p}/{s}-1}\le\varepsilon V^{{p}/{s}+a-1}+M_\varepsilon$.
By setting $\varepsilon=\hat{\alpha}/(2\hat{\beta})$,
$\tilde{\alpha}=\hat{\alpha}/{2}$ and
$\tilde{\beta}=\hat{\beta}M_{\varepsilon}$, we deduce that
$V^{{p}/{s}-1}(\hat{\beta}-\hat{\alpha}V^a)\le \tilde{\beta}-\tilde{\alpha}V^{{p}/{s}+a-1}$.
Hence, we derive from \eqref{p/s1} that
\[
V^{{p}/{s}+a-1}(\bar{X}_{k-1})\le
\frac{V^{{p}/{s}}(\bar{X}_{k-1})- \mathbb{E}\{V^{{p}/{s}}(\bar{X}_k)|\mathcal{F}_{k-1}\}}
{\tilde{\alpha}\gamma_k}+\frac{\tilde{\beta}}{\tilde \alpha}
\qquad\forall k\ge n_0.
\]
It follows that \eqref{conclusionrev} holds if
%
\begin{eqnarray}\label{topod}
&& \sup_{n\ge n_0+1} \Biggl(\frac{1}{H_n}\sum_{k=n_0+1}^{n}
\frac{\eta_k}{\gamma_k} \bigl(V^{{p}/{s}}(\bar{X}_{k-1})-
\mathbb{E}\{V^{{p}/{s}}(\bar{X}_k)|\mathcal{F}_{k-1}\} \bigr) \Biggr)
\nonumber\\[-8pt]
\\[-8pt]
&&\qquad < + \infty\qquad \mbox{a.s.}
\nonumber
\end{eqnarray}
We then prove \eqref{topod}. We decompose the above sum as follows:
\begin{eqnarray*}
&& \frac{1}{H_n}\sum_{k=n_0 +1}^{n}
\frac{\eta_k}{\gamma_k}\bigl(V^{{p}/{s}}
(\bar{X}_{k-1})-\mathbb{E}\{V^{{p}/{s}}(\bar{X}_k)|\mathcal{F}_{k-1}\}\bigr)
\\
&&\qquad =-\frac{1}{H_n}\sum_{k=n_0+1}^{n}
\frac{\eta_k}{\gamma_k}\Delta V^{{p}/{s}}(\bar{X}_k)
\\
&&\qquad \quad {}
+\frac{1}{H_n}\sum_{k=n_0+1}^{n}
\frac{\eta_k}{\gamma_k} \bigl(V^{{p}/{s}}(\bar{X}_{k})
- \mathbb{E}\{V^{{p}/{s}}(\bar{X}_k)|\mathcal{F}_{k-1}\} \bigr),
\end{eqnarray*}
where $\Delta V^{{p}/{s}}(\bar{X}_k)=V^{{p}/{s}}(\bar{X}_{k})-V^{{p}/{s}}(\bar{X}_{k-1})$.
First, an Abel's transform yields
\begin{eqnarray*}
-\frac{1}{H_n}\sum_{k=n_0+1}^{n}\frac{\eta_k}{\gamma_k}\Delta
V^{{p}/{s}}(\bar{X}_k)
&=& \frac{1}{H_n} \biggl(\frac{\eta_{n_0}}{\gamma_{n_0}}V^{{p}/{s}}(\bar{X}_{n_0})
- \frac{\eta_{n}}{\gamma_{n}}V^{{p}/{s}}(\bar{X}_{n}) \biggr)
\\
&&{} +\frac{1}{H_n} \Biggl(\sum_{k=n_{0+1}}^{n}
\biggl(\frac{\eta_k}{\gamma_k}-\frac{\eta_{k-1}}{\gamma_{k-1}}\biggr)
V^{{p}/{s}}(\bar{X}_{k-1}) \Biggr)
\\
&\le & \frac{\eta_{n_0}}{H_n\gamma_{n_0}}V^{{p}/{s}}(\bar{X}_{n_0}),
\end{eqnarray*}
where we used in the last inequality that $(\eta_n/\gamma_n)$ is
nonincreasing. Hence, since $H_n\stackrel{n\rightarrow+\infty}{\longrightarrow}+\infty$ and
$\frac{\eta_{n_0}}{H_n\gamma_{n_0}}V^{{p}/{s}}(\bar{X}_{n_0})\stackrel{n\rightarrow+\infty}{\longrightarrow}0$
a.s.,
%
\begin{equation}\label{erpe2}
\sup_{n\ge n_0} \Biggl(-\frac{1}{H_n}\sum_{k=n_0+1}^{n}
\frac{\eta_k}{\gamma_k}\Delta V^{{p}/{s}}(\bar{X}_{k}) \Biggr)
<+\infty\qquad \mbox{a.s.}
\end{equation}
Second, one denotes by $(M_n)_{n\in\mathbb{N}}$ the martingale
defined by
%
\begin{equation}\label{revmartingale}
M_n=\sum_{k=1}^{n}\frac{\eta_k}{H_k\gamma_k}
\bigl(V^{{p}/{s}}(\bar{X}_{k})
- \mathbb{E}\{V^{{p}/{s}}(\bar{X}_k)|\mathcal{F}_{k-1}\} \bigr).
\end{equation}
Let $\rho\in(1,2]$ and $(\pi_{k})$ be a sequence
of $\mathcal{F}_{k}$-measurable random variables such
that~\eqref{assumcpsrev} holds. We derive from the elementary\vadjust{\goodbreak}
inequality $|u+v|^\rho\le2^{\rho-1}(|u|^\rho+|v|^\rho)$ that
\begin{eqnarray*}
&& \mathbb{E} \bigl\{\bigl|V^{{p}/{s}}(\bar{X}_{k})
- \mathbb{E}\{V^{{p}/{s}}(\bar{X}_k)|\mathcal{F}_{k-1}\}\bigr|^{\rho} \bigr\}
\\
&&\qquad \le C\mathbb{E} \{|V^{{p}/{s}}(\bar{X}_{k})-\pi_{k-1}|^{\rho} \}
+ C\mathbb{E} \bigl\{\bigl|\mathbb{E} \bigl\{\bigl(\pi_{k-1}- V^{{p}/{s}}
(\bar{X}_k)\bigr)|\mathcal{F}_{k-1} \bigr\}\bigr|^{\rho} \bigr\}
\\
&&\qquad
\le C\mathbb{E} \{|V^{{p}/{s}}(\bar{X}_{k})-\pi_{k-1}|^{\rho} \},
\end{eqnarray*}
thanks to the Jensen inequality. Hence,
\hyperlink{eCps}{$\mathrm{(C_{p,s})}$} yields
$\sum_{k\ge1}\mathbb{E}\{|\Delta M_k|^\rho\}<+\infty$ a.s. Since
$\rho>1$, it follows from Chow's theorem (see \cite{hall}) that
$M_n\stackrel{n\rightarrow\infty}{\longrightarrow} M_\infty$ a.s. where
$M_\infty$ is finite a.s. Then, Kronecker's lemma yields
%
\begin{equation}\label{erpe}
\frac{1}{H_n}\sum_{k=n_0+1}^{n}\frac{\eta_k}{\gamma_k}
\bigl(V^{{p}/{s}}(\bar{X}_{k})-
\mathbb{E}\{V^{{p}/{s}}(\bar{X}_k)|\mathcal{F}_{k-1}\} \bigr)
\stackrel{n\rightarrow\infty}{\longrightarrow} 0\qquad \mbox{a.s.}
\end{equation}
Hence, \eqref{topod} follows from \eqref{erpe2} and \eqref{erpe}.
Finally, since
$\lim_{|x|\rightarrow+\infty}V^{{p}/{s}+a-1}(x)=+\infty$ when
$p/s+a-1>0$, we derive from a classical tightness criteria (see,
{e.g.}, \cite{duflo}, page~41)
that $(\bar{\nu}_n)_{n\ge1}$ is {a.s.} tight as soon as $p/s+a-1>0$.
\end{pf}

\subsection{A recursive stability relation}\label{4.1}

\begin{prop}\label{prop9}
Let $p>0$, $q\in[0,1]$ and $a\in(0,1]$. Assume
\hyperlink{eH}{$\mathrm{(H^1_p)}$},
\hyperlink{eH}{$\mathrm{(H_q^2)}$} and \ref{asS}. If, moreover,
$\mathbb{E}\{|U_1|^{2(p\vee1)}\}<+\infty$, then
\hyperlink{e15}{$\mathrm{(R_{a,p})}$} holds.
\end{prop}

The idea of the proof of Proposition \ref{prop9} is to
obtain an inequality of the following type:
\[
\mathbb{E}\{V^p(\bar{X}_{n+1})-V^p(\bar{X}_n)|\mathcal{F}_n\}\le\gamma
_{n+1}pV^{p-1}(\bar{X}_n)\Phi(\bar{X}_n)+R_n,
\]
where $\Phi=\langle\nabla
V,\tilde{b}\rangle+\phi_{p,q}(\sigma,\kappa,\pi,V)$ [see
\ref{asS}.2] and $R_n$ is asymptotically negligible
in a sense made clear in the proof. To this end, we begin by three
lemmas. In Lemma \ref{moment} we study the behavior of the
moments of $(Z_t)$ near 0. Then, in Lemma \ref{lemmetech}, we
state some properties of the derivatives of $V^p$ in terms of $p$
and in the last one (Lemma \ref{controljumpcomp}) we control the
contribution of the jump component on the conditional
expectation (conditioned by $\mathcal{F}_n$) of the increment
$V^p(\bar{X}_{n+1})-V^p(\bar{X}_n)$.
\begin{lemme}\label{moment}
\textup{(i)} Let $p>0$ such that \hyperlink{eH}{$\mathrm{(H^1_p)}$} holds. Then, for every
$h>0$, there exists a locally bounded function $\psi_h$ such that
%
\begin{equation}\label{bigjumprev}
\forall t\ge0 \qquad
\mathbb{E}\{|N^h_t|^{2p}\}= \int_{|y|>h}|y|^{2p}\pi (dy)
\bigl(t+\psi_h(t) t^2\bigr) .
\end{equation}
{\smallskipamount=0pt
\begin{longlist}[(iii)]
\item[(ii)] Let $q\in[0,1]$ such that
\hyperlink{eH}{$\mathrm{(H_q^2)}$} holds. Then, for every $h>0$,
\[
\cases{%
\displaystyle
\mathbb{E}\biggl\{\biggl|Y^h_t
+ t\int_{|y|\le h}y\pi(dy)\biggr|^{2q}\biggr\}
\le t\int_{|y|\le h}|y|^{2q}\pi(dy), & \quad if $q\le1/2$,\cr
\displaystyle
\mathbb{E}\{|Y^h_t|^{2q}\}\le C_q t
\int_{|y|\le h}|y|^{2q}\pi(dy), & \quad if $q\in(1/2, 1]$.}
\]
\item[(iii)] Let $p\in[1,+\infty)$ such that \hyperlink{eH}{$\mathrm{(H^1_p)}$} holds.
Then, there
exists $\eta>1$ such that, for every $T>0$, for every
$\varepsilon>0$, there exists $C_{\varepsilon,T,p}>0$ such that,
\[
\forall t\in[0,T]\qquad
\mathbb{E}\{|\hat{Z}_t|^{2p}\}\le t
\biggl(\int|y|^{2p}\pi(dy)+\varepsilon \biggr)
+C_{\varepsilon,T,p} t^\eta,
\]
where $(\hat{Z}_t)$ is the compensated jumps process defined by
$\hat{Z}_t=Z_t-t\int_{|y|>1}y\pi(dy)$. In particular,
$\mathbb{E}|\hat{Z}_t|^{2}= t\int|y|^{2}\pi(dy)$.
\end{longlist}}
\end{lemme}
\begin{Remarque}\label{roughness} 
In this lemma we obtain, in particular, a
control of the expansion of $t\mapsto\mathbb{E}\{|D_t|^r\}$ in the
neighborhood of 0 (where $D$ denotes one of the above jump
components and $r$, a positive number). We have the following type
of inequality: $\mathbb{E}\{|D_t|^{r}\}\le c_r t+O(t^\eta)$, where $c_r$
is a nonnegative real constant and $\eta>1$. In the first and in
the last inequality, we minimize this value because it has a
direct impact on the coefficients of the function $\phi_{p,q}$ and
then, on the mean-reverting assumption (see Lemma
\ref{controljumpcomp} for details). Note that we cannot have
$c_r=0$ in the inequalities of Lemma \ref{moment}. Indeed,
according to the Kolmogorov criterion, a L\'evy process $D$ that
satisfies $\mathbb{E}\{|D_t|^{r}\}\le Ct^\eta$ in the neighborhood of 0
is pathwise continuous [for the Brownian motion, $c_r=0$ as soon
as $r>2$ since $\mathbb{E}\{|W_t|^r\}=o(t^{{r}/{2}}$)]. When $p>1$,
this feature generates a specific contribution of the jump
component on the mean-reverting assumption
(\ref{asS}.2). This contribution appears in
$\phi_{p,q}$ where there is an additional term of order $2p$
coming only from the jump component.
\end{Remarque}
\begin{pf*}{Proof of Lemma \ref{moment}}
(i) $(N_t^h)_{t\ge0}$ is a compound Poisson process with
parameters $\lambda_h=\pi(|y|>h)$ and
$\mu^h(dy)={1_{\{|y|>h\}}\pi(dy)}/{\pi(|y|>h)}$. Hence, $(N_t^h)$
can be written as follows: $N_t^h={\sum}_{n\ge1} R_n 1_{T_n\le
t}$, where $(R_n)_{n\ge1}$ is a sequence of i.i.d. r.v. with law
$\mu^h$ and $ (T_n)_{n\in\mathbb{N}}$ is the sequence of the jump
times of a Poisson process with parameter ${\lambda_h}$
independent of $(R_n)_{n\ge1}$. We have
\begin{eqnarray}
\mathbb{E}\{|N_t^h|^{2p}\} &=& \sum_{n\ge1}\mathbb{E}
\Biggl\{\Biggl|\sum_{i=1}^n R_i\Biggr|^{2p}\Biggr\}
e^{-{\lambda_h} t}\frac{({\lambda_h} t)^n}{n!}
= {\lambda_h} t\mathbb{E}\{|R_1|^{2p}\}F_{{\lambda_h}}(t)
\nonumber\\
\eqntext{\displaystyle \mbox{where } F_{{\lambda}} (t)=e^{-{\lambda}
t}\sum_{n\ge0}\frac{\mathbb{E}\{|\sum_{i=1}^{n+1}
R_i|^{2p}\}}{\mathbb{E}\{|R_1|^{2p}\}}\frac{({\lambda} t)^n}{(n+1)!}.}
\end{eqnarray}
By the elementary inequality (this inequality will be
usually needed in the sequel for the control of the moments of
some sums of jumps)
%
\begin{equation}\label{eqelem2}
\forall a_1,\ldots,a_n\in\mathbb{R}^l,
\forall\alpha>0\qquad
\Biggl|\sum_{i=1}^n a_i\Biggr|^{\alpha}\le n^{(\alpha-1)_+}\sum_{i=1}^n |a_i|^{\alpha},
\end{equation}
used with $\alpha=2p$, we obtain
\[
\frac{\mathbb{E}\{|\sum_{i=1}^{n+1} R_i|^{{2p}}\}}{(n+1)!\mathbb{E}\{|R_1|^{2p}\}}
\le \frac{\mathbb{E}\{(n+1)^{({2p}-1)_+}\sum_{i=1}^{n+1}|R_i|^{2p}\}}
{(n+1)!\mathbb{E}\{|R_1|^{2p}\}}
= \frac{(n+1)^{({2p}-1)_+}}{n!}.
\]
It follows that $F_{\lambda_h}$ is an analytic function on $\mathbb{R}$
such that $F_{\lambda_h}(0)=1$. Therefore,
\[
F_{\lambda_h} (t)=1+t\psi_h(t)\quad\textnormal{with
}|\psi_h(t)|\le C(p,h,{\lambda_h})\quad\forall t\in
[0,T].
\]
Since $\mathbb{E}\{|R_1|^{2p}\}=\frac{1}{{\lambda_h}}
\int_{\{|y|>h\} }|y|^{2p}\pi(dy)$, the first equality is obvious.
{\smallskipamount=0pt
\begin{longlist}[(iii)]
\item[(ii)] If $\int_{|y|\le{h}}|y|^{2q}\pi(dy)<+\infty$ with $q\le1/2$,
then $Y^{h}$ has locally bounded variations and
${Y}^{h}_t+t\int_{|y|\le{h}} y\pi(dy)=\sum_{0<s\le t}
\Delta Y^{h}_t$. Inequality \eqref{eqelem2} with $\alpha=2q$ and the
compensation formula yield
\begin{eqnarray*}
\mathbb{E}\biggl\{\biggl|{Y}^{h}_t+t\int_{|y|\le{h}} y\pi(dy)\biggr|^{2q}\biggr\}
&\le & \mathbb{E}\Biggl\{\sum_{0<s\le t} |\Delta Y^{h}_t |^{{2q}}\Biggr\}
\\
&=& t\int_{y\le{h}}|y|^{2q} \pi(dy).
\end{eqnarray*}
Now, let $q\in(1/2,1]$. As $Y^{h}$ is a martingale, we derive
from the Burkholder--Davis--Gundy (BDG) inequality (see
\cite{bertoin}) that
\[
\mathbb{E}\{|Y^{h}_t|^{2q}\}\le C_q \mathbb{E}
\Biggl\{\Biggl(\sum_{0<s\le t}|\Delta Y^{h}_s|^{2}\Biggr)^{q}\Biggr\}.
\]
The second inequality follows from inequality \eqref{eqelem2} with
$\alpha=q$ and from the compensation formula.
\item[(iii)] One first considers case $p=1$. The
process $(M_t)$ defined by $M_t=|\hat{Z}_t|^2-t\int|y|^2\pi(dy)$
is a martingale. Then, in particular,
$\mathbb{E}\{|\hat{Z}_t|^2\}=t\int|y|^2\pi(dy)$. Suppose now
that $p>1$. In order to simplify the notation, we assume that
$T<1$. The BDG inequality yields
%
\begin{equation}\label{contpetit2}
\mathbb{E}\{|Y^{h}_t|^{2p}\}\le C_p \mathbb{E}
\Biggl\{\Biggl(\sum_{0<s\le t}|\Delta Y^{h}_s|^{2}\Biggr)^{p}\Biggr\}.
\end{equation}
\end{longlist}}%
For every integer $k\ge1$,
$M_{t,k}:=\sum_{0<s\le t}|\Delta Y^{h}_s)|^{2^k}-t\int_{\{|y|\le
h\}}|y|^{2^k}\pi(dy)$ is a martingale. By inequality
\eqref{eqelem2} and the BDG inequality applied to $(M_{t,k})$, we
obtain
\begin{eqnarray*}
&& \mathbb{E}\Biggl\{\Biggl(\sum_{0<s\le t}|\Delta Y^{h}_s|^{2^k}
\Biggr)^{{p}/{2^{k-1}}}\Biggr\}
\\
&&\qquad  \le C \biggl(\mathbb{E}\bigl\{|M_{t,k}|^{{p}/{2^{k-1}}}\bigr\}
+ \biggl(t\int_{\{|y|\le h\}}|y|^{2^k}\pi(dy) \biggr)^{{p}/{2^{k-1}}} \biggr)
\\
&&\qquad  \le  C \Biggl(\mathbb{E}\Biggl\{\Biggl(\sum_{0<s\le t}|\Delta
Y^{h}_s|^{2^{k+1}}\Biggr)^{{p}/{2^{k}}}\Biggr\}
+ \biggl(t\int_{\{|y|\le h\}}|y|^{2^k}\pi(dy) \biggr)^{{p}/{2^{k-1}}} \Biggr).
\end{eqnarray*}
Set $k_0=\inf\{k\ge1, 2^{k}\ge p\}$. Iterating the preceding
relation yields
\begin{eqnarray*}
\mathbb{E}\Biggl\{\Biggl(\sum_{0<s\le t}|\Delta Y^{h}_s|^{2}\Biggr)^{p}\Biggr\}
&\le & C\mathbb{E}\Biggl\{\Biggl(\sum_{0<s\le t}|\Delta
Y^{h}_s|^{2^{k_0+1}}\Biggr)^{{p}/{2^{k_0}}}\Biggr\}
\\
&&{} + C\sum_{k=1}^{k_0} \biggl(t\int_{\{|y|\le h\}}
|y|^{2^k}\pi(dy) \biggr)^{{p}/{2^{k-1}}}.
\end{eqnarray*}
By
construction, $p/2^{k_0}\le 1$. We then derive from inequality
\eqref{eqelem2} with $\alpha=p/2^{k_0}$, from the compensation
formula and from \eqref{contpetit2} that
%
\begin{eqnarray}\label{contpetit}
\mathbb{E}\{|Y_t^h|^{2p}\} &\le& C_p \mathbb{E}
\Biggl\{\Biggl(\sum_{0<s\le t}|\Delta Y^{h}_s|^{2}\Biggr)^{p}\Biggr\}
\nonumber\\[-8pt]
\\[-8pt]
&\le & C_pt\int_{\{|y|\le h\}}|y|^{2p}\pi(dy)+C_{p,h} t^{\eta_1}
\nonumber
\end{eqnarray}
with
$\eta_1=p/2^{k_0-1}>1$. We now consider $(\hat{Z}_t)$. For every
$h\in(0,+\infty)$, we have $\hat{Z}_t=Y_t^h+\hat{N}_t^h$, where
$\hat{N}_t^h=N_t^h-t\int_{|y|>h}y\pi(dy)$. Using the elementary
inequality,
%
\begin{equation}\label{elem}
\forall u,v \in\mathbb{R}_+, \forall \alpha\ge 1 \qquad
(u+v)^\alpha\le u^\alpha+\alpha 2^{\alpha-1}(u^{\alpha-1}v+v^\alpha),
\end{equation}
we derive from \eqref{bigjumprev} that
%
\begin{equation}\label{nth3}
\forall\alpha>1\qquad
\mathbb{E}\{|\hat{N}_t^h|^{\alpha}\}\le
t\int_{|y|>h}|y|^{\alpha}\pi(dy)+C_{\alpha,h} t^{\alpha\wedge2}.
\end{equation}
Using \eqref{elem} and the independence between $(\hat{N}_t^h)$
and $(Y_t^h)$ also yields
\begin{eqnarray*}
\mathbb{E}\{|\hat{Z}_t|^{2p}\} \le \mathbb{E}\{|\hat{N}_t^h|^{2p}\}
+C (\mathbb{E}\{|\hat{N}_t^h|\}^{2p-1}
\mathbb{E}\{|Y_t^{h}|\}+\mathbb{E}\{|Y_t^{h}|^{2p}\} ).
\end{eqnarray*}
Since $\mathbb{E}\{|Y_t^h|^2\}=t\int_{\{|y|\le h\}}|y|^2\pi(dy)$, we
derive from the Jensen inequality that $\mathbb{E}\{|Y_t^h|\}\le
C_h\sqrt{t}$. Hence, by \eqref{contpetit} and \eqref{nth3}, it
follows that, for every $h>0$ and $t\le T$,
\[
\mathbb{E}\{|\hat{Z}_t|^{2p}\}\le t\int_{\{|y|>
h\}}|y|^{2p}\pi(dy)+C^1_{p,h}t^{{3}/{2}\wedge\eta_1}+C^2_p
t\int_{\{|y|\le h\}}|y|^{2p}\pi(dy)
\]
with $\eta_1>1$,
$C_{p,h}^1>0$ and $C_p^2>0$. Let $\varepsilon$ be a positive
number. As $C_p^2$ does not depend on $h$, and $\int_{|y|\le
h}|y|^{2p}\pi(dy)\rightarrow0$ when ${h\rightarrow0}$, we can
choose $h_\varepsilon>0$ such that $C^2_p\int_{|y|\le
h_\varepsilon}|y|^{2p}\pi(dy)\le\varepsilon$. That yields the
announced inequality.
\end{pf*}
\begin{lemme} \label{lemmetech} 
Let $V$ be an EQ-function defined on $\mathbb{R}^d$. Then:
\begin{longlist}[(b)]
\item[(a)] If $p\in[0,1/2]$, $V^p$ is $\alpha$-H\"older for any
$\alpha\in[2p,1]$ and if $ p\in(0,1]$,
$ \nabla(V^p)$ is $\alpha$-H\"older for any $\alpha\in[2p-1,1]\cap (0,1]$.
\item[(b)] Let $x,$ $y\in\mathbb{R}^d$ and $\xi\in[x,x+y]$ and set
$\underline{v}=\min\{V(x), x\in\mathbb{R}^d\}.$ If $p\le1$,
%
\begin{equation}\label{tech}
\tfrac{1}{2}D^2(V^p)(\xi)y^{\otimes2}\le p\underline{v}^{p-1}\lambda_p|y|^2.
\end{equation}
If, moreover,
$|y|\le(1-\varepsilon)\frac{\sqrt{V}(x)}{[\sqrt{V}]_1}$ with
$\varepsilon\in(0,1]$, then,
%
\begin{equation}\label{INEQ2}
\tfrac{1}{2}D^2(V^p)(\xi)y^{\otimes2}\le
p\lambda_p\varepsilon^{2(p-1)} V^{p-1}(x)|y|^2.
\end{equation}
If $p> 1$,
%
\begin{equation}\label{INEQ3}
\qquad
\tfrac{1}{2} D^2(V^p)(\xi)y^{\otimes2}
\le p\lambda_p2^{(2(p-1)-1)_+}
\bigl(V^{p-1}(x)+\bigl[\sqrt {V}\bigr]_1 |y|^{2(p-1)} \bigr)|y|^2.
\end{equation}
\end{longlist}
\end{lemme}
\begin{pf}
Consider a continuous function $f\dvtx \mathbb{R}^d\mapsto \mathbb{R}$.
Let $\alpha\in(0,1]$ such that $|f|^{{1}/{\alpha}}$ is
Lipschitz. Then, $f$ is an $\alpha$-H\"older function. This
argument yields (a) (see~\cite{bib26} for details). Now,
let us pass to (b). We have
%
\begin{equation}
D^2(V^p)=pV^{p-1} \biggl(D^2V+(p-1)
\frac{\nabla{V}\otimes\nabla{V}}{V} \biggr),
\end{equation}
where $(\nabla{V}\otimes\nabla{V})_{i,j}=(\nabla
V)_i(\nabla V)_j$. Since $V^{p-1}\le\underline{v}^{p-1}$ if $p\le 1$,
we derive~\eqref{tech} from relations \eqref{sym} and
\eqref{cplamb}. For \eqref{INEQ2}, we consider $\xi=x+\theta y$
with $\theta\in[0,1]$ and $|y|\le (1-\varepsilon)\frac{\sqrt{V}(x)}{[\sqrt{V}]_1}$.
As $\sqrt{V}$ is a Lipschitz function,
\[
\sqrt{V}(\xi)\ge\sqrt{V}(x)-\bigl[\sqrt{V}\bigr]_1|y|
\ge \varepsilon\sqrt{V}(x) \quad \Longrightarrow\quad
V^{p-1}(\xi) \le \varepsilon^{2(p-1)}V^{p-1}(x).
\]
Hence, inequality \eqref{INEQ2} follows from \eqref{cplamb}. If
$p>1$,
\[
\sqrt{V}(\xi)\le\sqrt{V}(x)+\bigl[\sqrt{V}\bigr]_1|y|
\quad\Longrightarrow \quad
V^{p-1}(\xi)\le\bigl(\sqrt{V}(x)+\bigl[\sqrt{V}\bigr]_1|y|\bigr)^{2(p-1)}.
\]
We then derive \eqref{INEQ3} from \eqref{eqelem2}
[with $\alpha=2(p-1)$ and $n=2$] and from \eqref{cplamb}.
\end{pf}
\begin{lemme}\label{controljumpcomp}
Let $p\in(0,1)$, $q\in[0,1]$ and $a\in(0,1]$. Assume
\hyperlink{eH}{$\mathrm{(H^1_p)}$}, \hyperlink{eH}{$\mathrm{(H_q^2)}$} and
\ref{asS}.1. Then,
for every $\varepsilon>0$, there exists
$h_\varepsilon\in[0,+\infty]$, $T_\varepsilon>0$ and
$C_\varepsilon>0$ such that for every $x,z\in\mathbb{R}^d$, for every
$t\le T_\varepsilon$,
%
\begin{eqnarray}\label{eopz}
&& \mathbb{E}\bigl\{V^p\bigl(z+\kappa(x)Z_t^{h_\varepsilon}\bigr)
- V^p(z)\bigr\}
\nonumber\\[-8pt]
\\[-8pt]
&&\qquad \le t \biggl(p c_p\int|y|^{2p}\pi(dy)
1_{\{q\le p\}}\|\kappa(x)\|^{2p} +\varepsilon V^{p+a-1}(x)+C_\varepsilon
\biggr),
\nonumber
\end{eqnarray}
with $c_p$ given by \eqref{cplamb}, $h_\varepsilon\in(0,1]$ if
$p\le1/2<q$, $h_\varepsilon=0$ if $p,q\le1/2$ and
$h_\varepsilon=+\infty$ if $p\in(1/2,1)$.
\end{lemme}
\begin{pf}
Set $\Delta(z,x,U)=V^{p}(z+\kappa(x)U)-V^p(z)$. We first consider
the case $p\le1/2$ and $q>1/2$. Let $h\in(0,\infty)$. Since
$Z_t^h=Y_t^h+N_t^h$, we can decompose $\Delta(z,x,Z_t^h)$ as
follows:
\[
\Delta(z,x,Z_t^h)=\Delta\bigl(z+\kappa(x)N_t^h,x,Y_t^h\bigr)+\Delta(z,x,N_t^h).
\]
One controls each term of the right-hand side. On the one hand, as
$V^p$ is $2p$-H\"older with constant $[V^p]_{2p}=pc_p$ [see
\eqref{cplamb}], we deduce from Lemma \ref{moment}(i) that
%
\begin{eqnarray}\label{nth1}
\mathbb{E}\{\Delta(z,x,N_t^h)\}
&\le & pc_p\|\kappa(x)\|^{2p}\mathbb{E}\{|N_t^h|^{2p}\}
\nonumber\\[-8pt]
\\[-8pt]
&\le & pc_p\int_{|y|>h}|y|^{2p}\pi(dy)\|\kappa(x)\|^{2p}\bigl(t+\psi_h(t)
t^2\bigr),
\nonumber
\end{eqnarray}
where $\psi_h$ is a locally bounded function. On the other hand,
we set $\tilde{z}=z+\kappa(x)N_t^h$. By the Taylor formula,
\[
\Delta(\tilde{z},x,Y_t^h)
=\langle\nabla (V^p)(\tilde{z}),\kappa(x)Y_t^h\rangle
+ \langle\nabla(V^p)(\xi )-\nabla (V^p)(\tilde{z}),\kappa(x)Y_t^h\rangle
\]
with
$\xi\in[\tilde{z},\tilde{z}+\kappa(x)Y_t^h]$. As $(N_t^h)$ and
$(Y_t^h)$ are independent and $Y_t^h$ is centered,
$\mathbb{E}\{ \langle \nabla(V^p)(\tilde{z}),\kappa(x)Y_t^h\rangle\}=0$. By
Lemma \ref{lemmetech}, $V^{p-1}\nabla V=\nabla(V^p)/p$ is
\mbox{$(2q-1)$-H\"older} (because $2q-1\in[2p-1,1]\cap(0,1]$ in this case).
Then, it follows from Lemma~\ref{moment}(ii).2 that
%
\begin{eqnarray}\label{yth2}
\mathbb{E}\bigl\{\Delta\bigl(z+\kappa(x)N_t^h,x,Y_t^h\bigr)\bigr\}
&\le & p[V^{p-1}\nabla V]_{2q-1}\|\kappa(x)\|^{2q}
\mathbb{E}\{|Y_t^h|^{2q}\}
\nonumber\\[-8pt]
\\[-8pt]
&\le& C\|\kappa(x)\|^{2q} t\int_{|y|\le h}|y|^{2q}\pi(dy).
\nonumber
\end{eqnarray}
Let $\varepsilon>0$. First, by \ref{asS}.1,
$\|\kappa (x)\|^{2q}\le CV^{p+a-1}$. Then, using that\break
$\int_{|y|\le h}|y|^{2q}\pi(dy)\rightarrow0$ when $h\rightarrow0$, we can fix
$h_\varepsilon\in(0,1]$ such that
%
\begin{eqnarray}\label{yth1}
\mathbb{E}\bigl\{\Delta\bigl(z+\kappa(x)N_t^{h_\varepsilon
},x,Y_t^{h_\varepsilon}\bigr)\bigr\}
\le \frac{\varepsilon}{2} tV^{p+a-1}(x).
\end{eqnarray}
Second, since $\psi_{h_\varepsilon}$ is locally bounded, it
follows from \eqref{nth1} that there exists $C_\varepsilon^1$ such
that, for every $t\le1$,
$\mathbb{E}\{\Delta(z,x,N_t^{h_\varepsilon})\}\le C_\varepsilon^1
t\|\kappa(x)\|^{2p}$. Now, as $p<q$, for every $\delta>0$, there
exists $C_\delta^2>0$ such that $\|\kappa(x)\|^{2p}\le\delta
V^{p+a-1}+C_\delta^2$ [see \eqref{closeargumentrev} for similar
arguments]. Hence, setting
$\delta_\varepsilon=\varepsilon/(2C^1_\varepsilon)$ yields
%
\begin{equation}\label{nth2}
\mathbb{E}\{\Delta(z,x,N_t^{h_\varepsilon})\}\le
t\biggl (\frac{\varepsilon}{2} V^{p+a-1}(x)+C_\varepsilon \biggr)
\end{equation}
with $C_\varepsilon=C_\varepsilon^1 C_{\delta_\varepsilon}^2$.
Then, adding up \eqref{yth1} and \eqref{nth2} yields the result when $p\le\break 1/2<q$.

When $p,q\le1/2$, we deal with
$(\check{Z}_t)=({Z}^0_t)$. For every $h>0$,
$\check{Z}_t=\check{Y}_t^h+N_t^h$, where
$\check{Y}_t^h=Y_t^h+t\int_{\{|y|\le h\}}y\pi(dy)$. Hence, for
every $h>0$,
\[
\Delta(z,x,\check{Z}_t)=\Delta\bigl(z+\kappa(x)N_t^h,x,\check
{Y}_t^h\bigr)+\Delta(z,x,N_t^h).
\]
If $q\le p$, $\pi$ satisfies \hyperlink{eH}{$\mathrm{(H^2_p)}$}. Since $p\le1/2$,
$V^p$ is $2p$-H\"older. Therefore, by Lemma~\ref{moment}(ii).1,
\[
\mathbb{E}\bigl\{\Delta\bigl(z+\kappa(x)N_t^h,x,\check{Y}_t^h\bigr)\bigr\}
\le pc_p t\|\kappa(x)\|^{2p}\int_{|y|\le h}|y|^{2p}\pi(dy).
\]
By summing up
this inequality and \eqref{nth1}, we deduce $\eqref{eopz}$. When
$p< q\le1/2$, we use that $V^p$ is $2q$-H\"older (see Lemma
\ref{lemmetech})
and a proof analogous to the case $p\le1/2<q$ yields the result.

Finally, we consider the case $p>1/2$ where we deal with
$\hat{Z}_t=Z_t^\infty$. For every $h>0$, we have
$\hat{Z}_t=Y_t^h+\hat{N}_t^h$, where
$\hat{N}_t^h=N_t^h-\int_{\{|y> h\}}y\pi(dy)$. For every
$h>0$, $\Delta(z,x,\hat{Z}_t)$ can be written as follows:
%
\begin{equation}\label{whenpsup1rev}
\Delta(z,x,\hat{Z}_t)=\Delta\bigl(z+\kappa(x)\hat
{N}_t^h,x,{Y}_t^h\bigr)+\Delta(z,x,\hat{N}_t^h).
\end{equation}
One the one
hand, by the same process as that used for \eqref{yth2} and by
inequality \eqref{nth3}, we have
\begin{eqnarray*}
\mathbb{E}\{\Delta(z,x,\hat{N}_t^h)\}
&\le& p[V^{p-1}\nabla V]_{2p-1}\|\kappa(x)\|^{2p}\mathbb{E}\{|\hat{N}_t^h|^{2p}\}
\nonumber
\\
&\le & t\|\kappa(x)\|^{2p}
\biggl(\int_{|y|>h}|y|^{2p}\pi(dy)+C_h t^{2p-1} \biggr).
\end{eqnarray*}
On the other hand, by Lemma \ref{lemmetech}, $V^{p-1}\nabla V$ is
$(2(p\vee q)-1)$-H\"older. Hence, \eqref{yth2} is still valid in
this case if we replace $q$ with $p\vee q$. By using this control
for the first term of the right-hand side of \eqref{whenpsup1rev},
we obtain if $q\le p$
\begin{eqnarray*}
&& \mathbb{E}\{\Delta(z,x,\hat{Z}_t)\}
\\
&&\qquad
\le t\|\kappa(x)\|^{2p}  \biggl(\int_{\{|y|>h\}}|y|^{2p}\pi(dy)
+\int_{|y|\le h}|y|^{2p}\pi(dy)+C_h t^{2p-1} \biggr).
\end{eqnarray*}
The result follows in this
case. When $q>p$, the sequel of the proof is similar to the case
$p\le1/2<q$.
\end{pf}

\subsection[Proof of Proposition 2]{Proof of Proposition \textup{\protect\ref{prop9}}}\label{proofprop9}

For this proof, one needs to study separately
the $p<1$
and $p\ge1$ cases. We detail the first case. When $p\ge1$, we
briefly indicate the
process of the proof which is close to that of Lemma 3 of
\cite{bib3}.

\textbf{Case} $\bolds{p<1.}$ For $h>1$, we set
$\bar{Z}_{n}^h=\bar{Z}_{n}-\gamma_n\int_{\{1<|y|\le h\}}y\pi(dy)$
and for $h\in(0,1)$,
$\bar{Z}_{n}^h=\bar{Z}_{n}+\gamma_n\int_{\{h<|y|\le1\}}y\pi(dy)$.
If $q\le1/2$ (resp. $p>1/2$), we can take $h=0$ (resp.
$h=+\infty$). Thus, we can write
%
\begin{eqnarray}\label{955}
\qquad
\Delta \bar{X}_{n+1} :\!\!&=&\bar{X}_{n+1}-\bar{X}_{n}
= \sum_{k=1}^3\Delta \bar{X}^h_{n+1,k}\qquad \mbox{with }
\Delta\bar{X}^h_{n+1,1}=\gamma_{n+1}b^h(\bar{X}_n),\hspace*{-13pt}
\nonumber\\[-8pt]\hspace*{-13pt}
\\[-8pt]
\Delta \bar{X}^h_{n+1,2} &=& \sqrt{\gamma_{n+1}}\sigma(\bar{X}_{n})U_{n+1}
\quad\mbox{and}\quad
\Delta\bar{X}^h_{n+1,3}=\kappa(\bar{X}_n)\bar{Z}^h_{n+1}.\hspace*{-13pt}
\nonumber
\end{eqnarray}
The idea is to study the difference
$V^p(\bar{X}_{n+1})-V^p(\bar{X}_n)$ as the sum of three terms that
correspond to the above decomposition. For $k=1,2,3$, set
$\bar{X}^h_{n+1,k}=\bar{X}^h_n+\sum_{i=1}^k\Delta \bar{X}^h_{n+1,i}$.
\begin{longlist}
\item First term: There exists $n_1\in\mathbb{N}$ such that,
for every $n\ge n_1$
%
\begin{eqnarray}\label{un}
&& \mathbb{E}\{V^p(\bar{X}^h_{n+1,1})-V^p(\bar{X}_{n})|\mathcal{F}_n\}
\nonumber\\[-8pt]
\\[-8pt]
&&\qquad \le p\gamma_{n+1}
\frac{\langle\nabla V,b^h\rangle}{V^{1-p}}(\bar{X}_n)+C\gamma_{n+1}^2
V^{a+p-1}(\bar{X}_n).
\nonumber
\end{eqnarray}
Indeed, from Taylor's formula,
\begin{eqnarray*}
V^p(\bar{X}^h_{n+1,1})-V^p(\bar{X}_{n})
= p\gamma_{n+1}\frac{\langle\nabla V,b^h\rangle}{V^{1-p}}(\bar{X}_n)
+ \frac{1}{2}D^2 (V^p)(\xi^1_{n+1})(\Delta
\bar{X}^h_{n+1,1})^{\otimes2},
\end{eqnarray*}
where
$\xi^1_{n+1}\in [\bar{X}_{n}, \bar{X}_{n}+\gamma_{n+1}
b^h(\bar{X}_n)]$. Set $x=\bar{X}_n$ and $y=\gamma_{n+1}b^h(\bar{X}_n)$.
Since $\gamma_n\stackrel{n\rightarrow+\infty}{\longrightarrow}0$ and $|b^h|\le
C\sqrt{V}$ by \ref{asS}.1, there exists $n_1\in\mathbb{N}$ such
that, for \mbox{$n\ge n_1$}, $|y| \le\frac{\sqrt{V}(x)}{2[\sqrt{V}]_1}$ a.s.
Thus, we can apply the second inequality of Lemma~\ref{lemmetech}(b) with
$\varepsilon=1/2$ and deduce \eqref{un} from \ref{asS}.1.
\item Second term: For every $\varepsilon>0$, there exists
$n_{2,\varepsilon}\in\mathbb{N}$ such that, for every $n\ge
n_{2,\varepsilon}$,
%
\begin{equation}\label{deux}
\mathbb{E}\{V^p(\bar{X}^h_{n+1,2})-V^p(\bar{X}^h_{n+1,1})|\mathcal{F}_n\}
\le \varepsilon\gamma_{n+1}
V^{a+p-1}(\bar{X}_n)+C^1_\varepsilon\gamma_{n+1}.
\end{equation}
Let us
prove this inequality. Since $\mathbb{E}\{U_{n+1}|\mathcal{F}_n\}=0$, we
deduce from Taylor's formula that
\[
\mathbb{E}\{V^p(\bar{X}^h_{n+1,2})-V^p(\bar{X}^h_{n+1,1})|\mathcal{F}_n\}
= \tfrac{1}{2}\mathbb{E}\{D^2
(V^p)(\xi^2_{n+1})(\Delta\bar{X}_{n+1,2})^{\otimes2}|\mathcal{F}_n\}
\]
with $\xi^2_{n+1}\in [\bar{X}^h_{n+1,1};\bar{X}^h_{n+1,2}]$. Set
$x=\bar{X}^h_{n+1,1}$ and $y=\sqrt{\gamma_{n+1}}\sigma(x)U_{n+1}$.
By \ref{asS}.1, $\|\sigma(x)\|\le C_\sigma\sqrt{V}(x)$ because
$p+a-1\le1$. Then, the conditions of~\eqref{INEQ2} are satisfied
with $\varepsilon=1/2$ if $|U_{n+1}|\le\rho_{n+1}=1/(2C_\sigma
[\sqrt{V}]_1 \sqrt{\gamma_{n+1}})$. Therefore,
\begin{eqnarray*}
&& \mathbb{E}\bigl\{D^2 (V^p)(\xi^2_{n+1})(\Delta
\bar{X}_{n+1,2})^{\otimes2}1_{\{|U_{n+1}|\le
\rho_{n+1}\}}|\mathcal{F}_n\bigr\}
\\
&&\qquad \le C \gamma_{n+1}V^{p-1}(\bar{X}_n)
\operatorname{Tr}(\sigma\sigma^*)(\bar{X}_n)
\\
&&\qquad \le C\gamma_{n+1}V^{a+2(p-1)}(\bar{X}_n)
\end{eqnarray*}
since $\operatorname{Tr}(\sigma\sigma^*)\le CV^{p+a-1}$ when $p<1$. By
\eqref{tech} and \ref{asS}.1, we also have
\[
\mathbb{E}\bigl\{D^2 (V^p)(\xi^2_{n+1})(\Delta\bar{X}_{n+1,2})^{\otimes
2}1_{\{|U_{n+1}|>\rho_{n+1}\}}|\mathcal{F}_n\bigr\}
\le C\delta_{n+1}\gamma_{n+1} V^{p+a-1}(\bar{X}_n),
\]
where
$\delta_n=\mathbb{E}\{|U_{n}|^2 1_{\{|U_{n}|>\rho_{n}\}}\}$. Now, let
$\varepsilon>0$. First, since $a+2(p-1)<a+p-1$ when $p<1$, there
exists $C_\varepsilon>0$ such that $V^{p-1}\operatorname{Tr}(\sigma\sigma^*)\le
\varepsilon V^{a+p-1}+C_\varepsilon$ [see~\eqref{closeargumentrev}
for similar arguments]. Second, since
$\rho_n\rightarrow+\infty$, $\delta_n\rightarrow0$. Thus, there
exists $n_{2,\varepsilon}\in\mathbb{N}$ such that, for every $n\ge
n_{2,\varepsilon}$, $C \delta_n V^{p+a-1}\le\varepsilon V^{p+a-1}$.
The combination of these two arguments yields \eqref{deux}.
\item Third term: For every $\varepsilon>0$, there
exists $h_\varepsilon\in[0,\infty]$, $C^2_\varepsilon>0$ and
$n_{3,\varepsilon}$ such that, for all $n\ge n_{3,\varepsilon}$,
%
\begin{eqnarray}\label{troi3}
\qquad
&& \mathbb{E}\{V^p(\bar{X}^{h_\varepsilon}_{n+1,3})
-V^p(\bar{X}^{h_\varepsilon}_{n+1,2})|\mathcal{F}_n\}
\nonumber\\[-8pt]
\\[-8pt]
&&\qquad \le\gamma_{n+1}
\biggl (p c_p\int|y|^{2p}\pi(dy)1_{\{q\le p\}}\|\kappa(\bar{X}_n)\|^{2p}
+ \varepsilon V^{p+a-1}(\bar{X}_n)+C^2_\varepsilon \biggr)
\nonumber
\end{eqnarray}
with $h_\varepsilon\in(0,1]$ if $p\le1/2<q$, $h_\varepsilon=0$ if
$p,q\le1/2$ and $h_\varepsilon=+\infty$ if $p\in(1/2,1)$. This
step is a consequence of Lemma \ref{controljumpcomp}: since
$U_{n+1}$ and $\bar{Z}_{n+1}^h$ are independent, we have
\[
\mathbb{E}\{V^p(\bar{X}^{h_\varepsilon}_{n+1,3})-V^p(\bar
{X}^{h_\varepsilon}_{n+1,2})|\mathcal{F}_n\}
=\mathbb{E}\{G_{h_\varepsilon}(\bar{X}^{h_\varepsilon}_{n+1,2},\bar
{X}_n)|\mathcal{F}_n\},
\]
where $G_h(z,x)=\mathbb{E}\{V^p(z+\kappa(x)Z_t^{h})-V^p(z)\}$.
Then, Lemma \ref{controljumpcomp} yields \eqref{troi3}.
\end{longlist}

We can now prove the proposition. Let $\varepsilon>0$.
By adding \eqref{un}, \eqref{deux} and \eqref{troi3} and using
that $\gamma_{n}^2\le\varepsilon\gamma_n$ for sufficiently large
$n$ (since $\gamma_n\rightarrow0$), we obtain that there exists
$n_{\varepsilon}\in\mathbb{N}$, $h_\varepsilon>0$ and
$C_\varepsilon>0$
such that, for every $n\ge n_{\varepsilon}$,
%
\begin{eqnarray}\label{reconize}
&& \mathbb{E}\{V^p(\bar{X}_{n+1})|\mathcal{F}_n\}
\nonumber\\
&&\qquad \le  V^p(\bar{X}_n)+ \gamma_{n+1}
\bigl(\varepsilon V^{p+a-1}(\bar{X}_n)+ C_\varepsilon\bigr)
\nonumber\\[-8pt]
\\[-8pt]
&&\qquad \quad {} +\gamma_{n+1}pV^{p-1}(\bar{X}_n)
\nonumber\\
&&\qquad \quad\hspace*{10pt} {}\times
\biggl(\langle\nabla V,b^{h_\varepsilon}\rangle+1_{q\le
p}c_p\int|y|^{2p}\pi(dy)\|\kappa\|^{2p}V^{1-p} \biggr)
(\bar {X}_n).
\nonumber
\end{eqnarray}
When $p,q\le1/2$ (resp. $p>1/2$), $h_\varepsilon=0$ (resp.
$h_\varepsilon=\infty)$. We deduce that $\langle\nabla
V,b^{h_\varepsilon}\rangle= \langle\nabla V,\tilde{b}\rangle$ because
$\tilde{b}=b^0$ (resp. $\tilde{b}=b^\infty)$ when $p,q\le1/2$\break
(resp. $p>1/2$). We then recognize the left-hand side of
\ref{asS}.2 in \eqref{reconize}. When $p\le1/2<q$,
$h_\varepsilon\in(0,\infty)$ and $\langle\nabla
V,b^{h_\varepsilon}\rangle=\langle\nabla
V,\tilde{b}\rangle+\Phi_{h_\varepsilon}$, where $\Phi_h(x)=\langle
\nabla V(x),\kappa(x) \int_{\{h_\varepsilon<|y|\le1\}}y\pi(dy)\rangle$.
Therefore, by \ref{asS}.2, we obtain
\begin{eqnarray*}
&& \mathbb{E}\{V^p(\bar{X}_{n+1})|\mathcal{F}_n\}
\\
&&\qquad \le  V^p(\bar{X}_n)+ \gamma_{n+1}pV^{p-1}(\bar{X}_n)
\bigl(\beta-\alpha V^{a}(\bar{X}_n) \bigr)
\\
&&\qquad \quad {}+ \gamma_{n+1}\bigl(\varepsilon V^{p+a-1}(\bar{X}_n)+C_\varepsilon+1_{\{
p\le {1}/{2}<q\}}pV^{p-1}(\bar{X}_n)\Phi_{h_\varepsilon}(\bar{X}_n)\bigr).
\end{eqnarray*}
When $p,q\le1/2$ or $p>1/2$, we set $\varepsilon=p\alpha/2$ and
obtain \hyperlink{e15}{$\mathrm{(R_{a,p})}$} with
$\beta'=p\beta+C_\varepsilon/\underline{v}^{p-1}$ and
$\alpha'=p\alpha/2$. When $p\le1/2<q$, by \ref{asS}.1,
one checks that, for every $\varepsilon>0$, there exists
$\tilde{C}_\varepsilon>0$ such that
$V^{p-1}|\Phi_{h_\varepsilon}|\le\varepsilon V^{p+a-1}+\tilde{C}_\varepsilon$
and the result follows.

\textbf{Case} $\bolds{p\ge1.}$ Thanks to Taylor's formula,
\begin{eqnarray*}
V^p(\bar{X}_{n+1})&=& V^p(\bar{X}_{n})+\gamma_{n+1}
\langle\nabla{(V^p)}(\bar{X}_n),\Delta X_{n+1}\rangle
\\
&&{} + \tfrac{1}{2}D^2(V^p)(\xi_{n+1})(\Delta\bar{X}_{n+1})^{\otimes2},
\end{eqnarray*}
where $\xi_{n+1}\in[\bar{X}_n,\bar{X}_{n+1}]$ and
\[
\Delta\bar{X}_{n+1}
= \bar{X}_{n+1}-\bar{X}_n
= \gamma_{n+1}b^\infty(\bar{X}_n)
+\sqrt{\gamma_{n+1}}\sigma(\bar{X}_n)U_{n+1}
+\kappa(\bar{X}_n)\bar{Z}_{n+1}^\infty
\]
with
$\bar{Z}_{n}^\infty=\bar{Z}_{n}-\gamma_n\int_{\{|y|>1\}}y\pi(dy)$.
Using that $\tilde{b}=b^\infty$ in this case and that
$\mathbb{E}\{U_{n+1}|\mathcal{F}_n\}=\mathbb{E}\{\bar{Z}_{n+1}^{\infty
}|\mathcal{F}_n\}=0$ yields
\begin{eqnarray*}
\mathbb{E}\{V^p(\bar{X}_{n+1})|\mathcal{F}_n\}
&=&  V^p(\bar{X}_n)+p\gamma_{n+1}
V^{p-1}(\bar{X}_n)\langle\nabla{V},\tilde{b}\rangle(\bar{X}_n)
\\
&&{} +\tfrac{1}{2}\mathbb{E}\{D^2(V^p)(\xi_{n+1})
(\Delta\bar{X}_{n+1})^{\otimes2}|\mathcal{F}_n\} .
\end{eqnarray*}
The sequel of
the proof consists in studying the last term of this equality.
The main tools for this are the last inequality of Lemma
\ref{lemmetech} which provides a control of
$D^2(V^p)(\xi_{n+1})(\Delta\bar{X}_{n+1})^{\otimes2}$ and Lemma
\ref{moment}(iii), which gives a control of the moments of the
jump component (see \cite{bib26} for details or \cite{bib3} for a
similar proof).

\subsection[Consequences of Proposition 2]{Consequences of Proposition \textup{\protect\ref{prop9}}}\label{conseqprop}

In Proposition \ref{prop9}, we
established \hyperlink{e15}{$\mathrm{(R_{a,p})}$}. According to Lemma
\ref{explication}, it suffices now to prove
\hyperlink{eCps}{$\mathrm{(C_{p,s})}$}.
This property is established in Corollary \ref{cor} and is a
consequence of Proposition \ref{prop9} [under additional
assumptions on $(\gamma_n)$ and $(\eta_n)$ when $s<2$]. More
precisely, we first show in Lemma \ref{lemme4} that a
supermartingale property can be derived from
\hyperlink{e15}{$\mathrm{(R_{a,p})}$}
and that this property provides an $L^{p+a-1}$-control of the
sequence $(V(\bar{X}_n))$ [see \eqref{lpcontrol}]. Second, we
show in Corollary \ref{cor} that we can derive
\hyperlink{eCps}{$\mathrm{(C_{p,s})}$} from this lemma.
\begin{lemme}\label{lemme4}
Let $a\in(0,1]$ and $p>0$. Assume
\hyperlink{eH}{$\mathrm{(H^1_p)}$} and
\hyperlink{e15}{$\mathrm{(R_{a,p})}$}.
Let $(\theta_n)_{n\in\mathbb{N}}$ be a nonincreasing sequence
of nonnegative numbers such that
$\sum_{n\ge1} \theta_n\gamma_n < \infty$. Then, there exists
$n_0\ge0$, $\hat{\alpha}>0$ and $\hat{\beta}>0$ such that
$(S_n)_{n\ge n_0}$ defined by
\[
S_n=\theta_nV^p(\bar{X}_n)+\hat{\alpha}
\sum_{k=1}^n\theta_k \gamma_kV^{p+a-1}(\bar{X}_{k-1})
+ \hat{\beta}\sum_{k>n}\theta_k\gamma_k
\]
is a nonnegative $L^1$-supermartingale. In particular,
%
\begin{equation}\label{lpcontrol}
\qquad
\sum_{n\ge1}\theta_n\gamma_n
\mathbb{E}\{V^{p+a-1}(\bar{X}_{n-1})\}<+\infty\quad\mbox{and}\quad
\mathbb{E}\{V^p(\bar{X}_n)\}
\stackrel{n\rightarrow+\infty}{=}O\biggl(\frac{1}{\theta_n}\biggr).
\end{equation}
\end{lemme}
\begin{pf}
Since $b$, $\sigma$ and $\kappa$ have sublinear growth and
$\bar{Z}_n\in L^{2p}$ for every $n\ge1$, we can check by
induction that, for every $n\ge0$, $V^p(\bar{X}_n)$ is integrable.
Denote by $(\Delta_n)_{n\ge1}$ the sequence of martingale
increments defined by
$\Delta _n=V^p(\bar{X}_{n})-\mathbb{E}\{V^p(\bar{X}_{n})|\mathcal{F}_{n-1}\}$.
By \hyperlink{e15}{$\mathrm{(R_{a,p})}$},
there exists $n_0\in\mathbb{N}$ such that, for every $n\ge n_0$,
\begin{eqnarray*}
&& \theta_{n+1} V^p(\bar{X}_{n+1})
\\
&&\qquad \le  \theta_{n+1}\Delta_{n+1}
+\theta_{n+1}\mathbb{E}\{V^p(\bar{X}_{n+1})|\mathcal{F}_n\}
\\
&&\qquad \le  \theta_{n+1} \Delta_{n+1}
+\theta_{n+1}
\bigl(V^p(\bar{X}_n)+\gamma_{n+1}V^{p-1}
(\bar{X}_n)\bigl(\beta'-\alpha'V^{a}(\bar{X}_n)\bigr) \bigr).
\end{eqnarray*}
By the same argument as in \eqref{closeargumentrev}, one can find
$\hat{\alpha}>0$ and $\hat{\beta}>0$ such that
$V^{p-1}({\beta}'-{\alpha}'V^a)\le
\hat{\beta}-\hat{\alpha}V^{p+a-1}$. Since $(\theta_n)$ is
nonincreasing, we deduce that
\begin{eqnarray*}
&& \theta_{n+1}
\bigl(V^p(\bar{X}_{n+1})+\hat{\alpha}\gamma_{n+1}
V^{p+a-1}(\bar{X}_n)\bigr)
\\
&&\qquad
\le\theta_{n}V^p(\bar{X}_{n})+\theta_{n+1}
\Delta_{n+1}+\theta_{n+1}\gamma_{n+1}\hat{\beta}.
\end{eqnarray*}
Adding ``$\hat{\alpha}\sum_{k=1}^{n}\theta_k\gamma_k
V^{p+a-1}(\bar{X}_
{k-1})+\hat{\beta}\sum_{k>n+1}\theta_k\gamma_k$'' to both sides of
the inequality yields
\begin{eqnarray*}
S_{n+1}\le S_n+\theta_{n+1}\Delta_{n+1}
\quad\Longrightarrow\quad
\mathbb{E}\{S_{n+1}|\mathcal{F}_n\}\le S_n \qquad \forall n\ge n_0.
\end{eqnarray*}
Since $S_{n_0}\in L^1$, it follows that $(S_n)_{n\ge n_0}$ is a
nonnegative supermartingale and then, that
$\sup\mathbb{E}\{S_n\}<+\infty$. The result is obvious.
\end{pf}
\begin{Corollaire}\label{cor}
Let $a\in(0,1]$, $p>0$ and $q\in[0,1]$. Assume
\hyperlink{eH}{$\mathrm{(H^1_p)}$},
\hyperlink{eH}{$\mathrm{(H^2_q)}$} and \ref{asS}.
If $\mathbb{E}\{|U_1|^{2(p\vee1)}\}<+\infty$ and
$(\eta_n/\gamma_n)_{n\in\mathbb{N}}$
is nonincreasing,
%
\begin{equation}\label{tuc0}
\sum_{n\ge 1}\biggl(\frac{\eta_n}{H_n\gamma_n}\biggr)^{2}
\mathbb{E}\bigl\{\bigl|V^{{p}/{2}}(\bar{X}_n)
- V^{{p}/{2}}\bigl(\bar{X}_{n-1}+\gamma_n
\tilde{b}(\bar{X}_{n-1})\bigr)\bigr|^{2}\bigr\} < +\infty.
\end{equation}
In
particular, \hyperlink{eCps}{$\mathrm{(C_{p,2})}$} holds with $\rho=2$ and
$\pi_n=V^{{p}/{2}}(\bar{X}_{n}+\gamma_n
\tilde{b}(\bar{X}_{n}))$.

Furthermore, if conditions \eqref{cond1} and \eqref{hop} are
satisfied for $s\in(1,2)$,
%
\begin{eqnarray}\label{tuc}
&& \sum_{n\ge 1}\biggl(\frac{\eta_n}{H_n\gamma_n}\biggr)^{f_{a,p}(s)}
\mathbb{E}\bigl\{\bigl|V^{{p}/{s}}(\bar{X}_n
)-V^{{p}/{s}}\bigl(\bar{X}_{n-1}+\gamma_n\tilde{b}(\bar{X}_{n-1})\bigr)
\bigr|^{f_{a,p}(s)}\bigr\}
\nonumber\\[-8pt]
\\[-8pt]
&&\qquad < +\infty.
\nonumber
\end{eqnarray}
In particular, \hyperlink{eCps}{$\mathrm{(C_{p,s})}$} holds with $\rho=f_{a,p}(s)$
and $\pi_n=V^{{p}/{s}}(\bar{X}_{n}+\gamma_n
\tilde{b}(\bar{X}_{n}))$.
\end{Corollaire}
\begin{pf}
Let us begin the proof by two useful remarks. First,
\eqref{tuc0} is a particular case of \eqref{tuc} since
$f_{a,p}(2)=2$ and \eqref{hop} is always satisfied in this case.
Indeed, as $(\eta_n/\gamma_n)_{n\in\mathbb{N}}$ is nonincreasing,
so is $(\frac{1}{\gamma_n}(\frac{\eta_n}{H_n\sqrt{\gamma_n}})^2)$,
and
%
\begin{equation}\label{eirop}
\qquad
\sum_{n\ge1}\biggl(\frac{\eta_n}{H_n\sqrt{\gamma_n}}\biggr)^2
\le \frac{\eta_1}{\gamma_1}\sum_{n\ge1}\frac{\eta_n}{H^2_n}
\le \frac{\eta_1}{\gamma_1}\sum_{n\ge1}\frac{\Delta H_n}{H_n^2}
\le C \biggl(1+\int_{\eta_1}^{\infty}\frac{dt}{t^2} \biggr)<\infty,
\end{equation}
with $\Delta H_n=H_n-H_{n-1}$. Then, it suffices to prove
\eqref{tuc}. Second, by Lemma \ref{lemme4} applied with
$\theta_n=\frac{1}{\gamma_n}(\frac{\eta_n}{H_n\sqrt{\gamma
_n}})^{f_{a,p}(s)}$,
we have
%
\begin{eqnarray}\label{dij}
\sum_{n\ge 1} \biggl(\frac{\eta_n}{H_n\sqrt{\gamma_n}}\biggr)^{f_{a,p}(s)}
\mathbb{E}\{ V^{p+a-1}(\bar{X}_{n-1})\}<+\infty.
\end{eqnarray}
Hence, one checks that \eqref{tuc} holds as soon as
%
\begin{eqnarray}\label{tuo}
&& \mathbb{E}\bigl\{\bigl|V^{{p}/{s}}(\bar{X}_n)
-V^{{p}/{s}}\bigl(\bar {X}_{n-1}+\gamma_n \tilde{b}(\bar{X}_{n-1})
\bigr)\bigr|^{f_{a,p}(s)}\bigr\}
\nonumber\\[-8pt]
\\[-8pt]
&&\qquad \le C\gamma_n^{{f_{a,p}(s)}/{2}}\mathbb{E}
\{V^{p+a-1}(\bar{X}_{n-1})\}.
\nonumber
\end{eqnarray}
Thus, we only need to prove \eqref{tuo}. We inspect the
$p/s\le1/2$ and $p/s>1/2$ cases successively.

\textbf{Case} $\bolds{p/s\le1/2.}$
In this case, $f_{a,p}(s)=s$. We keep the notation
introduced in~\eqref{955}, with $h=1$ if $p\le1/2<q$, $h=0$ if
$p,q\le1/2$ and $h=+\infty$ if $p>1/2$, and derive from
\eqref{eqelem2} that
%
\begin{eqnarray}\label{decomprev5}
\quad
&& \bigl|V^{{p}/{s}}(\bar{X}_n)-V^{{p}/{s}}\bigl(\bar{X}_{n-1}
+\gamma_n\tilde{b}(\bar{X}_{n-1})\bigr)\bigr|^{s}
\nonumber\\[-8pt]
\\[-8pt]
&&\qquad \le  C|V^{{p}/{s}}(\bar{X}_{n,2}^h)-V^{{p}/{s}}(\bar{X}^h_{n,1})|^{s}
+C|V^{{p}/{s}}(\bar{X}_{n})-V^{{p}/{s}}(\bar{X}_{n,2}^h)|^{s}.
\nonumber
\end{eqnarray}
We study successively the two right-hand side members. First, by the
Taylor formula,
\[
|V^{{p}/{s}}(\bar{X}_{n,2}^h)-V^{{p}/{s}}(\bar{X}_{n,1}^h)|^{s}
\le C\gamma_n^{{s}/{2}} |\langle V^{{p}/{s}-1}\nabla
V(\xi_{n}^1), \sigma(\bar{X}_{n-1})U_n\rangle |^s
\]
with
$\xi_{n}^1\in[ \bar{X}_{n,1}^h;\bar{X}_{n,2}^h]$. The function
$V^{{p}/{s}-1}\nabla V$ is bounded. Hence, since
$\|\sigma\|^s\le C \operatorname{Tr}(\sigma\sigma^*)$ (because $s\le2$),
we derive from \ref{asS}.1 that
%
\begin{equation}\label{corequation1rev}
\mathbb{E} \{|V^{{p}/{s}}(\bar{X}_{n,2}^h)
-V^{{p}/{s}}(\bar{X}_{n,1}^h)|^{s}|\mathcal{F}_{n-1} \}
\le C\gamma_{n}^{{s}/{2}}V^{a+p-1}(\bar{X}_{n-1}).
\end{equation}
Second, since $U_n$ and $\bar{Z}_n^h$ are independent,
\begin{eqnarray}
\mathbb{E} \{|V^{{p}/{s}}(\bar{X}_{n})-V^{{p}/{s}}(\bar{X}_{n,2}^h)|^{s}\}
&=& \mathbb{E}\{\Upsilon_h(\bar{X}^h_{n,2},
\bar{X}_{n-1},\gamma_n)|\mathcal{F}_{n-1} \}
\nonumber\\
\eqntext{\mbox{where }
\Upsilon_h(z,x,\gamma)=\mathbb{E} \bigl\{\bigl|V^{{p}/{s}}
\bigl(z+\kappa(x)Z_\gamma^h\bigr)-V^{{p}/{s}}(z)\bigr|^s \bigr\}.}
\end{eqnarray}
By \eqref{decomprev5} and \eqref{corequation1rev}, one checks that
\eqref{tuo} holds if there exists $C>0$ such that, for every
$z,x\in\mathbb{R}^d$ and $\gamma\le\gamma_1$,
%
\begin{equation}\label{corequation2rev}
\Upsilon_h(z,x,\gamma)\le C\gamma^{{s}/{2}}V^{a+p-1}(x),
\end{equation}
where $h=1$ (resp. $h=0$, resp. $h=+\infty$) if $p\le1/2<q$ (resp.
if $p,q\le1/2$, resp. if $p>1/2$). Then, it suffices to prove
\eqref{corequation2rev}. First, when $p\le1/2<q$, we have
$Z_\gamma^1=Z_\gamma=Y_\gamma+N_\gamma$. On the one hand, since
$V^{p/s}$ is a $2p/s$-H\"older function [see Lemma
\ref{lemmetech}(a)], it follows from
\ref{asS}.1 and Lemma \ref{moment}(i) that
%
\begin{eqnarray}\label{grsautrev}
\mathbb{E} \bigl\{\bigl|V^{{p}/{s}}\bigl(z+\kappa (x)N_\gamma\bigr)
-V^{{p}/{s}}(z)\bigr|^s \bigr\}
&\le & C\|\kappa(x)\|^{2p}\mathbb{E}\{|N_\gamma|^{2p}\}
\nonumber\\[-8pt]
\\[-8pt]
&\le& C\gamma V^{p+a-1}(x).
\nonumber
\end{eqnarray}

On the other hand, $V^{p/s}$ is a $2q/s$-H\"older function
when $q/s\le1/2$ (because $2p/s\le2q/s\le1$ in this case).
Hence, using this property if $q/s\le1/2$ and
the Taylor formula if $q/s>1/2$ yields
%
\begin{eqnarray}\label{ptsautrev1}
&& \mathbb{E}\bigl\{\bigl|V^{{p}/{s}}\bigl(z+\kappa(x)({N}_\gamma +{Y}_\gamma)\bigr)
- V^{{p}/{s}}\bigl(z+\kappa(x){N}_\gamma\bigr)\bigr|^{s}\bigr\}
\nonumber\\[-8pt]
\\[-8pt]
&&\qquad \le C \cases{%
\|\kappa(x)\|^{2q}\mathbb{E}\{|{Y}_\gamma|^{2q}\}, & \quad if $q/s\le1/2$,\cr
\mathbb{E}\{ |\langle V^{{p}/{s}-1}\nabla V(\xi_2),
\kappa(x)Y_\gamma\rangle |^s\}, & \quad if $q/s>1/2$,}
\nonumber
\end{eqnarray}
with $\xi_2\in[ z+\kappa(x)N_\gamma,z+\kappa(x)(N_\gamma+Y_\gamma)]$. By Lemma
\ref{moment}(ii).2, $\mathbb{E}\{|Y_\gamma|^{2q}\}\le C\gamma$. It
follows from Jensen's inequality that
%
\begin{equation}\label{ptsautrev2}
\qquad
\mathbb{E}\bigl\{\bigl|V^{{p}/{s}}\bigl(z+\kappa(x)Z_\gamma\bigr)
-V^{{p}/{s}}\bigl(z+\kappa(x){N}_\gamma\bigr)\bigr|^{s}\bigr\}
\le C\gamma^{{s}/{(2q)}\wedge1}\|\kappa(x)\|^{s\wedge2q}.
\end{equation}
One checks that $\|\kappa\|^{s\wedge2q}\le
CV^{p+a-1}$ under \ref{asS}.1 and that
$\gamma+\gamma^{{s}/{(2q)}\wedge1}\le C\gamma^{s/2}$ for
every $\gamma\le\gamma_1$. Hence, summing up \eqref{grsautrev}
and \eqref{ptsautrev2} and using \eqref{eqelem2}
yields~\eqref{corequation2rev} (with $h=1$) when $p\le1/2<q$.

When $p,q\le1/2$ (resp. $p>1/2$), we have to check that
\eqref{corequation2rev} holds with $h=0$ (resp. with $h=+\infty$).
Then, we need to use a decomposition of the jump component adapted
to the value of $h$. We split up $Z_\gamma^0$ (resp.
$Z_\gamma^\infty$) as follows:
$Z_\gamma^0=\check{Y}_\gamma+N_\gamma$ with
$\check{Y}_\gamma=Y_\gamma+\gamma\int_{\{|y|\le1\}}y\pi(dy)$
[resp. $Z_\gamma^\infty=Y_\gamma+\hat{N}_\gamma$ with
$\hat{N}_\gamma=N_\gamma-\gamma\int_{\{|y|>1\}}y\pi(dy)$]. Then,
when $p,q\le1/2$ (resp. $p>1/2$), the idea is to replace
$Y_\gamma$ with $\check{Y}_\gamma$ (resp. $N_\gamma$ with
$\hat{N}_\gamma$) in the left-hand sides of \eqref{grsautrev} and
\eqref{ptsautrev1} and to derive some adapted controls from Lemma
\ref{moment} and inequality \eqref{nth3}. Since the
proof is close to that of the $p\le1/2<q$ case, we leave it to the
reader.

\textbf{Case} $\bolds{p/s> 1/2.}$
Since $p>1/2$, we use the notation introduced in \eqref{955}
with $h=+\infty$. We recall that
$\bar{X}_{n-1}+\gamma_n\tilde{b}(\bar{X}_{n-1})=\bar{X}_{n,1}^\infty$.
Then, applying the following inequality,
%
\begin{equation}
\forall u,v \ge0, \forall\alpha\ge1 \qquad
|u^\alpha -v^\alpha| \le C_\alpha (|u-v|u^{\alpha-1}+|u-v|^{\alpha} )
\end{equation}
with $u=\sqrt{V}(\bar{X}_n)$, $v=\sqrt{V}(\bar{X}_{n,1}^\infty)$
and $\alpha=(2p)/s$, we obtain
\begin{eqnarray*}
|V^{{p}/{s}}(\bar{X}_n)-V^{{p}/{s}}(\bar{X}_{n,1}^\infty )|
&\le & C\bigl|\sqrt{V}(\bar{X}_n)
-\sqrt{V}(\bar{X}_{n,1}^\infty)\bigr|V^{{p}/{s}-{1}/{2}}(\bar{X}_{n-1})
\\
&&{} +C\bigl|\sqrt{V}(\bar{X}_n)-\sqrt{V}(\bar{X}_{n,1}^\infty
)\bigr|^{2p/s}.
\end{eqnarray*}
We deduce from \ref{asS}.1 that
\begin{eqnarray*}
|\bar{X}_n-\bar{X}_{n,1}^\infty|\le \cases{%
CV^{(a+p-1)/(2p)}(\bar{X}_{n-1})
\bigl(\sqrt{\gamma_n}|U_{n}|+|\bar{Z}^\infty_n| \bigr), & \quad if $p< 1$,\cr
CV^{a/2}(\bar{X}_{n-1}) \bigl(\sqrt{\gamma_n}|U_{n}|+|\bar{Z}^\infty_n| \bigr),
& \quad if $p\ge1$.}
\end{eqnarray*}
Since $\sqrt{V}$ is Lipschitz, one then checks that
\begin{eqnarray}
&& |V^{{p}/{s}}(\bar{X}_n)-V^{{p}/{s}}(\bar{X}_{n,1}^\infty)|
\nonumber\\
&&\qquad
\le CV^r(\bar{X}_{n-1}) \bigl(\sqrt{\gamma_n}|U_n|+|\bar{Z}^\infty_n|
+\gamma_n^{p/s}|U_n|^{2p/s}+|\bar{Z}^\infty_n|^{2p/s} \bigr),
\nonumber\\
\eqntext{\mbox{where }
r= \cases{%
\biggl(\dfrac{p}{s}+\dfrac{a-1}{2p}\biggr)\vee\biggl(\dfrac{a+p-1}{s}\biggr),
& \quad if $p<1$,\cr
\biggl(\dfrac{p}{s}+\dfrac{a-1}{2}\biggr)\vee\dfrac{ap}{s},
& \quad if $p\ge 1$.}}
\end{eqnarray}
One derives from Lemma \ref{moment}(iii) and from the Jensen
inequality that, for $\alpha>0$,
$\mathbb{E}\{|\bar{Z}_n^\infty|^\alpha\}=O(\gamma_n^{({\alpha}/{2})\wedge1})$.
Therefore, since $2p/s\ge1/2$ and $f_{a,p}(s)\le2$, we have
\[
\mathbb{E} \bigl\{ \bigl(\sqrt{\gamma_n}|U_n|+|\bar{Z}^\infty
_n|+\gamma_n^{p/s}|U_n|^{2p/s}+|\bar{Z}^\infty_n|^{2p/s}
\bigr)^{f_{a,p}(s)} \bigr\}=O\bigl(\gamma_n^{{f_{a,p}(s)}/{2}}\bigr).
\]
Second, one deduces from the definition of $f_{a,p}$ that
$r f_{a,p}(s)\le a+p-1$. Therefore, inequality \eqref{tuo} follows.
\end{pf}

By Lemma \ref{explication}, Corollary \ref{cor} concludes the
proof of Proposition \ref{tension} and then, the part which is
concerned with the tightness of $(\bar{\nu}_n(\omega,dx))_{n\ge 1}$.
The only thing left to prove the theorem for Scheme (A) is
thus to identify the limit. This is the aim of the next section.

\section{ Identification of the weak limits of $(\bar{\nu}_n(\omega,dx))_{n\ge1}$}\label{section4}

In this section we show that every weak
limiting distribution of $(\bar{\nu}_n(\omega,dx))_{n\ge1}$ is
invariant for $(X_t)_{t\ge0}$. For this purpose, we will rely on
the Echeverria--Weiss theorem (see~\cite{bib4}, page~238, \cite{bib2}
and \cite{lemaire}). This is a criterion for invariance based on
the infinitesimal generator $A$ of $(X_t)$ defined by
\eqref{1110}. By the Echeverria--Weiss theorem, we know that if
$A(\mathcal{C}^2_K(\mathbb{R}^d))\subset\mathcal{C}_0(\mathbb{R}^d)$, a probability
$\nu$ is invariant for the SDE if for every $f\in\mathcal{C}^2_K(\mathbb{R}^d)$,
$\nu(Af)=0$. One can check that $A(\mathcal{C}^2_K(\mathbb{R}^d))\subset\mathcal{C}_0(\mathbb{R}^d)$ if
$\|\kappa(x)\|=o(|x|)$ when $|x|\rightarrow+\infty$
(and this condition cannot be improved in general). Hence, under this
condition on $\kappa$,
it follows that every weak limiting distribution of
$(\bar{\nu}_n)$ is invariant if for every
$f\in\mathcal{C}^2_K(\mathbb{R}^d),$ $\bar{\nu}_n(Af)\rightarrow0$.
The main result of this section is then the following proposition.
\begin{prop}\label{invariance}
Let $a\in(0,1]$, $p>0$ and $q\in[0,1]$. Assume
\hyperlink{eH}{$\mathrm{(H^1_p)}$},
\hyperlink{eH}{$\mathrm{(H^2_q)}$},
\textup{\ref{asS}.1}.
Assume that $\|\kappa(x)\|\stackrel{|x|\rightarrow+\infty}=o(|x|)$
and that $(\eta_n/\gamma_n)_{n\ge1}$ is nonincreasing. If, moreover,
%
\begin{eqnarray}\label{assuinv}
\sup_{n\ge1}\bar{\nu}_n\bigl(\|\kappa\|^{2q}
+\operatorname{Tr}(\sigma\sigma^*)\bigr) &<& \infty \quad \mbox{and}\quad
\nonumber\\[-8pt]
\\[-8pt]
\sum_{k\ge1}\frac{\eta_k^2}{H^2_k\gamma_k}\mathbb{E}
\{V^{a+p-1}(\bar{X}_{k-1})\} &<& +\infty,
\nonumber
\end{eqnarray}
then,
%
\begin{equation}
\forall f\in  \mathcal{C}^2_K(\mathbb{R}^d), \mbox{ a.s.}, \qquad
\int Af\,d\bar{\nu}_n\stackrel{n\rightarrow\infty}{\longrightarrow}0.
\end{equation}
Consequently, \textit{a.s.}, every weak limiting distribution of
$(\bar{\nu}_n(\omega,dx))_{n\ge1}$ is invariant for the SDE
\eqref{edss}.
\end{prop}
\begin{Remarque} 
This proposition is sufficient to conclude the proof
because the two assumptions in \eqref{assuinv} hold under the
assumptions of Theorem \ref{principal} (resp. Theorem
\ref{principal'}). Indeed, since
$\|\kappa\|^{2q}+\operatorname{Tr}(\sigma\sigma^*)\le CV^{{p}/{s}+a-1}$ with $s=2$
in Theorem \ref{principal} [resp. with $s$ satisfying \eqref{cond1} in
Theorem \ref{principal'}], the first is a consequence of
Proposition~\ref{tension}. Likewise, the second is a consequence
of Lemma \ref{lemme4} applied with
$\theta_n=(\eta_n/(H_n\gamma_n))^2$ [see \eqref{eirop}].
\end{Remarque}

\subsection[Proof of Proposition 3]{Proof of Proposition \textup{\protect\ref{invariance}}}

The proof of Proposition \ref{invariance} is built in two
successive steps that are represented by Propositions \ref{prop1}
and \ref{prop2}. In Proposition~\ref{prop1} we claim that
showing that $\bar{\nu}_n(Af)\rightarrow0$ {a.s.} is equivalent to
showing that
$1/H_n\sum_{k=1}^n(\eta_k/\gamma_k)\mathbb{E}\{f(\bar{X}_{k})
-f(\bar{X}_{k-1})|\mathcal{F}_{k-1}\}\rightarrow0$
{a.s.} Then, in Proposition \ref{prop2} we show that this last
term does tend to 0.
\begin{prop}\label{prop1}
Assume that the assumptions of Proposition
\ref{invariance} are fulfilled. Then,
for every $f\in\mathcal{C}^2_K(\mathbb{R}^d)$,
%
\begin{equation}
\qquad
\lim_{n\rightarrow\infty}\frac{1}{H_n}
\sum_{k=1}^{n}\eta_k \biggl(\frac{\mathbb{E}\{f(\bar{X}_k)-f(\bar
{X}_{k-1})|\mathcal{F}_{k-1}\}}{\gamma_k}-Af(\bar{X}_{k-1}) \biggr)
=0\qquad \mbox{a.s.}
\end{equation}
\end{prop}

We begin the proof by a technical lemma.
\begin{lemme} \label{souvent} 
Let $\Phi \dvtx\mathbb{R}^d\mapsto \mathbb{R}^l$ be
a continuous function with compact support,
$\Psi\dvtx \mathbb{R}^d\mapsto\mathbb{R}_+$, a locally bounded function,
$(h_1^{\theta})_{\theta\in[0,1]}$
and $(h_2^{\theta})_{\theta\in[0, 1]}$ two families of Borel
functions defined
on $\mathbb{R}^d\times\mathbb{R}_+$ with values in $\mathbb{R}^d$
satisfying the
following assumptions:
\begin{itemize}
\item There exists $\delta_0>0$ such that
%
\begin{equation}\label{assprop21}
\inf_{\theta\in[0,1],\gamma\in[0,\delta_0]}
\bigl(|h^\theta_1|(x,\gamma)+ |h^\theta_2|(x,\gamma) \bigr)
\stackrel{|x|\rightarrow+\infty}{\longrightarrow}+\infty.
\end{equation}
\item For every compact set $K$,
%
\begin{equation}\label{assprop23}
\sup_{x\in K,\theta\in[0,1]}|h^\theta_1(x,\gamma)-h^\theta_2(x,\gamma)|
\stackrel{\gamma\rightarrow0}{\longrightarrow}0.
\end{equation}
\end{itemize}
Then, for
every sequence $(x_k)_{k\in\mathbb{N}}$ of $\mathbb{R}^d$,
\[
\frac{1}{H_n}\sum_{k=1}^n \eta_k\sup_{\theta\in[0,1]}\|\Phi
(h^\theta_1(x_{k-1},\gamma_k))-\Phi(h^\theta_2(x_{k-1},\gamma_k))\|
\Psi(x_{k-1})\stackrel{n\rightarrow+\infty}{\longrightarrow}0.
\]
\end{lemme}
\begin{pf}
$\Phi$ has a compact support, therefore, we derive from
\eqref{assprop21} that there exists $M_{\delta_0}>0$ such that,
for every $|x|> M_{\delta_0}$, $\gamma\le\delta_0$ and
$\theta\in[0,1]$,
\[
\Phi(h^\theta_1(x,\gamma))=\Phi(h^\theta_2(x,\gamma))=0.
\]
Consider $\rho\mapsto
w(\rho,\Phi)=\sup\{\eta>0,\sup_{|x-y|\le\eta}|\Phi(x)-\Phi
(y)|\le \rho\}$. As $\Phi$ is uniformly
continuous, $w(\rho,\Phi)>0$ for every $\rho>0$. Thanks to
\eqref{assprop23}, for every $\rho>0$, there exists
$\delta_\rho\le\delta_0$ such that,
for every $\gamma\le\delta_\rho$, $\theta\in[0,1]$,
\[
\sup_{|x|\le M_{\delta_0}}|h^\theta_1(x,\gamma)-h^\theta
_2(x,\gamma)|\le w(\rho,\Phi).
\]
As $\gamma_k\stackrel{k\rightarrow+\infty}{\longrightarrow}0$, there exists
$k_\rho\in\mathbb{N}$ such that $\gamma_k\le\delta_\rho$ for
$k\ge k_\rho$. By using that
$H_n\stackrel{n\rightarrow+\infty}{\longrightarrow}+\infty$, we deduce that
\[
\limsup_{n\rightarrow+\infty}\frac{1}{H_n}\sum_{k=1}^n\eta_k
\sup_{\theta\in[0,1]}\|\Phi(h^\theta_1(x_{k-1},\gamma_k))-\Phi
(h^\theta_2(x_{k-1},\gamma_k))\|
\Psi(x_{k-1})\le\rho\overline{\Psi}_{\delta_0},
\]
where
$\overline{\Psi}_{\delta_0}:=\sup\{|\Psi(x)|,|x|\le
M_{\delta_0}\}<+\infty$ since $\Psi$ is locally bounded. The
result follows.
\end{pf}
\begin{pf*}{Proof of Proposition \ref{prop1}}
We have to inspect successively the $q\in(1/2,1]$ and
$q\in[0,1/2]$ cases.

\textbf{Case} $\bolds{q\in(1/2,1].}$ Let
$f\in \mathcal{C}^2_K(\mathbb{R}^d)$.
Decompose the infinitesimal generator as the sum of three terms defined by
\begin{eqnarray*}
A_1 f(x) &=& \langle\nabla f,b\rangle(x), \qquad
A_2 f(x)=\operatorname{Tr}(\sigma^* D^2 f\sigma)(x),
\\
A_3 f(x) &=& \int\bigl(f\bigl(x+\kappa(x)y\bigr)-f(x)-\langle\nabla
f(x),\kappa(x)y\rangle1_{\{|y|\le1\}}\bigr)\pi(dy).
\end{eqnarray*}
Set $\bar{X}_{k,1}=\bar{X}_{k-1}+\gamma_kb(\bar{X}_{k-1})$,
$\bar{X}_{k,2}=\bar{X}_{k,1}+\sqrt{\gamma_k}\sigma(\bar{X}_{k-1})U_k$
and
$\bar{X}_{k,3}=\bar{X}_{k,2}+\sqrt{\gamma_k}\kappa(\bar{X}_{k-1})\bar{Z}_k$.
We then part the proof into three steps:
\begin{eqnarray}
&& \mbox{\textbf{Step 1.}} \quad
\frac{1}{\gamma_k}\mathbb{E} \{f(\bar{X}_{k,1})-f(\bar{X}_{k-1})/\mathcal{F}_{k-1} \}
=A_1 f(\bar{X}_{k-1})+ R_1(\gamma_k,\bar{X}_{k-1})
\nonumber\\
\eqntext{\mbox{with } \displaystyle \frac{1}{H_n}\sum_1^n\eta_k
R_1(\gamma_k,\bar{X}_{k-1})\stackrel{n\rightarrow\infty}{\longrightarrow}0.}
\\
&& \mbox{\textbf{Step 2.}} \quad
\frac{1}{\gamma_k}\mathbb{E}
\{f(\bar{X}_{k,2})-f(\bar{X}_{k,1})|\mathcal{F}_{k-1} \}
=A_2 f(\bar{X}_{k-1})+ R_2(\gamma_k,\bar{X}_{k-1})
\nonumber\\
\eqntext{\mbox{with } \displaystyle
\frac{1}{H_n}\sum_1^n\eta_k
R_2(\gamma_k,\bar{X}_{k-1})\stackrel{n\rightarrow\infty}{\longrightarrow}0.}
\\
&& \mbox{\textbf{Step 3.}} \quad
\frac{1}{\gamma_k}\mathbb{E}\{f(\bar{X}_k)-f(\bar{X}_{k,2})|\mathcal{F}_{k-1}\}
=A_3 f(\bar{X}_{k-1})+ R_3(\gamma_k,\bar{X}_{k-1})
\nonumber\\
\eqntext{\displaystyle  \mbox{with }
\frac{1}{H_n}\sum_1^n\eta_k R_3(\gamma_k,\bar{X}_{k-1})
\stackrel{n\rightarrow\infty}{\longrightarrow}0.}
\end{eqnarray}
The combination of the three steps yields Proposition \ref{prop1}.
We refer to Proposition~4 of \cite{bib3} for steps 1 and 2 and
focus on the last step where the specificity of our jump L\'evy
setting appears. Since $\bar{X}_{k-1}$ is $\mathcal{F}_{k-1}$-measurable and
$\bar{Z}_k$, $U_k$ and $\mathcal{F}_{k-1}$ are independent, we have
\begin{eqnarray}
\mathbb{E}\bigl\{f\bigl(\bar{X}_{k,2}+\kappa(\bar{X}_{k-1})\bar
{Z}_k\bigr)|\mathcal{F}_{k-1}\bigr\}
=Q_{\gamma_k}f(\bar{X}_{k-1}),
\nonumber\\
\eqntext{\displaystyle \mbox{where }
Q_\gamma f(x)=\int_{\mathbb{R}^d}\mathbb{E}\bigl\{f\bigl(S_{x,\gamma,u}
+\kappa(x)Z_\gamma\bigr)\bigr\}\mathbb{P}_{U_1}(du),}
\end{eqnarray}
with $S_{x,\gamma,u}=x+\gamma b(x)+\sqrt{\gamma}\sigma(x)u$. Set
$V_t=S_{x,\gamma,u} +\kappa(x)Z_t$. Applying It\^o's formula to
$(f(V_t))_{t\ge0}$ yields
%
\begin{eqnarray}\label{tildeH}
f(V_t) &=& f(S_{x,\gamma,u})+ \int_0^t
\langle\nabla f(V_{s^-}),\kappa(x)\,dY_s\rangle
\nonumber\\[-8pt]
\\[-8pt]
&&{} +\sum_{0<s\le t}\tilde{H}^f \bigl(S_{x,\gamma,u}+\kappa
(x)Z_{s^-},x,\Delta Z_s \bigr)
\nonumber
\end{eqnarray}
\begin{equation} \label{tildeH2}
\hspace*{45pt}
\mbox{where, } \tilde{H}^f(z,x,y)=f\bigl(z+\kappa(x)y\bigr)
-f(z)- \langle\nabla f(z),\kappa(x)y\rangle1_{\{|y|\le1\}}.
\end{equation}
The process $(\int_0^t \langle\nabla f(V_{s^-}),\kappa(x)\,dY_s\rangle)$
is a true martingale
since $\nabla f$ is bounded. The compensation formula and a change of variable
yield
\begin{eqnarray*}
&& \mathbb{E}\bigl\{f\bigl(S_{x,\gamma,u} +\kappa(x)Z_\gamma\bigr)\bigr\}
\\
&&\qquad = \mathbb{E}\{f(V_{\gamma})\}
\\
&&\qquad = f(S_{x,\gamma,u})+\gamma\mathbb{E} \biggl\{\int_0^1 dv
\int \pi(dy)\tilde{H}^f\bigl(S_{x,\gamma,u}+\kappa(x)Z_{v\gamma},x,y\bigr) \biggr\}.
\end{eqnarray*}
Since $ A_3 f(x)=\int\pi(dy)\tilde{H}^f(x,x,y)=\mathbb{E} \{
\int_0^1
dv\int\pi(dy)\tilde{H}^f(x,x,y) \},$ it follows from the
previous inequality that
\[
\frac{1}{\gamma_k}\mathbb{E} \{f(\bar{X}_k)-f(\bar{X}_{k,2})|\mathcal{F}_{k-1} \}
=A_3 f(\bar{X}_{k-1})+ R_3(\gamma_k,\bar{X}_{k-1}),
\]
where,
\[
R_3(\gamma,x)=\int\mathbb{E} \biggl\{\int_0^1 dv\int\pi(dy)
\Delta\tilde{H}^f\bigl(S_{x,\gamma,u}+\kappa(x)Z_{v\gamma},x,x,y\bigr)
\biggr\} \mathbb{P}_{U_1}(du)
\]
with $\Delta\tilde{H}^f (z_1,z_2,x,y )=\tilde{H}^f(z_1,x,y)-\tilde{H}^f(z_2,x,y)$.
We upper-bound $R_3$ by two terms: $R_{3,1}$ and $R_{3,2}$ that
are associated to the small and big jumps components of
$(Z_t)$, namely,
\begin{eqnarray*}
&& R_{3,1}(\gamma,x)
\\
&&\qquad = \int\int_0^1 dv\int_{\{|y|\le1\}}\pi(dy)\mathbb{E}
\bigl |\Delta\tilde{H}^f
\bigl(S_{x,\gamma,u}+\kappa(x)Z_{v\gamma},x,x,y \bigr)\bigr |
\mathbb{P}_{U_1}(du),
\\
&& R_{3,2}(\gamma,x)
\\
&&\qquad = \int\int_0^1 dv\int_{\{|y|>1\}}\pi(dy)\mathbb{E}
\bigl |\Delta\tilde{H}^f
\bigl(S_{x,\gamma ,u}+\kappa(x)Z_{v\gamma},x,x,y \bigr)\bigr |
\mathbb{P}_{U_1}(du).
\end{eqnarray*}
We study successively $R_{3,1}$ and $R_{3,2}$. From Taylor's
formula, we have for every $y$ such that $|y|\le1$
\[
\bigl |\Delta\tilde{H}^f \bigl(S_{x,\gamma,u}+\kappa(x)Z_{v\gamma},
x,x,y \bigr) \bigr|
\le\tfrac{1}{2} R(Z,\gamma,x,u,v,y)|\kappa(x)y|^{2},
\]
where
\begin{eqnarray*}
&& R(Z,\gamma,x,u,v,y)
\\
&&\qquad =\sup_{\theta\in[0,1]}\bigl\|D^2
f \bigl(S_{x,\gamma,u} +\kappa(x)(Z_{v\gamma}+\theta y) \bigr)
-D^2 f\bigl(x +\theta\kappa(x)y\bigr)\bigr\|.
\end{eqnarray*}
By setting $\Phi=D^2f$, $\Psi(x)=\|\kappa(x)\|^2|y|^2$,
\[
h^\theta_1(x,\gamma)=S_{x,\gamma,u} +\kappa(x)(Z_{v\gamma}+\theta y)
\quad\mbox{and}\quad
h^\theta_2(x,\gamma)=x +\theta\kappa(x)y,
\]
we want to show that the assumptions of Lemma \ref{souvent} are
{a.s.} fulfilled for every \textit{fixed} $u,$ $v$ and $y$.

First, since $\kappa(x)\stackrel{|x|\rightarrow+\infty}{=}o(|x|)$,
there exists a continuous function $\varepsilon$ such that
$\kappa(x)=|x|\varepsilon(x)$ and $\varepsilon(x)\stackrel{|x|\rightarrow\infty}{\longrightarrow}0$.
Therefore, as $b$ and $\sigma$ have sublinear growth, one
checks that there exist
some positive real constants $C_1$ and $C_2$ such that
%
\begin{eqnarray}\label{56bis1}
\qquad\quad
\cases{%
|S_{x,\gamma,u} +\kappa(x)(Z_{v\gamma}+\theta y)|
\ge|x| \bigl(1-\gamma C_1-(|Z_{v\gamma}|+|y|)|\varepsilon(x)|\bigr)-C_2 ,\cr
|x+\theta\kappa(x)y|\ge|x| \bigl(1-|\varepsilon(x)||y| \bigr).}
\end{eqnarray}
Let $\delta_0$ be a positive number such that $1-\delta_0
C_1>0$. Since $(Z_t)$ is locally bounded (as a c\`adl\`ag
process) and $\varepsilon(x)\stackrel{|x|\rightarrow\infty}{\longrightarrow}0$,
there exists {a.s.} $M>0$,
\[
\inf_{|x|>M,\gamma\in[0,\delta_0]}
\bigl(1-\gamma C_1-(|Z_{v\gamma }|+|y|)|\varepsilon(x)| \bigr)>0.
\]
It follows that {a.s.},
\[
\inf_{\theta\in[0,1],\gamma\in[0,\delta_0]}
\bigl (|h^\theta_1|(x,\gamma)+|h^\theta_2|(x,\gamma)\bigr)
\stackrel{|x|\rightarrow\infty}{\longrightarrow}+\infty.
\]
Second, let $K$ be a compact set of $\mathbb{R}^d$. We check that
\eqref{assprop23} holds. We have
%
\begin{eqnarray}\label{56bis2}
&& \sup_{x\in K,\theta\in[0,1]}|h^\theta_1(x,\gamma)-h^\theta_2(x)|
\nonumber\\[-8pt]
\\[-8pt]
&&\qquad \le\sup_{x\in K}
\bigl(\gamma|b(x)|+\sqrt{\gamma}\|\sigma(x)\||u|+\|\kappa(x)\|
\,|Z_{v\gamma}| \bigr)
\stackrel{\gamma\rightarrow0}{\longrightarrow}0\qquad \mbox{a.s.}
\nonumber
\end{eqnarray}
because $b$, $\sigma$, $\kappa$ are locally bounded and
$\lim_{t\rightarrow0}Z_t=0$ a.s. Thus, by Lemma \ref{souvent},
for any sequence $(x_k)_{k\in\mathbb{N}}$ of $\mathbb{R}^d$, for every
$(u,v,y)\in\mathbb{R}^d\times[0,1]\times{ B}_d(0,1)$,
%
\begin{equation}\label{956bis}
\qquad
\frac{1}{H_n}\sum_{k=1}^{n}\eta_{k}\Delta\tilde{H}^f
\bigl(S_{{x_{k-1}},\gamma_k,u}+\kappa({x_{k-1}})Z_{v\gamma_k},
{x_{k-1}},{x_{k-1}},y\bigr)
\stackrel{n\rightarrow\infty}{\longrightarrow} 0 \qquad\mbox{a.s.}
\end{equation}
Now, since $\nabla f$ and $D^2 f$ are bounded, we derive from
Taylor's formula that, for every $z_1$, $z_2\in\mathbb{R}^d$,
\[
|\tilde{H}^f(z_2,x,y)-\tilde{H}^f(z_1,x,y)|1_{\{|y|\le1\}}
\le \cases{%
2\|\nabla f\|_\infty\|\kappa(x)\|\,|y|1_{\{|y|\le1\}},\cr
2\|D^2f\|_\infty\|\kappa(x)\|^2|y|^21_{\{|y|\le1\}}.}
\]
Then, for every $q\in(1/2,1]$,
%
\begin{equation}\label{956}
\quad
\bigl |\Delta\tilde{H}^f \bigl(S_{{x},\gamma,u}+\kappa(x)Z_{v\gamma},
{x},{x},y \bigr) \bigr|1_{\{|y|\le1\}}
\le C\|\kappa(x)\|^{2q}|y|^{2q}1_{\{|y|\le1\}},
\end{equation}
where
$C=2\max(\|\nabla f\|_\infty,\|D^2f\|_\infty)$. Therefore, by
assumption \hyperlink{eH}{$\mathrm{(H^2_q)}$}, we finally derive from
\eqref{956bis}, \eqref{956} and from the Lebesgue dominated
convergence theorem that
%
\begin{eqnarray}\label{956bisbis}
&& \frac{1}{H_n}\sum_{k=1}^{n}\eta_{k}
R_{3,1}(\gamma_k,x_{k-1})\stackrel{n\rightarrow\infty}{\longrightarrow} 0
\nonumber\\[-8pt]
\\[-8pt]
\eqntext{\displaystyle \mbox{if }
\sup_{n\in\mathbb{N}}\frac{1}{H_n}\sum
_{k=1}^{n}\eta_{k}\|\kappa(x_{k-1})\|^{2q}<\infty.}
\end{eqnarray}
We apply this result to $(x_k)=(\bar{X}_k)$. By
\eqref{assuinv},
$\sup_{n\in\mathbb{N}}\bar{\nu}_n(\|\kappa\|^{2q})<\infty$ a.s.
Hence, it follows that ${1}/{H_n}\sum_{k=1}^{n}\eta_{k}
R_{3,1}(\gamma_k,\bar{X}_{k-1})\stackrel{n\rightarrow\infty}{\longrightarrow} 0$ a.s.

Now, let us focus on $R_{3,2}$. Set
$\Delta f(z_1,z_2)=f(z_1)-f(z_2)$. Then,
\begin{eqnarray*}
R_{3,2}(\gamma,x) &=& \int\mathbb{E} \biggl\{\int_0^1
dv\int_{\{|y|>1\}}\pi(dy)\Delta f \bigl(x,S_{x,\gamma,u}
+\kappa (x)Z_{v\gamma} \bigr) \biggr\}\mathbb{P}_{U_1}(du)
\\
&&{}+ \int\mathbb{E} \biggl\{\int_0^1 dv\int_{\{|y|>1\}} \pi(dy)
\\
&&\hspace*{35pt}{}\times
\Delta f \bigl(S_{x,\gamma,u}+\kappa(x)
(Z_{v\gamma}+y),x+\kappa(x)y \bigr) \biggr\}\mathbb{P}_{U_1}(du).
\end{eqnarray*}
One proceeds as before. By using Lemma \ref{souvent}, one begins
by showing that, for any sequence $(x_k)_{k\in\mathbb{N}}$ , for
every $(u,v,y)\in[0,1]\times\mathbb{R}^d\times B_d(0,1)^c$, {a.s.},
%
\begin{eqnarray}\label{ur}
\qquad
&& \frac{1}{H_n}\sum_{k=1}^{n}\eta_{k}\Delta f
\bigl(x_{k-1},S_{x_{k-1},\gamma_k,u}+\kappa(x_{k-1})Z_{u\gamma_k} \bigr)
\stackrel{n\rightarrow\infty}{\longrightarrow} 0\quad\mbox{and}
\nonumber\\[-8pt]
\\[-8pt]
&& \frac{1}{H_n}\sum_{k=1}^{n}\eta_{k}\Delta
f \bigl(S_{x_{k-1},\gamma_k,u}+\kappa(x_{k-1})(Z_{v\gamma
_k}+y),x_{k-1}+\kappa(x_{k-1})y \bigr)
\stackrel{n\rightarrow\infty}{\longrightarrow}0.
\nonumber
\end{eqnarray}
By the dominated convergence theorem [which can be applied
because\break
\mbox{$\pi(|y|>1)<\infty$} and $f$ is bounded], we deduce that, for any
sequence $(x_k)_{k\in\mathbb{N}}$,
\[
\frac{1}{H_n}\sum_{k=1}^{n}\eta_{k}
R_{3,2}(\gamma_k,x_{k-1})
\stackrel{n\rightarrow\infty}{\longrightarrow} 0
\qquad\mbox{a.s.}
\]
This completes the proof of Step 3 when $q\in(1,2]$.

\textbf{Case} $\bolds{q\le1/2.}$
The reader can check that the
assumption $q\in(1/2,1]$ is used only once: when we want to apply
the dominated convergence theorem for $R_{3,1}$ [see \eqref{956}].
Since inequality \eqref{956} is not true when $q<1/2$, we need to
decompose the infinitesimal generator in a slightly different way:
\begin{eqnarray*}
A_1 f(x) &=& \langle\nabla f,b^0\rangle(x), \qquad
\\
A_2 f(x) &=& \operatorname{Tr}(\sigma^* D^2 f\sigma)(x),
\\
A_3 f(x) &=& \int\bigl(f\bigl(x+\kappa(x)y\bigr)-f(x)\bigr) \pi(dy).
\end{eqnarray*}
Note that this decomposition is only possible when $q\le1/2$. That
means that with the notation \eqref{955}, we decompose $\Delta
\bar{X}_k$ with $h=0$ and inspect the three induced steps. We do
not go into further details since the proof is similar to the
case \mbox{$q>1/2$}.
\end{pf*}
\begin{prop}\label{prop2}
Assume that the assumptions of Proposition \ref{invariance} are
fulfilled. Then,
\[
\lim_{n\rightarrow\infty}\frac{1}{H_n}\sum
_{k=1}^n\frac{\eta_k}{\gamma_k}\mathbb{E}\{f(\bar{X}_k)
-f(\bar{X}_{k-1})|\mathcal{F}_{k-1}\}=0 \qquad \mbox{a.s.}
\]
\end{prop}
\begin{pf}
We do not detail the proof of this proposition which is
an adaptation of Proposition 3 in \cite{bib3}.
\end{pf}

\section[Proof of the main theorems for Schemes (B) and (C)]{Proof of the main theorems for Schemes (\protect\ref{eB}) and (\protect\ref{eC})}\label{section5}

The aim of this section is to give a general
idea of the proof for Schemes (\ref{eB}) and (\ref{eC}) and to overcome the
main difficulties induced by the approximation of the jump
component. For Scheme (\ref{eA}), main theorems have been proven in two
successive steps. First, we focused on tightness results
(Proposition \ref{tension}) and then proved that every weak
limiting distribution is invariant for $(X_t)_{t\ge0}$
(Proposition \ref{invariance}). We follow the same process for
Schemes (\ref{eB}) and (\ref{eC}). We will successively explain for both
schemes why Proposition \ref{tension} and Proposition
\ref{invariance} remain valid.

\subsection{Almost sure tightness of $\bar{\nu}^{B}_n(\omega,dx)$
and $\bar{\nu}^{C}_n(\omega,dx)$}

The tightness result for Schemes (\ref{eB}) and (\ref{eC}) is strictly
identical to Proposition \ref{tension} [in particular, assumption
\eqref{condsupplementaire} is not necessary for tightness].
Looking carefully into the proof of this theorem for Scheme (\ref{eA})
shows that the properties of the jumps that we use are
the following: the
control of the moments of the jump components (Lemma~\ref{moment})
which is fundamental for Proposition \ref{prop9}, and independence
between $(\bar{Y}_n)_{n\in\mathbb{N}}$, $(\bar{N}_n)_{n\in\mathbb
{N}}$ and
$(U_n)_{n\in\mathbb{N}}$. We show in Lemma \ref{moment'} below that the
controls of Lemma \ref{moment} hold true for the moments of the
jump components of Schemes (\ref{eB}) and (\ref{eC}). Then, since Scheme (\ref{eB})
satisfies the independence properties, Proposition \ref{tension}
follows in this case. In Scheme (\ref{eC}), $(\bar{Y}^{C}_n)_{n\in\mathbb{N}}$
and $(\bar{N}^{C}_n)_{n\in\mathbb{N}}$ are no longer independent. It
raises several technical difficulties in the proof of Proposition
\ref{prop9} in case $p<1$, but the process of the proof is the
same. So, we only state a variant of Lemma \ref{moment} (see~\cite{bib26} for details).
\begin{lemme}\label{moment'}
Let $T_0$ be a positive number and
$T^n=\inf\{s> 0, |\Delta Z_t^n|>0\}$.
\begin{longlist}
\item Let $p>0$ such that
\hyperlink{eH}{$\mathrm{(H^1_p)}$} holds. Then, for every $t\le T_0$ and $h>0$,
\[
\mathbb{E}\{|N^h_{t\wedge T^n}|^{2p}\}\le t\int_{|y|>h}|y|^{2p}\pi
(dy)\qquad\mbox{if } p>0.
\]
\item Let $\tau$ be an $(\mathcal{F}_t)$-stopping time and
$q\in[0,1]$ such that
\hyperlink{eH}{$\mathrm{(H_q^2)}$} holds. Set $D_n^h=\{y,|y|\in(u_n,h]\}$
and $Y_t^{{h},n}=\sum_{0<s\le t}\Delta Y_s^{h}
1_{\{\Delta Y_s^{h}\in D_n^h\}}-t\int_{D_n^h}y\pi(dy)$. Then,
\[
\cases{%
\displaystyle
\mathbb{E}\biggl\{\biggl|Y^{h,n}_{t\wedge\tau}+(t\wedge
\tau) \int_{D_n^h} y\pi(dy)\biggr|^{2q}\biggr\}
\le t\int_{|y|\le h}|y|^{2q} \pi(dy), & \quad if $q\in[0,1/2]$,\cr
\displaystyle
\mathbb{E}\{|Y^{h,n}_{t\wedge\tau}|^{2q}\}\le C_q t
\int_{|y|\le h}|y|^{2q}\pi(dy), & \quad  if $q\in(1/2,1]$.}
\]
\item Let $p\ge1$ such that
\hyperlink{eH}{$\mathrm{(H^1_p)}$} holds. Set
$\hat{Z}^n_t=Z_t^n-t\int_{\{|y|>1\}}y\pi(dy)$. Then, there exists
$\eta>1$ such that, for every $T_0>0$, for every $\varepsilon>0$,
there exists \mbox{$C_{\varepsilon,T_0,p}>0$}, $n_0\in\mathbb{N}$ such that, for
every $t\ge T_0$ and $n\ge n_0$,
\[
\mathbb{E}\{|\hat{Z}^n_{t}|^{2p}\}\le
t \biggl(\int|y|^{2p}\pi(dy)+\varepsilon \biggr)
+ C_{\varepsilon,T_0,p} t^{\eta}
\]
and
\[
\mathbb{E}\{|\hat{Z}^n_{t\wedge T^n}|^{2p}\}\le t
\biggl(\int|y|^{2p}\pi(dy)+\varepsilon \biggr)
+C_{\varepsilon,T_0,p} t^{\eta}.
\]
\end{longlist}
\end{lemme}
\begin{pf}
The proof is left to the reader.
\end{pf}
\begin{Remarque}
In (iii), the control is only valid for $n$ sufficiently large
but that does not make any problem since
\hyperlink{e15}{$\mathrm{(R_{a,p})}$}
just needs to be valid for sufficiently large~$n$.
\end{Remarque}

\subsection{Identification of the limit of $(\bar{\nu}^{B}_n)_{n\in
\mathbb{N}}$ and $(\bar{\nu}^{C}_n)_{n\in\mathbb{N}}$}

The theorem which is obtained for
$(\bar{\nu}^{B}_n)_{n\in\mathbb{N}}$ and
$(\bar{\nu}^{C}_n)_{n\in\mathbb{N}}$ is strictly identical to
Proposition \ref{invariance} under the additional condition
\eqref{condsupplementaire} for Scheme (\ref{eC}). We recall that the
proof of Proposition~\ref{invariance} is based on two steps:
Propositions \ref{prop1} and \ref{prop2}. Proposition \ref{prop2}
is still valid without additional difficulties. However, the proof
of the analogous result to Proposition \ref{prop1} raises some new
difficulties. Denote by $A^{k,B}$ and $A^{k,C}$ the operators on
$\mathcal{C}_2^K(\mathbb{R}^d)$ with values in $\mathcal{C}_b(\mathbb
{R}^d,\mathbb{R})$
defined by
\begin{eqnarray*}
A^{k,B} f(x) &=& \langle\nabla f,b\rangle(x)+\tfrac{1}{2}
\operatorname{Tr}(\sigma^* D^2 f\sigma)(x)
+\int_{\{|y|\ge u_k\}}\tilde{H}^f(x,x,y)\pi(dy),
\\
A^{k,C} f(x) &=& A^{k,B} f(x)-\bigl(1-\alpha_k(\gamma_k)\bigr)
\int_{\{|y|\ge u_k\}}\tilde{H}^f(x,x,y)\pi(dy),
\end{eqnarray*}
where $\alpha_k(t)=\frac{1-e^{-\pi(|y|>u_k)t}}{\pi(|y|>u_k)t}$.
``$Af-A^{k,B}f$'' and ``$Af-A^{k,C}f$'' can be viewed as the
principal part of the weak error induced by the approximation in
Schemes (\ref{eB}) and (\ref{eC}) [$A^{k,B}$ is the infinitesimal generator of
$(X_t^k)$, where $(X_t^k)$ is solution to the SDE \eqref{edss}
driven by $(Z_t^k)$ instead of $(Z_t)$]. Thus, one may expect that
this error be negligible in the sense of our problem. This is the
aim of Lemma \ref{1810}.
\begin{lemme}\label{1810}
Assume \hyperlink{eH}{$\mathrm{(H^2_q)}$}. Let
$(x_k)_{k\in\mathbb{N}}$ be a sequence such that
%
\begin{equation}\label{1010}
\sup_{n\ge 1}\frac{1}{H_n}\sum_{k=1}^n\eta_k\|\kappa
(x_{k-1})\|^{2q}<\infty.
\end{equation}
Then, for every function $f\in\mathcal{C}_2^K(\mathbb{R}^d,\mathbb{R})$,
\[
\lim_{n\rightarrow+\infty}\frac{1}{H_n}\sum_{k=1}^n\eta_k
\bigl(Af(x_{k-1})-A^{k,B} f(x_{k-1}) \bigr)=0
\]
and if $\pi(D_n)\gamma_n\stackrel{n\rightarrow+\infty}{\longrightarrow}0$,
\[
\lim_{n\rightarrow+\infty}\frac{1}{H_n}\sum_{k=1}^n\eta_k
\bigl(Af(x_{k-1})-A^{k,C}f(x_{k-1}) \bigr)=0.
\]
\end{lemme}
\begin{pf}
Note that $A^{k,B}
f(x)-Af(x)=\int_{\{|y|<u_k\}}\tilde{H}^f(x,x,y)\pi(dy).$
When $q\ge1/2$, we deduce from Taylor's formula and the boundedness of
$\nabla f$ and $D^2f$
that there exists $C_q>0$ such that
\[
|\tilde{H}^f(x,x,y)|1_{\{|y|\le u_k\}}\le C_q \|\kappa(x)\|^{2q}
|y|^{2q} 1_{\{|y|\le u_k\}}.
\]
When $q\le1/2$, since $f$ is a $2q$-H\"older function,
\begin{eqnarray*}
|\tilde{H}^f(x,x,y)|1_{\{|y|\le u_k\}}
&\le & [f]_{2q}\|\kappa(x)\|^{2q} |y|^{2q} 1_{\{|y|\le u_k\}}
\\
&&{} +\sup_{x\in \operatorname{supp} f} |\nabla f(x)|\cdot \|\kappa(x)\|\,|y| 1_{\{|y|\le u_k\}}.
\end{eqnarray*}
By setting $v_{k,q}=\int_{\{|y|< u_k\}}|y|^{2q}\pi(dy)$, we have
\begin{eqnarray*}
|Af(x_{k-1})-A^{k,B} f(x_{k-1})|\le \cases{%
C  \bigl(v_{k,q}\|\kappa(x_{k-1})\|^{2q}+v_{k,1}\bigr), & \quad if $q\le1/2$,\cr
C v_{k,q}\|\kappa(x_{k-1})\|^{2q}, & \quad if $q\ge1/2$.}
\end{eqnarray*}
Since $v_{k,\alpha}\stackrel{k\rightarrow\infty}{\longrightarrow}0$ for every
$\alpha\ge q$ under assumption
\hyperlink{eH}{$\mathrm{(H^2_q)}$}, the first
result follows from \eqref{1010}. One deduces the second
inequality by checking that
\[
|A^{k,B}f(x)-A^{k,C}f(x)|\le C\pi(D_k)\gamma_k\bigl(1+\|\kappa(x)\|^{2q}\bigr).
\]\upqed
\end{pf}

Set
\begin{eqnarray*}
R_{3}^{B,k}(\gamma_k,\bar{X}^{B}_{k-1})
&=& \frac{\mathbb{E}\{f(\bar{X}^{B}_k)
- f(\bar{X}^{B}_{k-1})|\mathcal{F}^{B}_{k-1}\}}{\gamma_k}-A^{k,B}f(\bar{X}^{B}_{k-1}),
\\
R_{3}^{C,k}(\gamma_k,\bar{X}^{C}_{k-1})
&=& \frac{\mathbb{E}\{f(\bar{X}^{C}_k)-f(\bar{X}^{C}_{k-1})|\mathcal{F}^{C}_{k-1}\}}
{\gamma_k}-A^{k,C}f(\bar{X}^{C}_{k-1}).
\end{eqnarray*}
The rest of the proof then amounts to proving that
\[
\lim_{n\rightarrow\infty}
\frac{1}{H_n}\sum_{k=1}^{n}\eta_k
R_{3}^{B,k}(\gamma_k,\bar{X}^{B}_{k-1})=0
\]
and
\[
\lim_{n\rightarrow\infty}
\frac{1}{H_n}\sum_{k=1}^{n}\eta_k
R_{3}^{C,k}(\gamma_k,\bar{X}^{C}_{k-1})=0.
\]
We do not detail this proof based on the same approach as the
proof of Proposition~\ref{prop1} (see \cite{bib26} for more
details). However, we want to derive the main difficulties from
the proof. For Scheme (\ref{eB}), one deduces from the Ito formula that
\[
|R_{3}^{k,B}(\gamma,x)|\le\int\int_0^1 dv\int\pi(dy)
\mathbb{E} \bigl|\Delta\tilde{H}^f\bigl(S_{x,\gamma,u}+\kappa
(x)Z^{k}_{v\gamma},x,x,y\bigr)\bigr |\mathbb{P}_{U_1}(du).
\]
The right-hand term can be written
$R_{3,1}^{B,k}(\gamma,x)+R_{3,2}^{B,k}(\gamma,x)$, where
$R_{3,1}^{B,k}(\gamma,x)$ [resp. $R_{3,2}^{B,k}(\gamma,x)$] is
simply derived from $R_{3,1}(\gamma,x)$ [resp.
$R_{3,2}(\gamma,x)$], defined in the proof of Proposition
\ref{invariance}, by replacing $Z$ with $Z^{k}$. We focus on
$R_{3,1}^{B,k}$. One observes that the controls \eqref{56bis1}
and \eqref{56bis2} used for $R_{3,1}$ no longer work since the
jump component depends on $n$. An idea is to use the Skorokhod
representation theorem (see, e.g., \cite{rogers}) in order to
replace $(Z^k)$ by a uniformly controllable sequence.
\begin{lemme}\label{2'}
There exist a sequence of
càdlàg processes $(\tilde{Z}^n)$ and a càdlàg process $\tilde{Z}$
such that $\tilde{Z}^n\stackrel{\mathcal{L}}{=}Z^{n}$ for every $n\ge1$,
$\tilde{Z}\stackrel{\mathcal{L}}{=}Z$ and
$\tilde{Z}^n\rightarrow\tilde{Z}$ a.s. for the Skorokhod
topology. In particular,
%
\begin{eqnarray} \label{skorcontrol}
\sup_{n\in\mathbb{N}}\sup_{0\le s\le T}|\tilde{Z}_s^n|
&<& +\infty \qquad\forall  T>0\quad\mbox{and}\quad
\nonumber\\[-8pt]
\\[-8pt]
\limsup_{n\rightarrow +\infty,\gamma\rightarrow0}
\sup_{0\le s\le\gamma}|\tilde{Z}_s^n| &=& 0\qquad \mbox{a.s.}
\nonumber
\end{eqnarray}
\end{lemme}
\begin{pf}
$Z^n$ converges locally uniformly in $L^2$ toward
$Z$, hence, in distribution for the Skorokhod (Polish) topology.
Thanks to the Skorokhod representation theorem, there exists
$(\tilde{Z}^n)_{n\in\mathbb{N}}$ and $\tilde{Z}$ with $
\tilde{Z}^n\stackrel{\mathcal{L}}{=}Z^{n}$ and
$\tilde{Z}\stackrel{\mathcal{L}}{=}Z$ such that $\tilde{Z}^n$ tends a.s.
toward $\tilde{Z}$ for
Skorokhod topology. The assertion \eqref{skorcontrol} easily
follows from the continuity of $\alpha\mapsto\|\alpha\|_{\sup}$ and
$\alpha\mapsto\alpha(0)$ for the Skorokhod topology.
\end{pf}

Since $R_{3,1}^{B,k}$ only depends on the law of $Z^n$, one can
replace $Z^{n}$ with $\tilde{Z}^n$. Then, we use
\eqref{skorcontrol} as an alternative to the local boundedness and
the continuity at $t=0$ of $(Z_t)$ needed in \eqref{56bis1} and
\eqref{56bis2} respectively. A result analogous to \eqref{956bis}
follows. The idea is the same for $R_{3,2}^{B,k}$.

Finally, for Scheme (\ref{eC}), the result essentially follows
from the following remark:
\[
\sup_{0<s\le t}|Z_{s\wedge T^n}^n|\le\sup_{0<s\le t}|Z_s^n|.
\]
This means that the remainders in Scheme (C) are easier to control
than those of Scheme (B). For more details, we refer to
\cite{bib26}.

\section{A theoretical application}\label{section7}

The ``classical'' {a.s.} CLT due to
Brosamler \cite{bib6} and Schatte \cite{bib8} is the following
result. Let $(U_n)_{n\in\mathbb{N}}$ be a sequence of i.i.d.
random variables with values in $\mathbb{R}^d$ such that $\mathbb{E} U_1=0$
and $\Sigma_{U_1}=I_d$. Then,
\[
\mathbb{P}\mbox{-a.s}\qquad
\frac{1}{\ln n}\sum_{k=1}^n\frac{1}{k}\delta_{{1}/{\sqrt{k}} \sum_{i=1}^k
U_i} \Longrightarrow \mathcal{N}(0,I_d).
\]
This result is obviously connected with the central limit theorem
which expresses the fact that every square-integrable centered
random variable is in the domain of normal attraction of the
normal law. When the square-integrability no longer holds, Berkes,
Horvath and Khoshnevisan \cite{bib9} obtained an extension of
this result connected with the nonsquare-integrable attractive
laws which are stable laws [with index $\alpha\in(0,2)$].
We are going to show that we can deduce this extension from Theorem
\ref{principal}.

Let $c$ be a positive number and denote by $(Z^{\alpha,c}_t)_{t\ge
0}$ a symmetrical one-dimensional $\alpha$-stable process such
that the characteristic function $\phi$ of $Z^{\alpha,c}_1$
satisfies $\phi(u)=e^{-\rho|u|^\alpha}$, where
$\rho=2c\int_0^{+\infty}y^{-\alpha}\sin y \,dy$. Consider a
sequence $(V_n)_{n\in\mathbb{N}}$ of symmetrical i.i.d. random
variables such that, for $x>0$,
%
\begin{eqnarray} \label{01710}
\mathbb{P}(V_1\ge x) &=& \frac{c}{x^\alpha}+\delta(x)(x^{-\alpha}(\ln
x)^{-\gamma})
\nonumber\\[-8pt]
\\[-8pt]
\eqntext{\mbox{with } \gamma>0 \mbox{ and }
\delta(x)\stackrel{x\rightarrow+\infty}{\longrightarrow}0.}
\end{eqnarray}
By a
result of Gnedenko and Kolmogorov (see \cite{bib10}), we know that
\[
\frac{V_1+\cdots+V_n}{n^{1/\alpha}}\Longrightarrow Z^{\alpha,c}_1.
\]
Then, the following {a.s.} CLT holds:
\begin{theorem}\label{TCL}
Let $(\eta_k)_{k\in\mathbb{N}}$ be a nonincreasing sequence
with infinite sum such that $(k\eta_k)_{k\in\mathbb{N}}$ is
nonincreasing and set
$\nu=\mathcal{L}(Z^{\alpha,c}_1)$. Then, if $\gamma>\frac{1}{\alpha}$,
a.s.,
\[
\frac{1}{H_n}
\sum_{k=1}^n\eta_k\delta_{{(V_1+\cdots+V_k)}/{k^{1/\alpha}}}
\stackrel{(\mathbb{R})}{\Longrightarrow}\nu.
\]
In particular,
\[
\frac{1}{\ln n}\sum_{k=1}^n\frac{1}{k}
\delta_{{(V_1+\cdots+V_k)}/{k^{1/\alpha}}}
\stackrel{(\mathbb{R})}{\Longrightarrow}\nu\qquad \mbox{a.s.}
\]
\end{theorem}

In order to prove this theorem, we first need an almost
sure invariance principle due to Stout (see \cite{bib7} or
\cite{bib9}).
\begin{prop}\label{stout}
Let $(V_n)_{n\ge1}$ and $(\zeta_n)_{n\ge1}$ be
sequences of i.i.d. random variables such that
$\zeta_1\stackrel{\mathcal{L}}{=}Z^{\alpha,c}_1$ and $V_1$ is defined as
above. Then, if $\gamma>\frac{1}{\alpha}$, there exists a
probability space $(\hat{\Omega},\hat\mathcal{F},\hat{\mathbb{P}})$
and sequences of i.i.d. random variables $(\hat{V}_n)_{n\ge1}$ and
$(\hat{\zeta}_n)_{n\ge1}$ such that $\hat{V_1}\stackrel{\mathcal{L}}{=}V_1$,
$\hat{\zeta}_1\stackrel{\mathcal{L}}{=}\zeta_1$ and
%
\begin{equation}\label{invariancerevprinciple}
\sum_{i=1}^n\hat{\zeta}_i -\sum_{i=1}^n
\hat{V_i}\stackrel{n\rightarrow+\infty}{=}
o(n^{1/\alpha}(\ln n)^{-\rho})\qquad
\mbox{a.s. }\forall\rho\in\biggl(0,\gamma-\frac{1}{\alpha}\biggr).
\end{equation}
\end{prop}
\begin{pf*}{Proof of Theorem \ref{TCL}}
First, we assume that $V_1=\zeta_1\stackrel{\mathcal{L}}{=}Z^{\alpha,c}_1$. Set
\[
S_n=\frac{\zeta_1+\cdots+\zeta_{n+1}}{(n+1)^{1/\alpha}}
\qquad\forall n\ge0 .
\]
One easily checks that
$S_{n+1}=S_n-\frac{1}{\alpha}\gamma_{n+1}
S_n+\gamma_{n+1}^\frac{1}{\alpha}\zeta_{n+2}+R_{n+1}$
with $\gamma_n=\frac{1}{n+1}$ and $R_{n+1}=O(\gamma_{n+1}^2|S_n|)$. The
idea of the proof is to compare $(S_n)_{n\ge0}$ with the exact
Euler scheme with initial value $\zeta_1$ and step sequence
$(\gamma_n)$ associated with the SDE $\mathbf{(E_{\alpha,c})}$
defined by $dX_t=-\frac{1}{\alpha}X_{t^-}\,dt+dZ^{\alpha,c}_t$.
Since $(Z_t^{\alpha,c})_{t\ge0}$ is a self-similar process with
index ${1}/{\alpha}$ (see, e.g., \cite{bib11}), its Euler scheme
can be written
\[
\bar{X}_0=\zeta_1 \quad\mbox{and} \quad
\bar{X}_{n+1}=\bar{X}_n-\frac{1}{\alpha}\gamma_{n+1}
\bar{X}_n+\gamma_{n+1}^\frac{1}{\alpha}\zeta_{n+2} .
\]
As an
Ornstein--Uhlenbeck process driven by a symmetric stable law,
$(X_t)$ admits a unique invariant
measure $\nu$ and $\nu=\mathcal{L}(Z^{\alpha,c}_1)$ (see \cite{bib11}, page~188).
Since $\kappa$ is bounded, assumptions of Theorem \ref{principal}
are clearly fulfilled with $V(x)=1+x^2$, $a=1$ and for any
$p\in(0,\alpha/2)$ and $q\in(\alpha/2,1)$.
(In the rest of
the paper the initial value of the Euler scheme is supposed to be
constant, but it is obvious that Theorem \ref{principal} is still
true when $\bar{X}_0$ is a random variable satisfying
$\mathbb{E}\{|\bar{X}_0|^{2p}\}<+\infty$.)
Hence, it follows from Theorem \ref{principal} that
%
\begin{equation}\label{1700}
\frac{1}{H_n}\sum_{k=1}^n\eta_k\delta_{\bar{X}_{k-1}}
\stackrel{n\rightarrow+\infty}{\Longrightarrow}\nu\qquad\mbox{a.s.}
\end{equation}
Then, by using that $|f(S_k)-f(\bar{X}_k)|\le C(|S_k-\bar{X}_k|\wedge1)$
for every Lipschitz bounded function
$f$, one easily checks that Theorem \ref{TCL} holds with
$V_1=\zeta_1$ if
%
\begin{equation}\label{deltantendrev}
\Delta_n:=S_n-\bar{X}_n\stackrel{n\rightarrow+\infty}{\longrightarrow}0
\qquad \mbox{a.s.}
\end{equation}
Let us show \eqref{deltantendrev}. One first checks that
\[
\Delta_0=0 \quad\mbox{and}\quad
\Delta_{n}= \biggl(1-\frac{1}{\alpha(n+1)} \biggr)
\Delta_{n-1}+R_{n}\qquad \forall n\ge1.
\]
Setting $k_0=\inf\{k\ge0,k-1/\alpha>0\}$, we
deduce that, for every $n\ge k_0+1$,

\[
\Delta_n=\frac{\Delta_{k_0}}{c_n}+\frac{1}{c_n}\sum_{k=k_0+1}^n
c_kR_k\qquad\mbox{with }
c_n=\prod_{k=k_0+1}^n\biggl(1-\frac{1}{\alpha(k+1)}\biggr)^{-1}.
\]
One
observes that
\[
c_n=\exp \Biggl(-\sum_{k=k_0+2}^{n+1}
\ln\biggl(1-\frac{1}{\alpha k}\biggr) \Biggr)
= \exp \Biggl(\frac{1}{\alpha}\sum_{k=k_0+2}^{n+1}
\frac{1}k+O\biggl(\frac{1}{k^2}\biggr) \Biggr)
\stackrel{n\rightarrow+\infty}{\sim} C'n^{1/\alpha}.
\]
Then, ${\Delta_{k_0}}/{c_n}\rightarrow0$
a.s. Hence, \eqref{deltantendrev} holds if we check that
${1}/{c_n}\sum_{k=k_0+1}^n c_kR_k\rightarrow0$ a.s. First, if
$\alpha>1$, as $\zeta_1$ is integrable and
$R_k=O(S_{k-1}/(k+1)^2)$, we have
\begin{eqnarray*}
\sum_{k\ge1} \mathbb{E}\{|R_k|\}
&\le & C\sum_{k\ge1} \frac{\mathbb{E}\{|S_{k-1}|\}}{(k+1)^2}
\le C \sum_{k\ge1}\frac{k\mathbb{E}\{|\zeta_1|\}}{(k+1)^{1/(\alpha+2)}}
\\
&\le & C\sum_{k\ge1}\frac{1}{(k+1)^{1+1/\alpha}}<+\infty.
\end{eqnarray*}
We deduce that $\sum_{k\ge1} |R_k|<+\infty$
a.s. Since $c_n\stackrel{n\rightarrow+\infty}{\longrightarrow}+\infty$, we
derive from Kronecker's lemma that
\[
\frac{1}{c_n}\sum_{k=k_0+1}^n c_k R_k
\stackrel{n\rightarrow +\infty}{\longrightarrow}0
\quad\Longrightarrow\quad
\Delta_n\stackrel{n\rightarrow+\infty}{\longrightarrow}0\qquad
\mbox{a.s. if } \alpha>1.
\]
Second, if $\alpha\le1$, $\zeta_1$ has a moment of order $\theta$
for every $\theta<\alpha$. It follows from inequality
\eqref{eqelem2} that
\[
\mathbb{E}\{|R_k|^\theta\}\le
C\frac{k}{(k+1)^{2\theta+\theta/\alpha}}\mathbb{E}\{|\zeta
_1|^\theta\}\le\frac{C}{(k+1)^{\theta(2+1/\alpha)-1}}.
\]
Therefore, if $\theta$ satisfies $\theta(2+\frac1 \alpha)-1>1$,
that is, if $\frac{2\alpha}{2+\alpha}<\theta<\alpha$,
we have $\sum_{k\ge1}|R_k|^\theta<+\infty$ {a.s.} Hence, by
inequality \eqref{eqelem2} and Kronecker's lemma, it follows that
\[
\Biggl |\frac{1}{c_n}\sum_{k=k_0+1}^n c_k R_k \Biggr|^\theta
\le\frac{1}{c_n^\theta}\sum_{k=k_0+1}^n
c_k^\theta|R_k|^\theta\stackrel{n\rightarrow+\infty}{\longrightarrow}0
\qquad \mbox{a.s.}
\]
and the theorem is proved when $V_1=\zeta_1$. Now, consider
a sequence $(V_n)_{n\ge0}$ of i.i.d. symmetric random variables
satisfying \eqref{01710}. Since Theorem \ref{TCL} is true for
$(\zeta_n)_{n\ge1}$, it is also true for every sequence
$(\hat{\zeta}_n)_{n\ge1}$ of i.i.d. random variables satisfying
$\hat{\zeta}_1\stackrel{\mathcal{L}}{=}\zeta_1$. By taking
$(\hat{\zeta}_n)_{n\ge1}$ such that Proposition \ref{stout} holds,
we derive from \eqref{invariancerevprinciple} that there exists a
sequence of i.i.d. random variables $(\hat{V}_n)_{n\ge1}$ such\vadjust{\goodbreak}
that $V_1\stackrel{\mathcal{L}}{=}\hat{V_1}$ and
%
\begin{equation}\label{1710}
\frac{1}{H_n}\sum_{k=1}^n\eta_k
\delta_{{(\hat{V}_1+\cdots+\hat{V}_k)}/{k^{1/\alpha}}}
\stackrel{n\rightarrow+\infty}{\Longrightarrow}\nu
\qquad\mbox{a.s.}
\end{equation}
As $(V_n)_{n\ge1}$ and $(\hat{V}_n)_{n\ge1}$ are sequences of
i.i.d. random variables such that \mbox{$V_1\stackrel{\mathcal{L}}{=}\hat{V}_1$},
\eqref{1710} is also true for $(V_n)_{n\ge1}$.
\end{pf*}

\section{Simulations}\label{simulations}

\begin{figure}[b]

\includegraphics{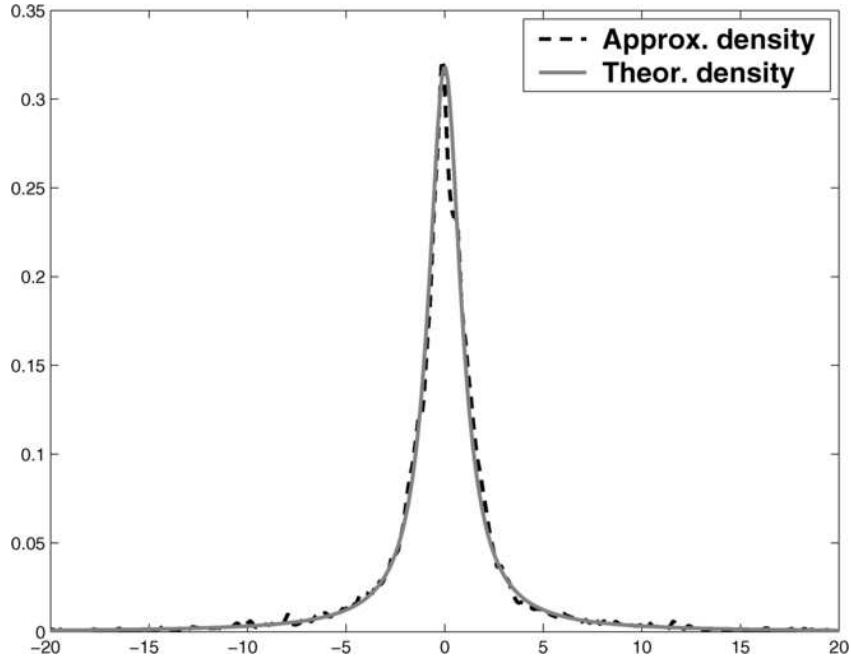}

\caption{Scheme \textup{(\protect\ref{eA})}, $t=12.5$.\label{f1}}
\end{figure}
\begin{figure}

\includegraphics{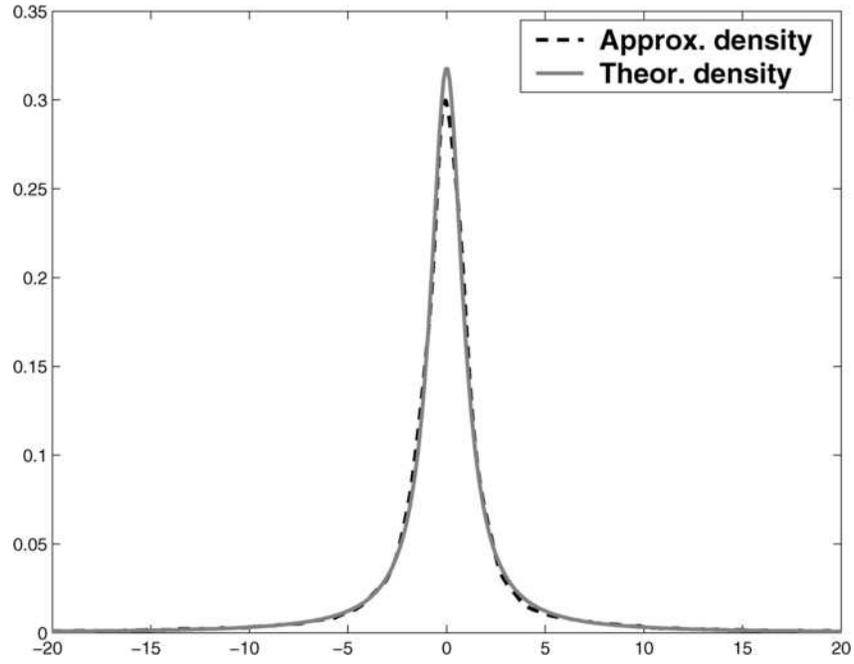}

\caption{Scheme \textup{(\protect\ref{eB})}, $t=16.6$.\label{f2}}
\end{figure}
\begin{figure}

\includegraphics{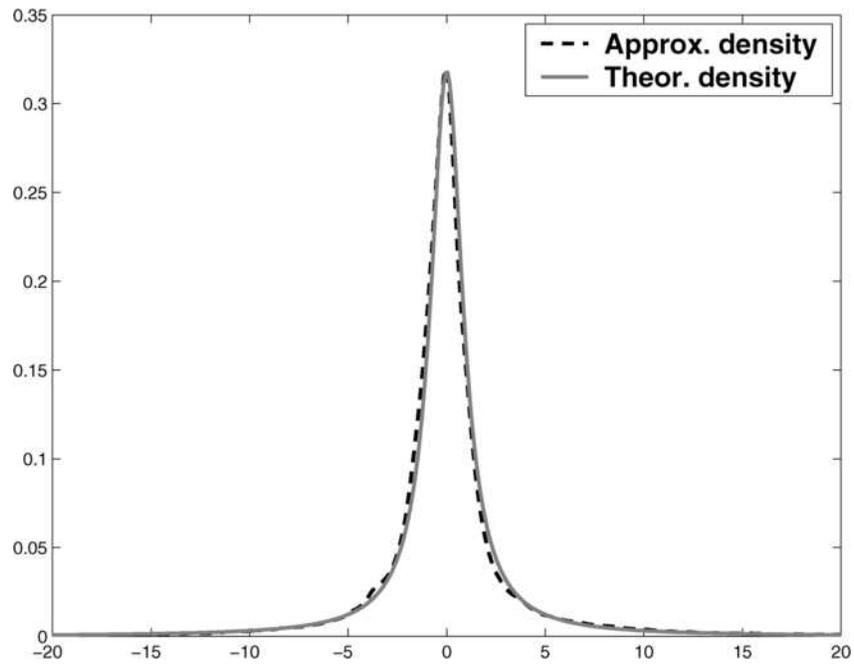}

\caption{Scheme \textup{(\protect\ref{eC})}, $t=16.4$.\label{f3}}
\end{figure}
\begin{figure}

\includegraphics{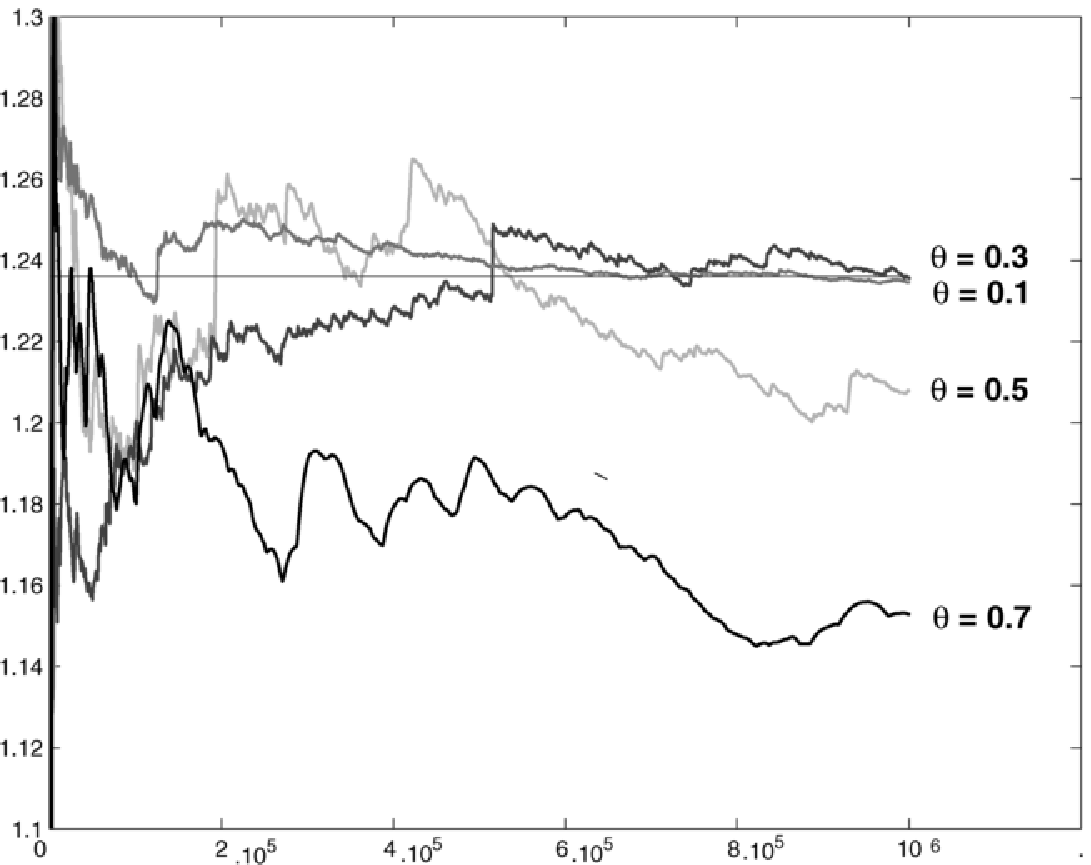}

\caption{Scheme \textup{(\protect\ref{eA})}.\label{f4}}
\end{figure}
\begin{figure}

\includegraphics{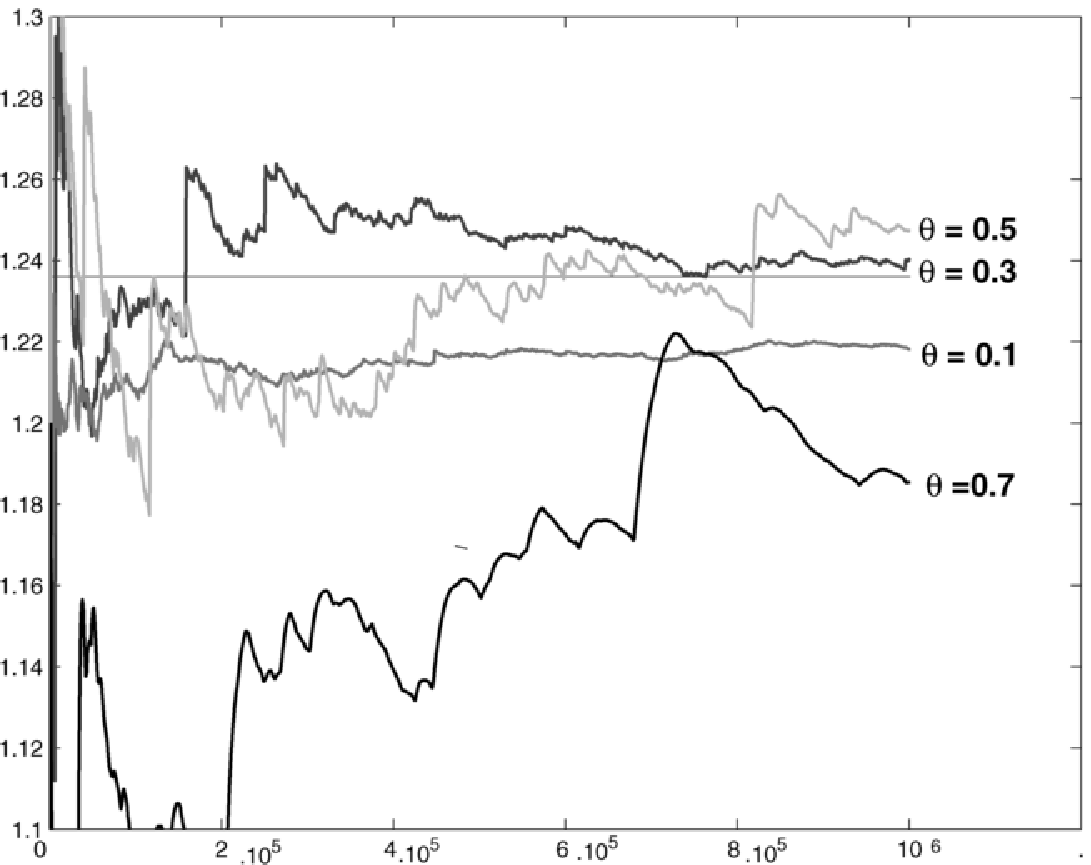}

\caption{Scheme \textup{(\protect\ref{eB})}.\label{f5}}
\end{figure}
\begin{figure}

\includegraphics{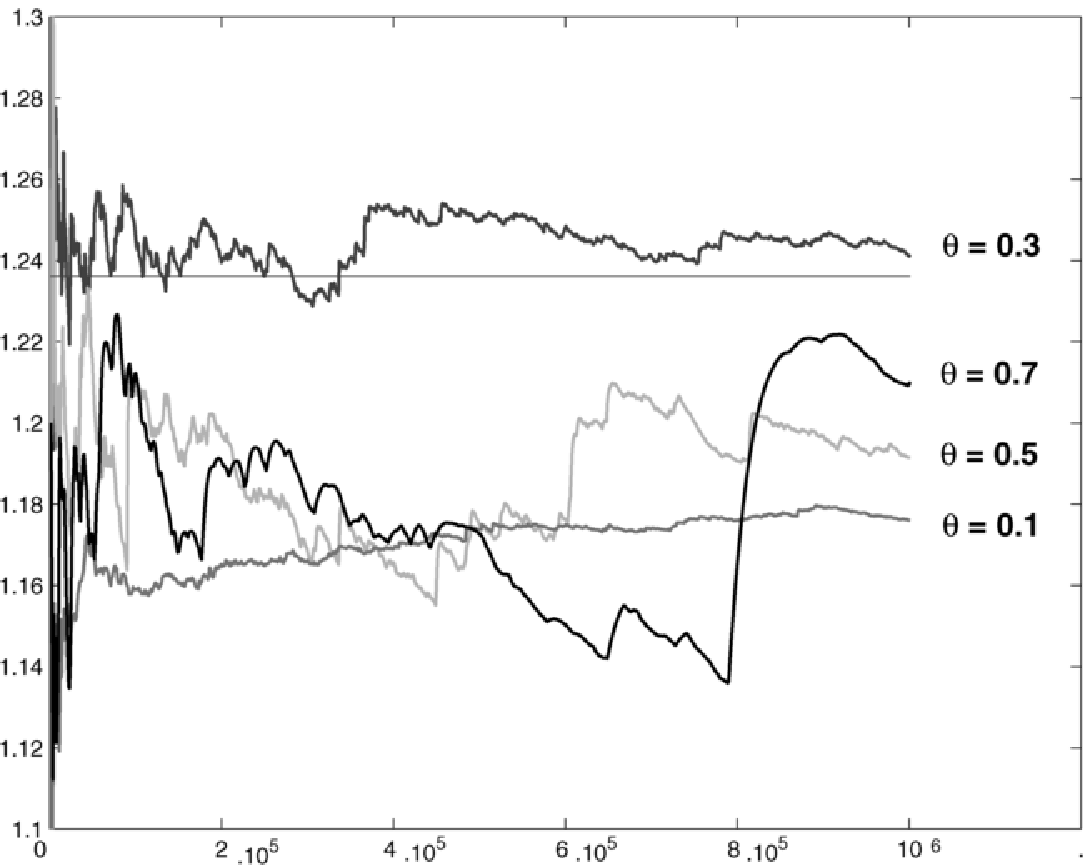}

\caption{Scheme \textup{(\protect\ref{eC})}.\label{f6}}
\end{figure}
\begin{figure}

\includegraphics{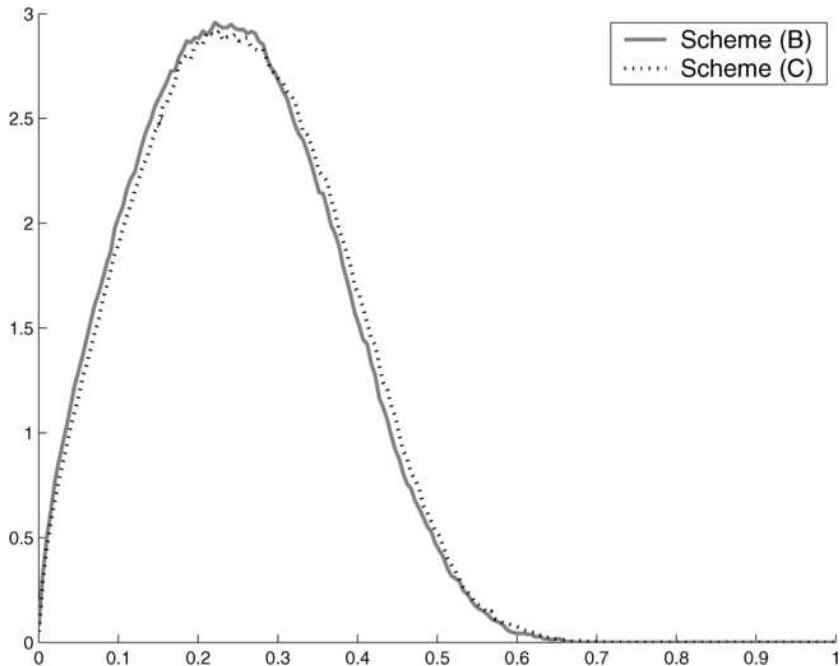}

\caption{Approximated density, $N=10^6$.\label{f7}}
\end{figure}

%
\begin{Exemple}
Denote by $(Z_t)_{t\ge0}$ a Cauchy process
with parameter 1 [with L\'evy measure defined by
$\pi(dy)=1/y^2\,dy$] and consider the Ornstein--Uhlenbeck process
solution to $dX_t=-X_{t^-}\,dt+dZ_t$ corresponding to
$(E_{1,1})$ defined in the previous subsection. The unique
invariant measure of $(X_t)_{t\ge0}$ is the Cauchy law (see
\cite{bib11}, page~188) and the assumptions of Theorem
\ref{principal} are fulfilled with $V(x)=1+x^2$, $a=1$ and every
$p\in(0,1/2)$ and $q\in(1/2,1)$. Therefore,
\[
\bar{\nu}_n(f), \bar{\nu}^B_n(f), \bar{\nu}^C_n(f)
\stackrel{n\rightarrow+\infty}{\longrightarrow}
\int\frac{f(x)}{\pi(1+x^2)}\,dx\qquad \mbox{a.s.}
\]
for every $f$ satisfying $f=O(|x|^{{1}/{2}-\varepsilon})$
with $\varepsilon>0$. In Figures \ref{f1}, \ref{f2} and \ref{f3},\vadjust{\goodbreak}
one compares the
theoretical density of the invariant measure with the density
obtained by convolution of each of the empirical measures by a
Gaussian kernel for $N=5.10^4$. We choose
$\eta_n=\gamma_n=1/\sqrt{n}$, $u_n=\sqrt{\gamma_n}$ [so that
$\pi(D_n)\gamma_n\rightarrow0$] and $t$ indicates the CPU time.
In order to have a more precise idea of the differences between
the three Euler schemes, we simulate and represent on Figures \ref{f4},
\ref{f5} and \ref{f6} the sequence $(\bar{\nu}_n(f))$ with $f(x)=|x|^{0.4}$, for
several choices of polynomial steps. We set
$\gamma_n=\eta_n=1/n^\theta$ and $u_n=\gamma_n$ (resp.
$u_n=\sqrt{\gamma}_n$) for Scheme (\ref{eB}) [resp. for Scheme~(\ref{eC})].
We observe that, among the tested steps, the best rate seems to be
obtained for $\theta=0.3$. Notably, in Schemes (\ref{eB}) and (\ref{eC}), we
see that, on the one hand, if the step decreases too slowly (e.g.,
$\theta=0.7$), so is the stabilization and, on the other
hand, if the steps decreases too fast (e.g.,
when $\theta =0.1$), there are not sufficient variations to correct the error.
\end{Exemple}
%
%
\begin{Remarque}
In \cite{panloup} we study the rate of convergence of these procedures in terms of
steps, weights and truncation thresholds. This enlightens these
first numerical illustrations.
\end{Remarque}
\begin{Exemple}
Now we deal with the following SDE:
\[
dX_t=(1-X_{t^-})\,dt-X_{t^-}\,dZ_t,
\]
where $(Z_t)_{t\ge0}$ is a drift-free subordinator with L\'evy
measure $\pi$ defined by
\[
\pi(dy)=\frac{f_{3/2,1/2}(y)}{y^2}\,dy,
\]
where $f_{a,b}$ is the density function of the $\beta(
a,b)$-distribution. This SDE models the dust generated by a
particular EFC process (see \hyperref[s1]{Introduction}) whose sudden
dislocations do not create dust, having parameters (according to
the notation of~\cite{bib17}):
\[
c_k=0, \qquad
c_e=1, \qquad
\nu_{\mathit{coag}}(dy)= f_{3/2,1/2}(y)\,dy.
\]
One checks that $\mathrm{(S_{1,1,{1}/{2}})}$ is satisfied with
$V(x)=1+x^2$. However, we do not have $\kappa(x)=o(|x|)$, but since
$\rm{supp}(\pi)$ is restrained to $[0,1]$ without singularities in
$0$ and $1$, we are able to show that assumption $\kappa(x)=o(x)$
is no longer necessary in this case. In Figure~\ref{f7} we represent the
approximation of the invariant measure obtained for Schemes (\ref{eB})
and (\ref{eC}) [we are not able to simulate Scheme (\ref{eA}) in that case].
\end{Exemple}

\section*{Acknowledgments}

Thanks to Gilles Pag\`es for extensive discussions and suggestions.

\printaddresses

\end{document}